\documentclass[11pt]{amsart}
\usepackage{graphicx}
\usepackage{amssymb}
\usepackage{epstopdf}

\newtheorem{lemma}{Lemma}[section] 
\newtheorem{propos}[lemma]{Proposition}
\newtheorem{example}[lemma]{Example}
\newtheorem{theorem}[lemma]{Theorem}
\newtheorem{cor}[lemma]{Corollary}
\newtheorem{corol}[lemma]{Corollary}
\newtheorem{defin}[lemma]{Definition}
\newtheorem{remark}[lemma]{Remark}

\newcommand{\CC}{{\hbox{{$\mathcal C$}}}}
\newcommand{\CD}{{\hbox{{$\mathcal D$}}}}

\newcommand{\CF}{\hbox{{$\mathcal F$}}}
\newcommand{\CG}{\hbox{{$\mathcal G$}}}
\newcommand{\CK}{\hbox{{$\mathcal K$}}}

\newcommand{\CH}{\hbox{{$\mathcal H$}}}

\newcommand{\CM}{\hbox{{$\mathcal M$}}}

\newcommand{\CR}{\hbox{{$\mathcal R$}}}

\newcommand{\C}{\mathbb{C}}
\newcommand{\R}{\mathbb{R}}

\newcommand{\Z}{\mathbb{Z}}

\newcommand{\del}{\partial}

\newcommand{\ev}{\mathrm{ev}}
\newcommand{\coev}{\mathrm{coev}}

\newcommand{\eps}{{\epsilon}}
\newcommand{\tens}{\mathop{\otimes}}
\newcommand{\la}{{\triangleright}}
\newcommand{\ra}{{\triangleleft}}

\newcommand{\id}{\mathrm{id}}

\newcommand{\<}{\langle}
\renewcommand{\>}{\rangle}

\newcommand{\varf}{{f^1}}
\newcommand{\bb}{{\mathrm{bb}}}
\newcommand{\YD}{{\buildrel {\rm x}\over \CM}}

\renewcommand{\o}{{}_{\scriptscriptstyle(1)}}
\renewcommand{\t}{{}_{\scriptscriptstyle(2)}}

\newcommand{\uo}{{}^{\scriptscriptstyle(1)}}
\newcommand{\ut}{{}^{\scriptscriptstyle(2)}}

\title{Bar categories and star operations}
\author{E.J.\,Beggs and S.\,Majid}
\address{EJB: Department of Mathematics, Swansea University\\ Singleton Park, Swansea SA2 8PP\\
SM: School of Mathematical Sciences, Queen Mary University of London\\ 327 Mile End Rd, London E1 4NS}
\dedicatory{Dedicated to Fred Van Oystaeyen, on the occasion of his sixtieth birthday}

\begin{document}
\maketitle

\begin{abstract}
We introduce the notion of `bar category' by which we mean a monoidal category equipped with additional structure formalising the notion of complex conjugation. Examples of our theory include bimodules over a $*$-algebra, modules over a conventional $*$-Hopf 
algebra and modules over a more general object which call a `quasi-$*$-Hopf algebra' and for which examples include the standard quantum groups $u_q(\mathfrak{g})$ at $q$ a root of unity (these are well-known not to be a usual $*$-Hopf algebra). We also provide examples of strictly quasiassociative bar categories, including modules over `$*$-quasiHopf algebras' and a construction based on finite subgroups $H\subset G$ of a finite group. Inside a bar category one has natural notions of `$\star$-algebra' and  `unitary object' therefore extending these concepts to a variety of new situations. We study braidings and duals in bar categories and $\star$-braided groups (Hopf algebras) {\em in} braided-bar categories. Examples include the transmutation $B(H)$ of a quasitriangular $*$-Hopf algebra and the quantum plane $\C_q^2$ at certain roots of unity $q$ in the bar category of  $\widetilde{u_q(su_2)}$-modules. We use our methods to provide a natural quasi-associative $C^*$-algebra structure on the octonions ${\mathbb O}$ and on a coset example. In the appendix we  extend the Tannaka-Krein reconstruction theory to bar categories in relation to $*$-Hopf algebras. 
\end{abstract}

\section{Introduction}

One of the key ingredients in an operator algebra approach to noncommmutative geometry is that of complex conjugation and its extension to adjoints, $*$-involutions on algebras and other forms of structure. Here commutative unital $C^*$-algebras correspond by the  theorem of Gelfand and Naimark to compact Hausdorf spaces, and many geometric ideas extend from this point of view to noncommutative $C^*$-algebras. Even if one is not concerned with $C^*$-completions, one may be concerned with the properties of $*$-algebras at an algebraic level and associated concepts. Thus, for example, there is an established abstract notion of $*$-Hopf algebras (also known as
Hopf $*$-algebras)
 over $\C$ which in some cases can be completed to a compact quantum group \cite{woron87}. When one introduces constructions in $*$-noncommutative geometry one needs $*$ to be respected and this is usually possible in the context of Hilbert spaces and adjoints. However, even in this case the notion of bar category that we introduce is useful as  a book keeping device to make  
explicit the idea of conjugate objects and antilinear maps in various contexts relating to bimodules and conjugate representations.

The framework of bar category moreover allows such concepts to be extended from the familiar case 
of bimodules over  $*$-algebras or (co)modules over a $*$-Hopf algebra to other situations were the corresponding notions of conjugation would otherwise not be clear. The most important of these we believe is to standard quantum groups $\C_q[G]$ or 
$u_q(\mathfrak{g})$ at $|q|=1$. These are not standard $*$-Hopf algebras but we show that their modules still form a bar catgeory. We also introduce a theory of $*$-structures for Drinfeld's quasiHopf algebras \cite{driquasi} such that the category of modules is a bar category, and a construction for obtaining these by twisting a $*$-Hopf algebra. Another application is to quasi-associative algebras \cite{albmajoct} where without a categorical framework it would not be clear how to define concepts such as conjugation of bimodules or `skew-linear inner products' to define Hilbert spaces. Guided by our examples we actually introduce two notions in the categorical setting, namely a general bar category and a stronger notion which we call a `strong bar category' in which the bar operation is more strongly involutive.

The related idea of star operation $\star$ for objects {\em in} a bar category is introduced. Quantum planes at $|q|=1$, the octonions and our coset quasi-associative algebra  example are all shown to be a star algebras in a very natural way, though the star operations in the first and last case do not square to the identity. The quantum plane example in fact turns out to be a braided $\star$-Hopf algebra (in a braided bar category). The coset example is closely related to
ideas in \cite{BHMnonass}, which is concerned with applications to mathematical physics. Other motivation comes from earlier attempts at such notions in the specific context of braided linear spaces \cite{majquastpoin,majstarbr}. 

One of the original motivations for this work was
noncommutative Riemannian geometry, where we require $*$-structures in noncommutative differential forms \cite{connesbook}, projective modules etc. and this will appear in a sequel.
It is also related to the `type A' and `type B' morphisms mentioned in 
\cite{shombeggs3}.

 An outline of the paper is as follows. The basic definition is given in Section~2, with elementary examples such as the bar category {\bf Vect} of vector spaces with conjguation and the bimodule, module and comodule categories. Section~3 contains new examples of bar categories beyond the standard settings. Section~4 covers braiding in bar categories. Section~5 covers $\star$-algebras and $\star$-Hopf algebras in (braided) bar categories.  Section 6 concludes with the natural notion of dual objects in a bar category. A main result here is that in a bar category a left-dual implies a right-dual. 
 
 In the appendix we show that a $*$-Hopf algebra can be reconstructed from
its (bar) category of modules, showing that in some sense, this
is the `correct' notion of bar category \cite{majquasitann}.

The reader should note that this material is very different,
both in definition and application, from the idea of
*-autonomous category given in  \cite{barr}. For one thing, there $*$ is a  functor
rather than a natural transformation. The $*$  functor in \cite{barr} is contravariant,
rather than the covariant bar functor described here.
However, an idea of formulating bar as a functor did appear in the work of Baez on 2-Hilbert spaces
\cite{baez2hilbert}. 

As regards notation, we uniformly use $*$ as a superscript to denote an 
antilinear operation, and $\star$, which is never a superscript, to denote a certain
 morphism in a category. The letter $\Upsilon$ is a capital upsilon. The notation used for
various sorts of $*$ algebras which are also (quasi)Hopf algebras is hardly 
optimal, but is mostly fixed by past usage in individual cases.  The reader should
 not confuse `quasi-$*$-Hopf algebras' and `$*$-quasiHopf algebras', both of which we introduce in the paper.

The paper was written during the visit July-December 2006 of the authors to the  Isaac Newton Institute. We thank the institute for their support.

\section{Bar categories}

Here we give the abstract definition that we propose, and the elementary `model' examples that justify the definitions. We recall that a monoidal category $(\CC,\tens,\Phi,1_\CC,l,r)$ means a category, with a functor $\tens:\CC\times\CC\to \CC$, a natural equivalence $\Phi:((\ \tens\ )\tens\ )\to (\  \tens(\ \tens\ ))$ subject to Mac Lane's pentagon coherence identity and an identity object $1_\CC$ and associated natural isomorpisms $l:\id\to \id\tens 1_{\CC}$ and $r:\id\to 1_\CC\tens\id$ compatible with $\Phi$. We refer to \cite{Mac} for details. We recall that a monoidal functor $F:\CC\to \CD$ between monoidal categories means a functor  $F$ and a natural equivalence $f: F(\ )\tens F(\ )\to F\circ\tens$ with
\[ f_{X,Y\tens Z}\circ(\id\tens f_{Y,Z})\circ\Phi_{F(X),F(Y),F(Z)}=F(\Phi_{X,Y,Z})\circ f_{X\tens Y,Z}\circ (f_{X,Y}\tens\id)\]
\[ F(1_\CC)=1_\CD,\quad f_{1_{\CC},X}\circ l_{F(X)}=F(l_X),\quad f_{X,1_{\CC}}\circ r_{F(X)}=F(r_{X})\]
for all objects $X,Y,Z$ of $\CC$. Note that in a monoidal category the unit object is unique up to isomorphism and consequently is it conventional without loss of generality to assume strict equality in the last line above. If one wishes to be more precise it would be better to say that there is an isomorphism $f^1$ with
\[ \varf:1_\CD\to F(1_\CC),\quad   f_{1_{\CC},X}\circ (\varf\tens\id)\circ l_{F(X)}=F(l_X),\quad f_{X,1_{\CC}}\circ (\id\tens\varf)\circ r_{F(X)}=F(r_{X})\]
for all $X$. Here $f^1$ if it exists is determined uniquely from $(F,f)$ and we say that $F$ is monoidal with {\em induced unit isomorphism} $\varf$.

A natural equivalence means that the natural transformation in question is invertible. We recall that a natural transformation means that the collection of morphisms is functorial with respect to morphisms in the domain of the relevant functors.

\subsection{The definition of a bar category}

We denote typical objects of a bar category by $X,Y,Z$ etc. We denote by $\mathrm{flip}$ the
functor from ${\CC}\times {\CC}$ to ${\CC}\times {\CC}$ sending
$(X,Y)$ to $(Y,X)$.

\begin{defin} \label{barcatdef}
A bar category is a tensor category $({\CC},\tens,1_{\CC}, l,r,\Phi)$
together with the following data:

\noindent
1)\quad A functor $\mathrm{bar}:{\CC}\to {\CC}$ (written 
as $X\mapsto \overline{X}$).

\noindent
2)\quad A natural equivalence $\mathrm{bb}$ between the 
identity and the $\mathrm{bar}\circ\mathrm{bar}$ functors on ${\CC}$.

\noindent
3)\quad An invertible morphism $\star:1_{\CC} \to\overline{1_{\CC}}$.

\noindent
4)\quad A natural equivalence $\Upsilon$ (capital upsilon) between
$\mathrm{bar}\circ\tens$ and $\tens\circ(\mathrm{bar}\times \mathrm{bar})
\circ\mathrm{flip}$
from ${\CC}\times {\CC}$ to ${\CC}$.
\smallskip These are required to obey the rules:

\noindent
a)\quad The following compositions are $r_{\bar X}$ and $l_{\bar X}$ respectively:
\begin{eqnarray*}
&&\overline{X} \stackrel{\overline{l_X}}\longrightarrow \overline{X\tens 1_{\CC}}
\stackrel{\Upsilon_{X,1_{\CC}}}\longrightarrow \overline{ 1_{\CC}} \tens
\overline{X} \stackrel{\star^{-1}\tens\id}\longrightarrow 1_{\CC} \tens
\overline{X}\ ,
\cr
&&\overline{X} \stackrel{\overline{r_X}}\longrightarrow \overline{1_{\CC}\tens X}
\stackrel{\Upsilon_{1_{\CC},X}}\longrightarrow 
\overline{X}\tens\overline{ 1_{\CC}} 
 \stackrel{\id\tens \star^{-1}}\longrightarrow \overline{X} \tens 1_{\CC} \ .
\end{eqnarray*}

\noindent
b)\quad The following condition is satisfied by $\Upsilon$
and the associator $\Phi$:
\begin{eqnarray*}
\Phi_{\overline{Z},\overline{Y},\overline{X}}
(\Upsilon_{Y, Z}\tens\id)\Upsilon_{X,Y\tens Z}\overline{\Phi_{X,Y,Z}}\,
&=&(\id\tens \Upsilon_{X,Y})\Upsilon_{X\tens Y, Z}\ .
\end{eqnarray*}

\noindent
c)\quad $\overline{\star}\star=\mathrm{bb}_{1_{\CC}}:
1_{\CC}\to \overline{\overline{1_{\CC}}}$.

\noindent
d)\quad $\overline{\mathrm{bb}_X}=\mathrm{bb}_{\bar X}:\overline{X}\mapsto 
\overline{\overline{\overline{X}}}$\ .
\end{defin}

The natural transformation property for $\Upsilon$ is explicitly that, given morphisms $\theta:X\to B$ and $\phi:Y\to C$, we have
\begin{eqnarray*}
\Upsilon_{B,C}\, (\overline{\theta\tens\phi}) \,=\,
(\overline{\phi}\tens\overline{\theta})\, \Upsilon_{X,Y} : \overline{X\tens Y} \to
\overline{C}\tens\overline{B}\ .
\end{eqnarray*}
The natural transformation property of $\mathrm{bb}$ is similarly stated for any morphism $\theta$ as
\[ \mathrm{bb}_{B}\circ\theta=\overline{\overline{\theta}}\circ\mathrm{bb}_X.\]

\begin{remark}
We can define $\CC^{\rm op}$ to be the same monoidal category as $\CC$ with opposite tensor product, the same unit object $1_\CC$ and the induced remaining structures. Then a bar category means:

1)  a monoidal category $\CC$ equipped with a monoidal functor ${\rm bar}:\CC\to \CC^{\rm op}$ with associated natural equivalence $\Upsilon^{-1}$ and induced unit isomorphism $\star:1_\CC\to \overline{1_\CC}$ 

2) a natural equivalence $\mathrm{bb}:\id\to{\rm bar}\circ{\rm bar}$ of functors such that 4),c),d) hold. Note that this is as functors not as monoidal functors (see Example~\ref{Copbar}), otherwise it would  be better to write $\mathrm{bb}:\id\to {\rm bar}^{\rm op}\circ {\rm bar}$ here. 
\end{remark}

\begin{defin}\label{strictbarcat} In any bar category we define the associated
natural `coboundary' isomorphism $\Xi:\tens\to\tens$ by 

\begin{picture}(100,80)(-100,6)
\put(140,63){\vector(1,0){45}}
\put(60,23){\vector(1,0){125}}
\put(38,55){\vector(0,-1){20}}
\put(207,35){\vector(0,1){20}}
\put(20,60){$X\tens Y$}
\put(100,60){$X\tens Y$}
\put(50,63){\vector(1,0){45}}
\put(20,20){$\overline{\overline{X\tens Y}}$}
\put(191,60){$\overline{\overline{X}}\tens \overline{\overline{Y}}$}
\put(191,20){$\overline{\overline{Y}\tens \overline{X}}$}

\put(135,70){$\mathrm{bb}_{X}\tens\mathrm{bb}_{Y}$}
\put(65,70){$\Xi_{X,Y}$}
\put(107,30){$\overline{\Upsilon_{X,Y}}$}
\put(-1,43){$\mathrm{bb}_{X\tens Y}$}
\put(213,43){$\Upsilon_{\overline{Y},\overline{X}}$}
\end{picture}

The bar category is called {\em strong} if $\Xi$ is the identity. 
\end{defin} 

The terminology here is similar to that for Hopf algebras\cite{MajBook} or group cohomology.  Similarly, we have:

\begin{propos}\label{Xi} $\Xi$ is a `cocycle' in the sense $\Xi_{X,1_\CC}=\id$, $\Xi_{1_\CC,X}=\id$ and
\[ \Phi_{X,Y,Z}(\Xi_{X,Y}\tens\id)\Xi_{X\tens Y,Z}=(\id\tens\Xi_{Y,Z})\Xi_{X,Y\tens Z}\Phi_{X,Y,Z}.\]
\end{propos}
\begin{proof} For the first part, we apply the definition and axioms 4a) of a bar category. For  the last part we expand $\Xi_{X\tens Y,Z}$ using the above in terms of $\mathrm{bb}_{X\tens Y}\tens \mathrm{bb}_Z$ and replace the latter in favour of $\Xi_{X,Y}$. Doing the same on the other side for $\Xi_{X,Y\tens Z}$ and comparing, using the axiom 4 b) of a monoidal category twice gives the second result. \end{proof}

Just as a braiding is a natural general notion for additional structure on a monoidal category, while the special case of  symmetric monoidal is of significant interest in its own right, similarly the coboundary natural isomorphism is useful (and so we allow it in the general definition of a bar category) but the strong case is adequate for many basic examples.

\subsection{Basic Example: Vector spaces over $\C$}

Given a complex vector space $V$ we define $\overline{ V}$ to be the same Abelian group as $V$ but with a conjugate action of $\C$ so 
\[ \lambda.\overline{v}=\overline{\overline{\lambda}.v},\quad v\in V, \lambda\in \C.\]
Here $\overline v$ denotes the same element $v\in V$ but viewed in $\overline{V}$.  This defines a functor from the category of complex vector spaces to itself. Objects are sent as above, and a morphism 
$\phi:V\to W$ is sent to $\overline{\phi}:\overline{V}\to \overline{W}$ with 
$\overline{\phi}(\bar v)=\overline{\phi(v)}$ for all $v\in V$. 

Vector spaces over any field form a monoidal category with $\Phi,l,r$ trivial and unit object the field. In our case the unit object is $1_{\bf Vect}=\C$ and $\star$ is given by complex conjugation which we also denote as usual. One may easily verify all the axioms with the trivial choices
\begin{eqnarray}\label{serth}
\mathrm{bb}_V(v) \,=\,\overline{\overline{v}}\ ,\quad \Upsilon_{W,V}(\overline{w\tens v})
\,=\, \overline{v}\tens \overline{w}\ .
\end{eqnarray}
In this way the category {\bf Vect} of complex vector spaces becomes a strong bar category.

\subsection{Basic Example: Bimodules over a star algebra}\label{bimodd}

Let $A$ be an algebra. The
category ${}_A{\CM}_A$ of $A$-bimodules and $A$-bimodule maps forms a monoidal category using $\tens_A$. 
The left/right action on a tensor product of bimodules is given by multiplication
on the  left/right most factor. The associator is trivial and the  identity object is $A$. 

Now suppose that $A$ is a $*$-algebra  over $\C$ (or over some other field with involution) in the sense of an antilinear order reversing involution $a\mapsto a^*$.  For an $A$-bimodule $E$, define another $A$-bimodule $\overline{E}$ as follows: As abelian groups 
$\overline{E}$ is identified with $E$ (as for the bar-category {\bf Vect}; as there we use an overline 
  to distinguish $\overline{e}\in \overline{E}$ from $e\in E$.)

Define left and right actions of $A$ on 
$\overline{E}$ by $a.\overline{e}=\overline{e.a^*}$ and  $\overline{e}.a=\overline{a^*.e}$. 
The functor $\mathrm{ bar}:{}_A{\CM}_A \to {}_A{\CM}_A$
 is defined on objects by
$E\mapsto \overline{E}$ and on morphisms by $\phi\mapsto\overline{\phi}$, where
$\overline{\phi}(\overline{e})=\overline{\phi(e)}$. 

Considering $A$ as an $A$-bimodule, define the morphism 
\[ \star:A\to \overline{A},\quad a\mapsto\overline{a^*}.\]
Similarly, define $\Upsilon_{F,E}: \overline{F\tens_A E} \to\overline{E} \tens_A\overline{F}$,
 by 
$\Upsilon_{F,E}(\overline{f\tens_A e})=\overline{e}\tens_A\overline{f}$.
 It is easily seen that $\Upsilon_{F,E}$ is an
invertible morphism in ${}_A{\CM}_A$:
\begin{eqnarray*}
\Upsilon((\overline{f\tens_A e}).a) \,=\, \Upsilon(\overline{a^*.f\tens_A e}) \,=\,
\overline{e}\tens_A \overline{a^*.f} \,=\, \overline{e}\tens_A \overline{f}.a\ ,
\end{eqnarray*}
and likewise for the left multiplication. 

Finally, the isomorphism $\mathrm{bb}_E:E\to \overline{\overline{E}}$ is given by
$\mathrm{bb}_E(e)=\overline{\overline{e}}$, and
this is a natural transformation between the identity and the bar${}^2$ functor. 
Hence the monoidal category of bimodules becomes
a strong bar category in a natural way. The $*$-structure appears naturally as the $\star$ isomorphism on $A$ viewed as unit object of the category of $A$-bimodules and its opposite.

\subsection{Basic Example: (co)modules over $*$-Hopf algebras}\label{modbar}

Let $H$ be a Hopf algebra over a field $k$ in the usual sense, with coproduct $\Delta:H\to H\tens H$, counit $\eps:H\to k$, and antipode $S:H\to H$. We refer to \cite{MajBook} for details and we use the `Sweedler notation' $\Delta h=h\o\tens h\t$ for all $h\in H$. The category ${\CM}_H$ of right $H$-modules is given the tensor product action
of $H$ on tensor products, i.e.\ 
\begin{eqnarray*}
(v\tens w)\ra h \,=\, v\ra h_{(1)}\tens w \ra h_{(2)}\ .
\end{eqnarray*}

A natural notion over $\C$ (or over any field with involution) is that $H$ is a $*$-algebra and that
$\Delta$ is a $*$-algebra map
\[\Delta\circ *=(*\tens *)\circ\Delta.\]
{}From this it follows that $S$ is necessarily invertible and 
\[ (Sh)^*=S^{-1}(h^*)\ .\]
We say that $H$ in this case is a {\em $*$-Hopf algebra}. The notion has been used notably by
Woronowicz as an algebraic notion that can in some cases be completed to a compact quantum group\cite{woron87}. There are other possibilities (see Section~\ref{fliphopf} and~\ref{kfkfk} later).

In this case ${\CM}_H$ can be made into a strong bar
category as follows: As in \ref{bimodd} we identify $\overline{V}$ with $V$ as sets,
using overline to distinguish $\overline{v}\in \overline{V}$ from $v\in V$. 
Then the action $\ra$ of $H$ on $\overline{V}$ is given by
\begin{eqnarray*}
\overline{v}\ra h \,=\, \overline{v\ra S^{-1}(h^*)}\ .
\end{eqnarray*}
We take the same trivial $\Upsilon,\mathrm{bb}$ as in (\ref{serth}) for {\bf Vect}  and the same $\star$ on the unit object $\C$ given by complex conjugation.  We check that $\bb$ is a right module map:
 \begin{eqnarray*}
 \mathrm{bb}(v)\ra h &=& \overline{\overline{v}}\ra h \,=\, 
 \overline{\overline{v}\ra S^{-1}(h^*)} \cr
 &=&  \overline{\overline{ v\ra S^{-1}((S^{-1}(h^*))^*)}} \cr
 &=& \overline{\overline{ v\ra h}} \,=\, \mathrm{bb}(v\ra h)\ .
\end{eqnarray*}
We also check that $\Upsilon$ is a right module map,
\begin{eqnarray*}
\Upsilon((\overline{w\tens v})\ra h) &=& 
\Upsilon(\overline{(w\tens v)\ra S^{-1}(h^*)}) \cr
&=& \Upsilon(\overline{ w\ra S^{-1}(h^*)_{(1)}\tens v\ra S^{-1}(h^*)_{(2)}}) \cr
&=&  \Upsilon(\overline{ w\ra S^{-1}(h^*_{(2)})\tens v\ra S^{-1}(h^*_{(1)})}) \cr
&=& \overline{ v\ra S^{-1}(h^*_{(1)})} \tens \overline{ w\ra S^{-1}(h^*_{(2)})} \cr
&=& \overline{ v}\ra h_{(1)} \tens \overline{ w}\ra h_{(2)} \cr
&=& \Upsilon(\overline{w\tens v}) \ra h\ .
\end{eqnarray*}

Analogous results hold for the category ${}_H{\CM}$ of left $H$-modules. Here tensor product of modules is given now by the left action of $\Delta h$. We specify the left action on $\overline{V}$ by\begin{eqnarray*}
h\la \overline{v} \,=\, \overline{S^{-1}(h^*)\la v}\ .
\end{eqnarray*}
and again the trivial forms for $\star,\Upsilon,\mathrm{bb}$.

Since Hopf algebras are in a certain sense `self dual' under arrow reversal there are analogous results for comodule categories. Thus the category ${\CM}^H$ of right $H$-comodules (say) has objects
which are comodules $(V,\rho_V)$ where 
\[ \rho_V:V\to V\tens H,\quad \rho(v)=v_{[0]}\tens v_{[1]}\]
is a coaction in an explicit notation. Again we identify $V$ and $\overline{V}$ as sets, and define
$\star,\mathrm{bb},\Upsilon$ as in for {\bf Vect}. The coaction on $\overline{V}$
is given by
\begin{eqnarray*}
\rho(\overline{v}) \,=\, \overline{v_{[0]}}\tens v_{[1]}^*\  .
\end{eqnarray*}
We check that $\bb$ is a right comodule map:
\begin{eqnarray*}
\rho\,\mathrm{bb}(v) &=& \rho(\overline{\overline{v}}) \,=\,
\overline{(\overline{v})_{[0]}} \tens (\overline{v})_{[1]}^* \cr
&=& \overline{\overline{v_{[0]}}} \tens v_{[1]}^{**}\,=\, 
 \overline{\overline{v_{[0]}}} \tens v_{[1]}\end{eqnarray*}
 We also check that $\Upsilon$ is a right comodule map,
\begin{eqnarray*}
(\Upsilon\tens\id)\rho(\overline{w\tens v}) &=& 
\Upsilon(\overline{(w\tens v)_{[0]}}) \tens
((w\tens v)_{[1]})^* \cr
&=& \Upsilon(\overline{w_{[0]}\tens v_{[0]}}) \tens (w_{[1]}v_{[1]})^* \cr
&=& \overline{v_{[0]}}\tens \overline{w_{[0]}} \tens v_{[1]}^*w_{[1]}^*\ ,\cr
\rho\Upsilon(\overline{w\tens v}) &=& \rho(\overline{v}\tens\overline{w}) \cr
&=& \overline{v_{[0]}}\tens \overline{w_{[0]}} \tens v_{[1]}^*w_{[1]}^*\ ,
\end{eqnarray*}

In this way ${\CM}^H$ becomes a strong bar category when $H$ is a $*$-Hopf algebra. Similarly for the category ${}^H{\CM}$ of left comodules.

\subsection{Functors between bar categories}\label{subsectfunc}

\begin{defin} \label{recdef}
Given two bar categories ${\CC}$ and ${\CD}$, a functor from
${\CC}$ to ${\CD}$ is a monoidal functor $F:{\CC}\to{\CD}$
with the additional properties:

\medskip \noindent(1)\quad There is a natural equivalence $\mathrm{fb}$ between the functors $\mathrm{bar}\, F$ and $F\,\mathrm{bar}$, i.e.\ 
$\mathrm{fb}_Y:\overline{F(Y)}\to F(\overline{Y})$.

\smallskip \noindent(2)\quad 
For the unit, $F(\overline 1)=\overline{F(1)}$, and for the
 $\star:1\to\overline{1}$ morphism,
$F(\star)=\star$ (more precisely,  the latter should $F(\star)f^1=\overline{f^1}\star$ if we do not suppress $f^1$)

\smallskip \noindent(3)\quad $F(\mathrm{bb}_Y)=
\mathrm{fb}_{\bar Y}\circ \overline{\mathrm{fb}_Y}\circ \mathrm{bb}_{F(Y)}$.

\smallskip \noindent(4)\quad The following diagram commutes, where
we recall that the definition of a monoidal functor involves a natural equivalence
$f_{X,Y}:F(X)\tens F(Y)\to F(X\tens Y)$ which respects the associator:

\begin{picture}(100,80)(-50,6)

\put(60,63){\vector(1,0){45}}
\put(60,23){\vector(1,0){45}}
\put(38,55){\vector(0,-1){20}}
\put(10,60){$\overline{F(X\tens Y)}$}
\put(10,20){$F(\overline{X\tens Y})$}
\put(111,60){$\overline{F(X)\tens F(Y)}$}
\put(111,20){$F(\overline{Y}\tens\overline{ X})$}

\put(65,70){$\overline{f^{-1}_{X,Y}}$}
\put(67,30){$F(\Upsilon)$}
\put(3,43){$\mathrm{fb}_{X\tens Y}$}

\put(175,63){\vector(1,0){35}}
\put(165,23){\vector(1,0){45}}
\put(227,55){\vector(0,-1){20}}

\put(233,43){$\mathrm{fb}_{Y}\tens \mathrm{fb}_{X}$}
\put(175,30){$f_{\overline{Y},\overline{X}}^{-1}$}
\put(213,20){$F(\overline{Y})\tens F(\overline{ X})$}
\put(213,60){$\overline{F(Y)}\tens \overline{F(X)}$}
\put(185,70){$\Upsilon$}
\end{picture}
\end{defin}

\begin{propos} \label{kxzj}
Given a bar functor $(F,f)$ between bar categories ${\CC}$
and ${\CD}$ as in Proposition~\ref{recdef}, the following diagram commutes:

\begin{picture}(100,80)(-100,6)

\put(60,63){\vector(1,0){45}}
\put(60,23){\vector(1,0){45}}
\put(38,55){\vector(0,-1){20}}
\put(127,55){\vector(0,-1){20}}
\put(20,60){$\overline{F(X)}$}
\put(20,20){$ \overline{F(\overline{\overline{X}})}$}
\put(111,60){$F(\overline{X})$}
\put(111,20){$\overline{\overline{F(\overline{X})}}$}

\put(75,70){$\mathrm{fb}$}
\put(67,30){$\overline{\mathrm{fb}^{-1}}$}
\put(-1,43){$\overline{F(\mathrm{bb})}$}
\put(133,43){$\mathrm{bb}$}

\end{picture}

\end{propos}
\begin{proof} By using  (3) from  Definition~\ref{recdef} and (d) from
\ref{barcatdef}, the diagram can be rewritten as: 

\begin{picture}(100,80)(-100,6)

\put(60,63){\vector(1,0){45}}
\put(60,23){\vector(1,0){45}}
\put(38,55){\vector(0,-1){20}}
\put(127,55){\vector(0,-1){20}}
\put(20,60){$\overline{F(X)}$}
\put(20,20){$ \overline{\overline{\overline{F(X)}}}$}
\put(111,60){$F(\overline{X})$}
\put(111,20){$\overline{\overline{F(\overline{X})}}$}

\put(75,70){$\mathrm{fb}$}
\put(75,30){$\overline{\overline{\mathrm{fb}}}$}
\put(15,43){$\mathrm{bb}$}
\put(133,43){$\mathrm{bb}$}

\end{picture}

\noindent
This commutes as $\mathrm{bb}$ is a natural transformation between
the identity and $\mathrm{bar}^2$. \end{proof}

\medskip

Note that in axiom (1) of a functor between bar categories we should more precisely write that $\mathrm{fb}:\mathrm{bar}F\to F^{\rm op}\mathrm{bar}$ is a natural equivalence where $F^{\rm op}:\CC^{\rm op}\to \CD^{\rm op}$ is the same functor $F$ on objects but has $f^{op}_{X,Y}=f_{Y,X}$.

\begin{example}\label{Copbar}. If  $\CC$ is a bar category then $\CC^{\rm op}$ is also, with $\mathrm{bar}^{\rm op}$ as the bar functor for this. Then $\mathrm{bar}:\CC\to \CC^{\rm op}$ becomes automatically a bar functor between bar categories. This is the content of the axiom 4 d) of a bar category. \end{example}

To explain the categorical picture  further, let us note that for any natural  equivalence $\mathrm{b}:F\to G$ of two functors $F,G:\CC\to \CD$ between monoidal categories, with associated associated natural equivalences $f,g$ respectively, we can  define a `coboundary' $\del(\mathrm{b})_{X,Y}:F(X)\tens F(Y)\to F(X)\tens F(Y)$ by
\[ g_{X,Y}(b_X\tens b_Y)\del(\mathrm{b})_{X,Y}=b_{X\tens Y}f_{X,Y},\quad \forall X,Y\in\CC.\]
It is a 2-cocycle in a similar manner to Proposition~\ref{Xi}. We say that the natural equivalence is {\em monoidal} if $\del(\mathrm{b})=\id$.

For a general bar category we have already considered $\mathrm{bb}:\id\to \mathrm{bar}^{\rm op}\circ\mathrm{bar}$ with coboundary $\Xi=\del(\mathrm{bb})$. The case of a strong bar category is when $\mathrm{bb}$ is a monoidal equivalence of monoidal functors, not just an equivalence of functors.  By contrast, the meaning of axiom (4) of a bar functor between bar categories is that we ask for $\mathrm{fb}$ to be a monoidal equivalence.

\begin{propos}\label{Xipreserved}
If $F:\CC\to \CD$ is a bar functor between bar categories, the coboundary of $\CC$ maps to the coboundary of $\CD$, i.e. $F(\Xi_{X,Y})f_{X,Y}=f_{X,Y}\Xi_{F(X),F(Y)}$ for all $X,Y\in \CC$.
\end{propos}
\begin{proof} The proof is tedious `diagram filling' which we leave to the reader. The conceptual explanation is that the natural equivalence $\mathrm{fb}$ does not have any coboundary of its own, but connects that of $\CC$ to that of $\CD$. \end{proof}

This is similar to the idea that a tensor functor between braided monoidal categories is required to send the braiding of one to the braiding  of the other. 

\begin{corol}For any bar category, $\Upsilon_{X,Y}\overline{\Xi_{X,Y}}=\Xi_{\overline Y,\overline X}\Upsilon_{X,Y}$ for all objects $X,Y$.\end{corol}
\begin{proof} This follows form the above proposition applied  to Example~\ref{Copbar} above. It can also be proven easily enough by direct computation. \end{proof}

\medskip
Clearly, a concrete example of a bar functor is the forgetful functor to {\bf Vect} from any of the above four bar categories associated to a $*$-Hopf algebra $H$. This forgets the (co)action of $H$. In the appendix we prove the converse that if a bar functor to {\bf Vect} when restricted to objects with duals is representable, it factors through such a construction for some $*$-Hopf algebra which may be reconstructed.

\subsection{$\star$-Objects in a bar-category} \label{jkhdsz}

In any bar category it makes sense to ask if an object $X$ is isomorphic to $\overline{X}$. 

\begin{defin}\label{psdhvsd}
An object $X$ in a bar category is called a star object if it is equipped with  a morphism
$\star_X:X \to \overline{X}$ with the following properties:

\noindent a)\quad $\overline{\star_X} \circ \star_X=\mathrm{bb}_X$.

\noindent b)\quad $(\star_X\tens\id)r_X=r_{\bar X}\,\star_X$
and $(\id\tens\star_X)l_X=l_{\bar X}\,\star_X$. 
 
\end{defin}

 Clearly $1_\CC$ is always an example. 

\begin{example} In the bar category {\bf Vect} a $\star$ object structure on a vector space $V$ is given by an antilinear map $J:V\to V$ such that $J^2=\id$.  Here $\star(v)=\overline{J(v)}$ is $J$ viewed as a linear map $V\to \overline{V}$. 
We will usually write $J(v)=v^*$. 
\end{example}

\begin{example}
In ${}_H{\CM}$ (see \ref{modbar}) a $\star$ structure is given in the same way as in {\bf Vect} but with $J$ covariant in the sense $J(h\la v)=(Sh)^*\la J(v)$. This is the notion of a `$*$-action' or unitary action in \cite{MajBook} (where $J$ is denoted by $*$).  
\end{example}

\begin{example}
In ${}_A{\CM}_A$ (see \ref{bimodd}) we again
 use the identification of $\overline{X}$ as a bimodule with the Abelian group $X$ but reverse bimodule structures. Then 
$\star_X:X\to \overline{X}$ takes the form $e\mapsto \overline{J(e)}$ where $J^2=\id$ and now $J:X\to X$ is a skew bimodule map in the sense $J(a.e)=J(e).a^*$ and $J(e.a)=a^*.J(e)$ for $a\in A$ and $e\in X$. \end{example}

We can also write some other related definitions in any bar category. 

\begin{defin} \label{inndeff}
An inner product on an object $X$ a bar category means a morphism  $\<\ ,\ \>:\overline{X}\tens X\to 1_\CC$ such that 
\[ \star\circ\<\ ,\ \>=\overline{\<\ ,\ \>}\circ\Upsilon^{-1}_{\bar X,X}\circ(\id\tens\bb_X)\]
(This is just the usual conjugate symmetry for many examples.)
It will also be convenient to have the opposite definition with the bar on the
second copy of $X$, which we denote 
$\<\ ,\ \>_r:X\tens \overline{X}\to 1_\CC$ with the corresponding symmetry condition
\begin{eqnarray*}
\overline{\<,\>_r}\,\Upsilon^{-1}(\mathrm{bb}\tens\id) \,=\, \star\,
\<,\>_r : X\tens\overline{X}\to 
\overline{\mathbb{C}}\ .
\end{eqnarray*}
\end{defin}

\begin{example} In {\bf Vect} this becomes a usual sequilinear inner product $\<\ ,\ \>$ on a vector space $V$ that is antilinear in its first input and which is symmetric in the skew sense
\[ \overline{\<v,w\>}=\<w,v\>,\quad v,w\in V\]
where the bar on the left denotes usual complex conjugation (it is $*$ in {\bf Vect}). One should add a positivity requirement to have a (pre)-Hilbert space.
\end{example}

\begin{example} In bimodules
 ${}_A\CM_A$ over a $*$-algebra $A$, an $A$-valued inner product $\<\ ,\ \>$ on a bimodule $E$ obeys
\[ \<e,f\>^*=\<f,e\>,\quad \<e,f\>a=\<e,f.a\>,\quad  \<a^*.e,f\>=\<e,a.f\>\]
The first two conditions are those of a right (pre) Hilbert module if one adds a positivity requirement. Given the first, the second condition can also be written as $a\<e,f\>=\<e.a^*,f\>$. The last condition applies in our bimodule setting and says that additional left action of $A$ is adjointed by $\<\ ,\ \>$ in a way that agrees with $*$ (a `unitary' representation in some sense). From our point of view it is the condition that the map $\overline{E}\times E\to A$ descends to $\overline{E}\tens_A E\to A$.
\end{example}

\begin{example}
In the case ${}_H\CM$ for a $*$-Hopf algebra it means a sesquilinear inner product as for {\bf Vect} which is covariant under the action of $H$ on $\overline{V}\tens V$. A short computation shows that this last is equivalent to 
\[ \<h^*\la v,w\>=\<v, h\la w\>,\quad h\in H,\ v,w\in V\]
i.e. that the action of $H$ is adjointed by $\<\ ,\ \>$ in a way that agrees with $*$ (i.e. `unitary' in some sense).
\end{example}

\begin{defin} $B$ is a star algebra in the bar category ${\CC}$ if there
is a product $\mu:B\tens B\to B$ making $B$ into an algebra
in ${\CC}$, and if $B$ is also a star object so that
\begin{eqnarray*}
\overline{\mu}\,\Upsilon^{-1}_{B,B}(\star_B\tens \star_B)\,=\,\star_B\mu:B\tens B\to\overline{B}\ .
\end{eqnarray*}
If $\eta:1_\CC\to B$ is a unit for the product, we require $\overline{\eta}\circ\star=\star_B\circ \eta$.
\end{defin}

It is intended here that the product is associative with respect to a possibly non-trivial associator $\Phi$ in the category. 

\begin{example}
In the case of {\bf Vect} this means that $B$ is an ordinary $*$-algebra with 
$\star(b)=\overline{b^*}$, i.e.\
 the antilinear map $J$ above is antimultiplicative with respect to the product (and in this context denoted $*$).   An example is of course bounded operators $B(\CH)$ on a Hilbert space $\CH$, with $*$ given by the adjoint operation.
\end{example}

\begin{example} An algebra $B$ in ${}_A\CM_A$ means a bimodule $B$ and a product map $\mu:B\tens_AB\to B$ which need be associative only when one takes the tensor product over $A$.  A unit means a bimodule map $A\to B$, determined by the image of the identity, and forming a unit for $\mu$.  We provide the notion of a $\star$-involution on such an algebraic structure.
\end{example}

\begin{example}
In the case  ${}_H\CM$, an algebra $B$ means an $H$-module algebra in the sense $h\la(bc)=(h\o\la b)(h\t\la c)$ and $h\la 1=\eps(h)1$. We require further that $B$ is a $*$-algebra in the usual sense as for {\bf Vect} and that $(h\la b)^*=(Sh)^*\la b^*$ as above. This is the usual notion for a $*$-Hopf algebra acting on a $*$-algebra \cite{MajBook}. Here
we again have $\star(b)=\overline{b^*}$. 
\end{example}

\begin{defin} \label{dkjshbcvb} Suppose that 
 $B$ is a  star algebra in a bar category  ${\CC}$ and that the object $H$
 has an inner product $\<,\>_r:H\tens\overline{H}\to 1_{\mathbb{C}}$. A right action
 $\ra:H\tens B\to H$ is  compatible with the inner product
 if
\begin{eqnarray*}
\<,\>_r(\ra\tens\id)\,=\,
\<,\>_r(\id\tens\overline{\ra}\,\Upsilon^{-1}(\star_B\tens\id))\Phi
:( H \tens B)\tens\overline{H}\to  1_{{\CC}}\ .
\end{eqnarray*}
Alternatively the  inner product $\<,\>:\overline{H}\tens H\to 1_{\mathbb{C}}$
is compatible with the left action  $\la:B\tens H\to H$ if
\begin{eqnarray*}
\<,\>(\id\tens\la)\,=\,
\<,\>(\overline{\la}\,\Upsilon^{-1}(\id\tens\star_B)\tens\id)\Phi^{-1}
:\overline{H}\tens ( B \tens H)\to  1_{{\CC}}\ .
\end{eqnarray*}
\end{defin}

There is obviously a dual notion of a coalgebra $C$ in a bar category $\CC$ as a coalgebra with coproduct $\Delta:C\to C\tens C$ and counit $\eps:C\to 1_\CC$ which are coassociative and counital in in the monoidal category, now with $C$ a $\star$-object such that
\[ \Upsilon_{C,C}\circ\overline{\Delta}\circ\star_C=(\star_C\tens\star_C)\circ\Delta,\quad \star\circ\eps=\overline{\eps}\circ\star_C.\]
This is not the only possible equation, we could add a braiding or similar map. 

\begin{example} A coalgebra in ${}_A\CM_A$  means a `coring' \cite{BrzWis}. Hence the above is a natural notion of $\star$-structure on a coring.
\end{example}

This completes the basic  definitions and `elementary examples' in the sense that one does not need the notion of a bar category to formulate these examples. Rather it is a useful idea to unify various notions and to keep track of the book-keeping.

\section{New examples of bar categories}

Until now we have tested our ideas out on `standard' settings where one knows in any case how
conjugation etc.\ should proceed. In that case our bar notion serves only to unify some slightly different contexts. In this section we give some genuinely new examples where the definitions would not be
so clear without a categorical context. The first  one is only slightly new (it is
 known in some form) but included for completeness.

\subsection{Example: Crossed or Drinfeld-Yetter modules over $*$-Hopf 
algebras}\label{YDbar}

For any Hopf algebra the category of right `crossed' or 
Drinfeld-Yetter modules $\YD{}_H^H$ for $H$ consists of right $H$-modules and
right $H$-comodules satisfying the consistency condition
\begin{eqnarray} \label{ydprop}
\rho(v\ra h)=v_{[0]} \ra
h_{(2)}
\tens S(h_{(1)})\,v_{[1]}\,h_{(3)}\quad\forall v\in V\ ,\ \forall a\in H\ .
\end{eqnarray}
In the case of $H$ finite-dimensional a crossed module is just the same as a right module for the Drinfeld quantum double Hopf algebra $D(H)=H\bowtie H^{*\rm op}$ and its properties tend to be in line with those known for this.  In the finite-dimensional case it was shown \cite{majquastpoin} that $D(H)$ for $H$ a $*$-Hopf algebra inherits a natural $*$-Hopf algebra structure. In the general case:

\begin{propos} Let $H$ be a $*$-Hopf algebra. The category of right $H$-crossed module or `Drinfeld-Yetter' modules is a strong bar category with $\overline{ V}$ defined as for $\CM_H$ and $\CM^H$ separately, i.e. $\overline{V}=V$ as an Abelian group and with action and coaction on $\overline{V}$ given by
\begin{eqnarray*}
\rho(\overline{v}) \,=\, \overline{v_{[0]}}\tens v_{[1]}^*\ ,\quad
\overline{v} \ra h \,=\, \overline{v\ra S^{-1}(h^*)}\ .
\end{eqnarray*}
\end{propos}
\begin{proof}
We only have to check that  $\overline{V}$ is again a crossed or Drinfeld-Yetter module, i.e.\ 
\begin{eqnarray*}
\rho(\overline{v} \ra h) &=& \rho(\overline{v\ra S^{-1}(h^*)}) \cr
&=& \overline{(v\ra S^{-1}(h^*))_{[0]}} \tens ((v\ra S^{-1}(h^*))_{[1]})^* \cr
&=& \overline{v_{[0]}\ra S^{-1}(h^*)_{(2)}} \tens 
(S(S^{-1}(h^*)_{(1)})v_{[1]} S^{-1}(h^*)_{[3]})^* \cr
&=&  \overline{v_{[0]}\ra S^{-1}(h^*_{(2)})} \tens 
(h^*_{(3)}v_{[1]} S^{-1}(h^*_{[1]}))^* \cr
&=& \overline{v_{[0]}}\ra h_{(2)} \tens S(h_{[1]})v_{[1]}^*h_{(3)}\ .
\end{eqnarray*}
The rest of the axioms hold as inherited from {\bf Vect}. \end{proof}

\subsection{Example: Bimodules over flip-$*$-Hopf algebras}\label{fliphopf}

Sometimes one has Hopf algebras which are simply not $*$-Hopf algebras but for which 
$*$ still behaves well with respect to the coproduct. The most common `nonstandard' variant
is:

\begin{defin}
A flip-$*$-Hopf algebra is a Hopf algebra which is a $*$-algebra and for which the
coproduct and counit obey
\[ \Delta\circ *=(*\tens *)\tau\circ\Delta,\quad \overline{\eps(h)}=\eps(h^*)\ ,\]
where $\tau$ is transposition. In this case it follows
that $S(h)^*=S(h^*)$. 
\end{defin}

A simple commutative example is the polynomial algebra $\C[u^i{}_j]$
($1\le i,j\le n$) with the relation $\det(u)=1$, coproduct, counit and $*$
\[ \Delta u^i{}_j=\sum_k u^i{}_k\tens u^k{}_j,\quad \eps u^i{}_j=\delta^i{}_j,\quad (u^i{}_j)^*=u^j{}_i.\]
This is a non-standard real form of $\C[SL_n]$ which does not correspond to any group but to the set of hermitian matrices of determinant 1. For $n=2$ it is the $*$-algebra of regular functions on the mass hyperboloid in Minkowski space $\R^{1,3}$ (and it is the $q=1$ limit of the hyperboloid in an approach to $q$-deformed Minkowski space based on braided groups \cite{MajBook}.)

\begin{propos} For any Hopf algebra one has a category  ${}_H{\CM}_H$ of $H$-bimodules with tensor
product defined by $\Delta$. If $H$ is a flip-$*$-Hopf algebra this can made into a strong bar 
category by $\overline{V}=V$ as Abelian groups, with right and left actions on $\overline{V}$ given by
\begin{eqnarray*}
h\la \overline{v} \,=\, \overline{v\ra h^*}\ ,\quad
 \overline{v}\ra h \,=\, \overline{h^*\la v   }\ .
 \end{eqnarray*}
 and  $\mathrm{bb}$ and $\Upsilon$ as  in (\ref{serth}).
 \end{propos}
\begin{proof} We need to check that $\mathrm{bb}$ and $\Upsilon$ give morphisms,
leaving some of this calculation to the reader. 
\begin{eqnarray*}
h\la\mathrm{bb}(v) &=& h\la  \overline{\overline{v}} \cr
&=&  \overline{\overline{v}\ra h^*} \,=\, \overline{\overline{h^{**}\la v}}\ , \cr
h\la\Upsilon(\overline{v\tens w}) &=& h\la(\overline{w}\tens \overline{v}) \cr
&=& h_{(1)}\la\overline{w}\tens h_{(2)}\la\overline{v} \cr
&=& \overline{w\ra h_{(1)}^{\phantom{(1)}*}}\tens
 \overline{v\ra h_{(2)}^{\phantom{(2)}*}} \cr
 &=& \overline{w\ra h^{*}_{\phantom{*}(2)}}\tens
 \overline{v\ra h^{*}_{\phantom{*}(1)}}\ ,\cr
 \Upsilon(h\la(\overline{v\tens w})) &=&  \Upsilon(\overline{(v\tens w)\ra h^*})\cr
 &=& \Upsilon(\overline{v\ra h^{*}_{\phantom{*}(1)}\tens w
 \ra h^{*}_{\phantom{*}(2)}})\cr
  &=& \overline{w\ra h^{*}_{\phantom{*}(2)}}\tens
 \overline{v\ra h^{*}_{\phantom{*}(1)}}\ .
\end{eqnarray*}
\end{proof}

\subsection{Example: Modules over quasi-$*$ Hopf algebras} \label{kfkfk}

We recall \cite{MajBook} that a 2-cocycle on a Hopf algebra $H$ is an invertible  element $\CF\in H\tens H$ such that
\[ \CF_{12}(\Delta\tens \id)(\CF)=\CF_{23}(\id\tens\Delta)(\CF),\quad (\eps\tens\id)\CF=1\]
Given such an element, one may conjugate the coproduct of $H$ to a new Hopf algebra $H^{\CF}$ which an equivalent category of modules and we will do so later. At present we give a different application of such a cocycle and will denote it by $\CG$ rather than by $\CF$ for this reason and to avoid confusion. Our construction is not directly related to, but in the same spirit as,
 the notion of quasi-$*$ Hopf algebra introduced in \cite{majquastpoin} to describe $q$-deformed inhomogeneous quantum groups for real $q$. Our version will be more applicable to $|q|=1$.

\begin{defin} \label{qudef}
 A quasi-$*$ Hopf algebra is a quadruple $(H,*,\CG,\gamma)$
  where $H$ is a Hopf algebra and $*$-algebra, $\CG\in H\tens H$ is a 2-cocycle and $\gamma\in H$ is an invertible element such that $(S\gamma)^*=\gamma$ and 
 \[
\Delta(h^*)={\CG}^{-1}(\Delta\ h)^{*\tens *}{\CG},\quad
\eps(h )^*=\eps(h^*),\quad S^{-1}(h^*)=\gamma^{-1}(Sh)^*\,\gamma\, \]
for all $h\in H$. The quasi-$*$ Hopf algebra is called {\em strong} if in addition we have 
$(\gamma\tens \gamma)\Delta\gamma^{-1} =(*\tens *)(S\tens S)(\CG_{21})\CG$. 
\end{defin}

\begin{remark}
If $(*\tens *)(S\tens S)(\CG)=\CG_{21}^{-1}$ we say that the quasi-$*$ Hopf algebra 
is of real type. 
In this case the strongness condition reduces to $\gamma$  grouplike, i.e.\
$\Delta\gamma=\gamma\tens\gamma$. 
\end{remark}

\begin{theorem}
The category ${}_H{\CM}$ of left modules of a quasi-$*$ Hopf algebra
$(H,*,\CG,\gamma)$ can be made into a bar category as follows: 
Again, for $V$ an object in ${}_H{\CM}$, we define $\overline{V}$ to be 
$V$ as a set, with 
\begin{eqnarray*}
h\la \overline{v} \,=\, \overline{(Sh)^*\la v}\ .
\end{eqnarray*}
We define $\Upsilon_{V,W}:\overline{V\tens W}\to \overline{V} \tens \overline{W}$
and $\mathrm{bb}_V:V\to\overline{\overline{V}}$ by 
\begin{eqnarray}\label{quasiupsdef}
\Upsilon_{V,W}(\overline{v\tens w}) &=& \overline{{\CG}^{(2)}\la w}
\tens \overline{{\CG}^{(1)}\la v}\ ,\cr
\mathrm{bb}_V(v) &=& \overline{\overline{\gamma\la v}}\ .
\end{eqnarray}
Then the coboundary map $\Xi$ in Definition~\ref{strictbarcat} is
\begin{eqnarray*}
\Xi(v\tens w) \,=\, 
 \gamma^{-1}(S\tilde{\CG}^{(2)})^*{\CG}^{(1)}
\gamma_{(1)}\la v
 \tens \gamma^{-1}(S\tilde{\CG}^{(1)})^*{\CG}^{(2)}
\gamma_{(2)}\la w\ , 
\end{eqnarray*}
where $\tilde \CG$ is a second independent copy of $\CG$. The bar category ${}_H\CM$ is strong precisely when the quasi-$*$ Hopf algebra $H$ is. \end{theorem}
\begin{proof} We first show that 
$\mathrm{bb}_V:V\to\overline{\overline{V}}$
defined in (\ref{quasiupsdef}) is a morphism in ${}_H{\CM}$.
\begin{eqnarray*}
g\la \mathrm{bb}_V(v) &=& g\la \overline{\overline{\gamma\la v}} \cr
&=& \overline{(Sg)^*\la \overline{\gamma\la v}} \cr
&=& \overline{\overline{(S((Sg)^*))^*\,\gamma\la v}}\ ,
\end{eqnarray*}
so to be a morphism we require $(S(Sg)^*))^*=\gamma g \gamma^{-1}$ which is equivalent to the condition displayed. Moreover, the condition 4 d) of Definition~\ref{barcatdef} translates as $(S\gamma)^*=\gamma$. Next we show that 
$\Upsilon_{V,W}:\overline{V\tens W}\to \overline{W} \tens \overline{V}$
defined in (\ref{quasiupsdef}) is a morphism in ${}_H{\CM}$.
\begin{eqnarray*}
\Upsilon(g\la (\overline{v\tens w})) &=&\Upsilon( \overline{S(g)^*\la(v\tens w)}) \cr
&=& \Upsilon(\overline{S(g)^*_{\phantom{*}(1)}\la v\tens S(g)^*_{\phantom{*}(2)}\la w})\\
 &=& \overline{ \CG^{(2)} (Sg)^*_{\phantom{*}(2)}\la w}\tens \overline{\CG^{(1)}(Sg)^*_{\phantom{*}(1)}\la v}\\
g\la\Upsilon(\overline{v\tens w}) &=& g\la(\overline{ {\CG}^{(2)}\la w} \tens
\overline{ {\CG}^{(1)}\la v} ) \cr
&=& \overline{(Sg_{(1)})^* {\CG}^{(2)}\la w} \tens
\overline{(Sg_{(2)})^* {\CG}^{(1)}\la v} \ .
\end{eqnarray*}
Comparing the two, we require $\Delta (Sg)^*=\CG^{-1}((*S\tens *S)\tau\Delta g)\CG$ which is equivalent to the condition stated. Next we show that 
the compatibility condition between $\Upsilon$ and the associator
is satisfied.
\begin{eqnarray*}
(\Upsilon_{V,W}\tens\id)\Upsilon_{U,V\tens W}(\overline{u\tens(v\tens w)}) &=& 
\Upsilon_{V,W}(\overline{\CG^{(2)}_{\phantom{(2)}(1)}\la v  \tens \CG^{(2)}_{\phantom{(2)}(2)}\la w })
\tens \overline{\CG^{(1)}\la u} \cr
&=& 
(\overline{\tilde\CG^{(2)}\CG^{(2)}_{\phantom{(2)}(2)}\la w}  \tens 
\overline{\tilde\CG^{(1)}\CG^{(2)}_{\phantom{(2)}(1)}\la v })
\tens \overline{\CG^{(1)}\la u} \ ,\cr
(\id\tens\Upsilon_{U,V})\Upsilon_{U\tens V,W}
(\overline{(u\tens v)\tens w}) &=& \overline{\CG^{(2)}\la w}\tens \Upsilon_{U,V}
(\overline{\CG^{(1)}_{\phantom{(1)}(1)}\la u  \tens \CG^{(1)}_{\phantom{(1)}(2)}\la v }) \cr
&=& \overline{\CG^{(2)}\la w}\tens 
(\tilde\CG^{(2)}\overline{\CG^{(1)}_{\phantom{(1)}(2)}\la v}  \tens
 \overline{\tilde\CG^{(1)}\CG^{(1)}_{\phantom{(1)}(1)}\la u })\ ,
\end{eqnarray*}
As the underlying category has trivial associator, this is implied by
 the cocycle condition for ${\CG}\in H\tens H$.

Finally, we check the equation for 
 $\Xi$.
Begin with
\begin{eqnarray*}
\Upsilon\circ\overline{\Upsilon}\Big(\overline{\overline{v\tens w}}\Big) &=& 
\Upsilon\Big(\overline{\overline{{\CG}^{(2)}\la w}
\tens \overline{{\CG}^{(1)}\la v}}\Big) \cr
&=& \overline{\tilde{\CG}^{(2)}\la (\overline{{\CG}^{(1)}\la v})}
 \tens\overline{\tilde{\CG}^{(1)}\la (\overline{{\CG}^{(2)}\la w})} \cr
 &=& 
 \overline{\overline{(S\tilde{\CG}^{(2)})^*{\CG}^{(1)}\la v}}
 \tens\overline{\overline{(S\tilde{\CG}^{(1)})^*{\CG}^{(2)}\la w}}\ ,
\end{eqnarray*}
and use
\begin{eqnarray*}
\mathrm{bb}(v\tens w) &=& \overline{\overline{\gamma\la(v\tens w)}} \cr
&=& \overline{\overline{\gamma_{(1)}\la v\tens \gamma_{(2)}\la w}} 
\end{eqnarray*}\end{proof}

\medskip
To give a simple class of examples, we recall \cite{Dri} that a Hopf algebra $H$ is quasitriangular if there  is an invertible element $\CR\in H\tens H$ such that 
\[ (\Delta\tens\id)\CR=\CR_{13}\CR_{23},\quad (\id\tens\Delta)\CR=\CR_{13}\CR_{12},\quad\tau\circ\Delta=\CR(\Delta\ )\CR^{-1}\]
(which then obeys the YBE or braid relations in $H^{\tens 3}$). In this case the elements
\begin{eqnarray} \label{uvdef}
u=(S\CR\ut)\CR\uo\ ,\quad v=\CR\uo (S\CR\ut)
\end{eqnarray}
implement $S^2$ and $S^{-2}$ respectively by conjugation. A ribbon Hopf algebra is a quasitriangular one for which the central element $vu$ has a square root $\nu\in H$ (the ribbon element) such that
\[ \Delta \nu=\CR^{-1}\CR^{-1}_{21}(\nu\tens \nu) \ ,\quad \eps\nu=1, \quad S\nu=\nu\ .\]
An exposition with required identities can be found in \cite{MajBook}.

\begin{propos} \label{ufudsn77} A quasitriangular Hopf algebra $(H,\CR)$ which is a flip-$*$-Hopf algebra  can be viewed as a
quasi-$*$-Hopf algebra with $\CG=\CR$ and $\gamma=u^{-*}$ or $\gamma=v^*$. Then 
the map $\Xi$ on ${}_H{\CM}$ is given by the action of $(\gamma^{-1}\tens\gamma^{-1})(({*\tens *})(\CR_{21}))\CR\, \Delta\gamma$. 
\end{propos}
\begin{proof} Since we have assumed both  parts of 
\begin{eqnarray*}
(\Delta g)^{*\tens *}\,=\,\tau\Delta(g^*)\ ,\quad 
\tau\Delta g = {\CR}(\Delta g){\CR}^{-1}\ ,
\end{eqnarray*}
where $\tau$ is transposition, it is clear that $\CG=\CR$ plays the role required.  It remains only to  find $\gamma$. We start by noting
 that for a flip-$*$-Hopf algebra, $S(h^*)=(Sh)^*$, so the condition for $\gamma$ becomes 
$\gamma^{-*}g\gamma^*= S^2(g)$. 
{}From general analysis on quasi-triangular Hopf algebras \cite{MajBook},
this is satisfied by $\gamma=u^{-*}$ or $\gamma=v^*$. 
The form for $\Xi$ is given by using the fact that 
$(S\tens S)(\CR)=\CR$ for any quasitriangular Hopf algebra. \end{proof}

\begin{propos} \label{ufudsn} A ribbon Hopf algebra $(H,\CR,\nu)$ which is a flip-$*$-Hopf algebra such that ${\CR}^{*\tens *}={\CR}_{21}^{-1}$ can be viewed as a
strong quasi-$*$-Hopf algebra with $\CG=\CR$ and $\gamma=v^{-1}\nu$.
\end{propos}
\begin{proof} As in \ref{ufudsn77} we choose $\CG=\CR$, but we need a different $\gamma$
to make $\Xi$ the identity. We need a grouplike $\gamma$ with the property that 
$\gamma^{-*}g\gamma^*= S^2(g)$. 
 Since 
${\CR}^{*\tens *}={\CR}_{21}^{-1}$, we have 
$u^*=u^{-1}$ and $v^*=v^{-1}$ hence if we have a ribbon element, 
$\gamma=v^{-1}\nu$ is grouplike and still implements $S^2$ as required. \end{proof}

\begin{example}\label{anyhfhfh} The anyon Hopf algebra \cite{MajBook} 
for $q$ a primitive $l$-th root of unity ($l=2k+1>1$ an odd integer)
has a single invertible generator $K$ for which 
\begin{eqnarray*}
K^l \,=\,1\ ,\quad \Delta K=K\tens K ,\quad 
\epsilon(K)=1\quad S(K)=K^{-1}.
\end{eqnarray*}
This Hopf algebra is ribbon, with quasitriangular and associated elements
\begin{eqnarray}\label{kfjsdv}
{\CR}_K \,=\,\frac{1}{l}\sum_{a,b=0}^{l-1}q^{-2ab}\, K^{a}\tens K^{b},\quad \nu=u=v &=& \frac{1}{l}\sum_{a,r=0}^{l-1}q^{-2a(a-r)}\, K^{r} 
\end{eqnarray}
We take the star structure to be $K^*=K^{-1}$, and then
${\CR}^{*\tens *}={\CR}^{-1}={\CR}^{-1}_{21}$ and we have a quasi-$*$ Hopf algebra as an example of \ref{ufudsn} with $\gamma=1$.

The bar category of modules of this quasi-$*$ Hopf algebra is that of anyonic or $\Z_l$-graded complex vector spaces (here $K\la v=q^{|v|}v$ acting in an element of degree $|v|$). Morphisms are degree-preserving linear maps and
\[ \bb_V(v)=\overline{\overline{v}},\quad \Upsilon_{V,W}(\overline{v\tens w})=q^{2|v||w|}\overline{w}\tens\overline{v}\]
on homogeneous elements. 
\end{example}

\begin{example} \label{hjsdgfc}
As another example of the structure given in \ref{ufudsn},
take the example of ${u_q(su_2)}$, where $q$ is a primitive $l$th
root of unity (for odd integer $l>1$). The algebra is defined as generated by
$1$, $K^{\pm 1}$, $E$ and $F$, subject to relations
\begin{eqnarray*}
K\, E=q^2\, E\, K\ ,\quad K\, F=q^{-2}\, F\, K\ ,\quad
[E,F]=\frac{K-K^{-1}}{q-q^{-1}}\ ,\quad
E^l=0=F^l\ ,\quad K^{l}=1\ .
\end{eqnarray*}
The coproduct on the generators is given by
\begin{eqnarray*}
\Delta E=E\tens K+1\tens E\ ,\quad 
\Delta F=F\tens 1+K^{-1}\tens F\ ,\quad 
\Delta K=K\tens K\ .
\end{eqnarray*}
The counit is defined by
\begin{eqnarray*}
\epsilon(1)=\epsilon(K^{\pm 1})=1\ ,\ \epsilon(E)=\epsilon(F)=0\ ,
\end{eqnarray*}
and the antipode by
\begin{eqnarray*}
S(K^{\pm 1})=K^{\mp 1}\ ,\ S(E)=-E\,K^{-1}\ ,\ S(F)=-K\, F\ . 
\end{eqnarray*}
This Hopf algebra is known to be ribbon with quasitriangular and related structures
\begin{eqnarray*}
{\CR} &=& {\CR}_K
\sum_{n=0}^{l-1} \frac{(q-q^{-1})^n}{[n;q^{-2}]!}\, E^n\tens F^n\ ,\cr
\nu &=& \frac{1}{l}\bigg(\sum_{m=0}^{l-1}q^{2\, m^2}\bigg)\sum_{n,a=0}^{l-1}
\frac{(q-q^{-1})^n}{[n;q^{-2}]!}\, q^{-(l+1)(n-a-1)^2/2}\, E^n\, K^a\, F^n\ ,
\end{eqnarray*}
where ${\CR}_K$ is defined in (\ref{kfjsdv}).
It is also known that this is made into a flip-$*$-Hopf algebra by
\begin{eqnarray*}
(q^{\pm H})^*=q^{\mp H}\ ,\ E^*=-F\ ,\ F^*=-E\ .
\end{eqnarray*} 
Applying the star to ${\CR}$ gives ${\CR}^{*\tens *}=(\CR^{-1})_{21}$ as required to apply the above results. Hence we have a strong quasi-$*$ Hopf algebra with
$\gamma=v^{-1}\nu$. 
\end{example}

\begin{example} \label{hjsdgfcent}
As another example of the structure given in Proposition~\ref{ufudsn},
take the example of $\widetilde{{u_q(su_2)}}$, where $q$ is a primitive $l$th
root of unity (for odd integer $l>1$). This is a central extension of ${u_q(su_2)}$
by an invertible generator $\tilde K$ satisfying the anyonic relations. 
Then we have a new quasitriangular structure:
\begin{eqnarray*}
\tilde {\CR} \,=\,{\CR}_{\tilde K}\,  {\CR}_K\,
\sum_{n=0}^{l-1} \frac{(q-q^{-1})^n}{[n;q^{-2}]!}\, E^n\tens F^n\ ,
\end{eqnarray*}
\end{example}
and the same element $\gamma$ as in 
\ref{hjsdgfc}.

\subsection{Example: Modules over $*$-quasi Hopf algebras and twisting}
We first recall the Drinfeld twist theory. Take ${}_H{\CM}$ to be the category of left modules
of a Hopf algebra $H$.  For a cochain ${\CF}\in H\tens H$ we can twist the coproduct
by $\Delta_{{\CF}}={\CF}\,\Delta(\ ){\CF}^{-1}$ to give a quasiHopf algebra \cite{driquasi} which we denote $H^{\CF}$. For simplicity we assume that $(\id\tens\epsilon){\CF}=
(\epsilon\tens\id){\CF}=1$. This construction can be thought of as twisting the tensor product in ${}_H{\CM}$ to a new category $\CC={}_{H^{\CF}}\CM$ connected by a monoidal `twisting functor' $F:{}_H\CM\to \CC$ \cite{majquasitann,MajBook}. Here $F$ is the identity on objects and
morphisms, $ {\CC}$ has a tensor product,
and the tensor products on the two categories are related by
\begin{eqnarray*}
f_{X,Y}:F(X)\tens F(Y) \to F(X\tens Y)\ ,\quad F(x)\tens F(y)\mapsto 
F({\CF}^{-(1)}\la x \tens {\CF}^{-(2)}\la y)\ .
\end{eqnarray*}
Here since $F$ is the identity on objects we extend its notation to the identity operation on elements, which is a useful convention to allow us to give concise formulae. Thus  $F(x)$ denotes $x\in X$ viewed in $F(X)$. In terms of this notation,
$h\la F(x)=F(h\la x)$. A brief calculation shows that, even though ${}_H{\CM}$
has trivial associator, the new category ${\CC}$ in general
needs an associator $\Phi:(F(X)\tens F(Y))\tens F(Z)\to 
F(X)\tens (F(Y)\tens F(Z))$, given by the action of the following element
$\phi\in H^{\tens 3}$:
\begin{eqnarray*}\label{sjhdcvsjhvg}
\phi\,=\, (\id\tens{\CF})\, ((\id\tens\Delta){\CF})
\,((\Delta\tens\id){\CF}^{-1})\,({\CF}^{-1}\tens\id)\ .
\end{eqnarray*}
The condition that this is the identity is exactly that ${\CF}$ is a cocycle,
but we shall keep to the general case for the moment. 

Now we suppose that $H$ is a $*$-Hopf algebra, and thus that 
${}_H{\CM}$ becomes a bar category with the notations given earlier in \ref{modbar}. 
The problem is how to make ${}_{H^{{\CF}}}{\CM}$ into a bar category. 
We use the following conventions for operations on 
${}_{H^{{\CF}}}{\CM}$:
\begin{eqnarray} \label{kjhdzv}
h\la \overline{F(x)}&=&\overline{ (T_{\CF} h)\la F(x)},\cr
\mathrm{bb}(F(x)) &=&  \overline{\overline{F(\gamma\la x)}}\ , \cr
\mathrm{fb}(\overline{F(x)}) &=& F(\overline{\varphi\la x})\ ,\cr
\Upsilon(\overline{F(x)\tens F(y)}) &=& 
\overline{F({\CG}^{(2)}\la y)} \tens \overline{F({\CG}^{(1)}\la x)}\ ,
\end{eqnarray}
for some $\varphi\in H$ and ${\CG}\in H\tens H$, and some algebra homomorphism $T_{\CF}$ that determines the module structure on barred objects in the target category just as $*S$ does in ${}_H\CM$. 

\begin{theorem}\label{mdvfmasm}
Suppose that $H$ is a $*$-Hopf algebra, and ${\CF}\in H\tens H$ is 
an invertible element with $(\epsilon\tens\id){\CF}=
(\id\tens\epsilon){\CF}=1$. Then ${}_{H^{{\CF}}}{\CM}$ is
 a strong bar category such that the twisting functor from ${}_{H}{\CM}$ to
 ${}_{H^{{\CF}}}{\CM}$  is
a morphism of bar categories if and only if there
exists invertible $\varphi\in H$ such that $\epsilon(\varphi)=1$. In this case
\begin{equation*}T_{\CF}(h)=\varphi^{-1}(S h)^*\varphi,\quad 
\gamma\,=\, \varphi^{-1}\,(S\varphi^{-1})^*\ ,\quad 
{\CG}\,=\,  (\varphi^{-1}\tens\varphi^{-1})(*S\tens *S)(\CF_{21})(\Delta\varphi)\CF^{-1}.\end{equation*}
\end{theorem}
\begin{proof} 
We calculate $\mathrm{bb}$ from part (3) of Definition~\ref{recdef}) as
\begin{eqnarray*}
\mathrm{bb}_{F(X)} (F(x))
 &=&  \overline{\mathrm{fb}_X^{-1}}\circ \mathrm{fb}_{\bar X}^{-1}\circ F(\mathrm{bb}_X)
 (F(x)) \cr
 &=&  \overline{\mathrm{fb}_X^{-1}}\circ \mathrm{fb}_{\bar X}^{-1}
 (F(\overline{\overline{x}}))\cr
  &=&  \overline{\mathrm{fb}_X^{-1}}
(\overline{ F(\varphi^{-1}\la \overline{x})}) \cr
&=&  \overline{\mathrm{fb}_X^{-1}}
(\overline{ F(\overline{(S\varphi^{-1})^*\la x})})  \cr
&=& 
\overline{ \overline{F(\varphi^{-1}\,(S\varphi^{-1})^*\la x)}}.
\end{eqnarray*}
Similarly, we determine $\Upsilon$ from part (4) of Definition~\ref{recdef} as
\begin{eqnarray*}
\Upsilon(\overline{F(x)\tens F(y)}) &=& 
(\mathrm{fb}^{-1}\tens \mathrm{fb}^{-1})\, f^{-1}\,F(\Upsilon)\,\mathrm{fb}\,
\overline{f}\,(\overline{F(x)\tens F(y)}) \cr
 &=& 
(\mathrm{fb}^{-1}\tens \mathrm{fb}^{-1})\, f^{-1}\,F(\Upsilon)\,\mathrm{fb}\,
(\overline{F({\CF}^{-(1)}\la x \tens {\CF}^{-(2)}\la y)}) \cr
 &=& 
(\mathrm{fb}^{-1}\tens \mathrm{fb}^{-1})\, f^{-1}\,F(\Upsilon)\,
(F(\overline{\varphi_{(1)}\,{\CF}^{-(1)}\la x \tens
\varphi_{(2)}\, {\CF}^{-(2)}\la y})) \cr
 &=& 
(\mathrm{fb}^{-1}\tens \mathrm{fb}^{-1})\, f^{-1}\,
(F(\overline{\varphi_{(2)}\, {\CF}^{-(2)}\la y} \tens
\overline{\varphi_{(1)}\,{\CF}^{-(1)}\la x })) \cr
 &=& 
(\mathrm{fb}^{-1}\tens \mathrm{fb}^{-1})\, 
(F({\CF}^{(1)}\la\overline{\varphi_{(2)}\, {\CF}^{-(2)}\la y}) \tens
F({\CF}^{(2)}\la\overline{\varphi_{(1)}\,{\CF}^{-(1)}\la x })) \cr
 &=& 
\overline{F(\varphi^{-1}\,
(S{\CF}^{(1)})^*\,\varphi_{(2)}\, {\CF}^{-(2)}\la y)} \tens
\overline{F( \varphi^{-1}\,
(S{\CF}^{(2)})^*\,\varphi_{(1)}\,{\CF}^{-(1)}\la x )} \ .
\end{eqnarray*}
Finally, that $\mathrm{fb}_X$ is a morphism $\overline{F(X)}\to F(\overline{X})$ determines $T_{\CF}$. We can then use $\mathrm{fb}$ to transfer over the strong bar category structure of ${}_H\CM$ to one of ${}_{H^{\CF}}\CM$. It is instructive, although not necessary, to verify the strong bar category structure explicitly. For example,  \ref{barcatdef} (d) requires that $\gamma\,=\, T_{\CF}(\gamma)$ and one may check that this always holds as $T_{\CF}(\gamma)=\varphi^{-1} (S\gamma)^*\varphi=\varphi^{-1} (S(\varphi^{-1}(S\varphi^{-1})^*))^*\varphi=\gamma$. Similarly, the condition in \ref{barcatdef} with $\Upsilon$ and the associators gives
the equality of the following expressions:
\begin{eqnarray*}
&& \kern-30pt \Phi\,(\Upsilon\tens\id)\,\Upsilon\,(\overline{F(x)\tens( F(y)\tens F(z)) }) \cr
&=& \Phi\,(\Upsilon\tens\id)\,(
\overline{ 
F({\CF}^{(1)}\,{\CG}^{(2)}_{\phantom{(1)}(1)}{\CF}^{-(1)} \la y)\tens 
F({\CF}^{(2)}\,{\CG}^{(2)}_{\phantom{(1)}(2)}{\CF}^{-(2)} \la z) }
\tens \overline{F({\CG}^{(1)}\la x) }) \cr
&=& \Phi\,((
\overline{ 
F(\tilde {\CG}^{(2)}\,
{\CF}^{(2)}\,{\CG}^{(2)}_{\phantom{(1)}(2)}{\CF}^{-(2)} \la z) } \tens
\overline{ 
F(\tilde {\CG}^{(1)}\,{\CF}^{(1)}\,{\CG}^{(2)}_{\phantom{(1)}(1)}
{\CF}^{-(1)} \la y)}
) \tens \overline{F({\CG}^{(1)}\la x) }
) \ ,\cr
&&\kern -30pt (\id\tens\Upsilon)\,\Upsilon\,\overline{\Phi^{-1}}\,
(\overline{F(x)\tens( F(y)\tens F(z)) }) 
 \cr 
 &=& (\id\tens\Upsilon)\,\Upsilon\,
 (\overline{(F(\phi^{-(1)}\la x)\tens F(\phi^{-(2)}\la y))\tens F(\phi^{-(3)}\la z) }) \cr
  &=& (\id\tens\Upsilon)\,(\overline{F({\CG}^{(2)}\,\phi^{-(3)}\la z)}\tens 
  \overline{
  F({\CF}^{(1)}\,{\CG}^{(1)}_{\phantom{(1)}(1)}{\CF}^{-(1)} 
  \phi^{-(1)}\la x)\tens F(
 {\CF}^{(2)}\,{\CG}^{(1)}_{\phantom{(1)}(2)}{\CF}^{-(2)}\phi^{-(2)}\la y)
  }) \cr
  &=& \overline{F({\CG}^{(2)}\,\phi^{-(3)}\la z)}\tens (\overline{F(
 \tilde {\CG}^{(2)}\,{\CF}^{(2)}\, {\CG}^{(1)}_{\phantom{(1)}(2)}{\CF}^{-(2)}
 \phi^{-(2)}\la y)} \tens\overline{
   F( \tilde {\CG}^{(1)}\,{\CF}^{(1)}\,
   {\CG}^{(1)}_{\phantom{(1)}(1)}{\CF}^{-(1)} 
  \phi^{-(1)}\la x)
} )\ ,
\end{eqnarray*}
which amounts to the identity
\begin{eqnarray}\label{shjvccn}
&& \kern-30pt T_{\CF}(\phi^{(3)})\,{\CG}^{(1)} 
\tens 
T_{\CF}(\phi^{(2)})\,\tilde {\CG}^{(1)}\,{\CF}^{(1)}\,{\CG}^{(2)}_{\phantom{(1)}(1)}
{\CF}^{-(1)} 
 \tens 
 T_{\CF}(\phi^{(1)})\,\tilde {\CG}^{(2)}\,
{\CF}^{(2)}\,{\CG}^{(2)}_{\phantom{(1)}(2)}{\CF}^{-(2)} \nonumber\\ 
&=& 
  \tilde {\CG}^{(1)}\,{\CF}^{(1)}\,
   {\CG}^{(1)}_{\phantom{(1)}(1)}{\CF}^{-(1)} 
  \phi^{-(1)}
\tens 
\tilde {\CG}^{(2)}\,{\CF}^{(2)}\, {\CG}^{(1)}_{\phantom{(1)}(2)}{\CF}^{-(2)}
 \phi^{-(2)} 
 \tens {\CG}^{(2)}\,\phi^{-(3)}  \ . \nonumber\\ 
\end{eqnarray}
and one may check that this is always verified for the $\CG, T_{\CF}$ stated. Similarly, an explicit  formula for $\Xi$ is deduced from
\begin{eqnarray*}
\Upsilon\, \overline{\Upsilon}\,\mathrm{bb}\,(F(x)\tens F(y)) &=&
\Upsilon\, \overline{\Upsilon}\,(\overline{\overline{
F({\CF}^{(1)}\,\gamma_{(1)}\,{\CF}^{-(1)}\,\la x) \tens 
F({\CF}^{(2)}\,\gamma_{(2)}\,{\CF}^{-(2)}\,\la y)  }}) \cr
&=& \Upsilon\, (\overline{
\overline{ F({\CG}^{(2)}\,{\CF}^{(2)}\,\gamma_{(2)}\,{\CF}^{-(2)}\,\la y)  }
\tens
\overline{ F({\CG}^{(1)}\,{\CF}^{(1)}\,\gamma_{(1)}\,{\CF}^{-(1)}\,\la x)  }
}) \cr
&=& \overline{\overline{ 
F((T_{\CF} \tilde{\CG}^{(2)})\,{\CG}^{(1)}\,{\CF}^{(1)}\,
\gamma_{(1)}\,{\CF}^{-(1)}\,\la x)  }} \cr && \tens
\overline{\overline{ 
F(T_{\CF}(\tilde{\CG}^{(1)})\,{\CG}^{(2)}\,{\CF}^{(2)}\,
\gamma_{(2)}\,{\CF}^{-(2)}\,\la y)  }}\ .
\end{eqnarray*}
as the action $\Xi(x\tens y)=\Xi^{(1)}\la x\tens \Xi^{(2)}\la y$ of the element
\begin{equation}\label{Xiformula}\Xi= (\gamma^{-1}\tens\gamma^{-1})((T_{\CF}\tens T_{\CF})(\CG_{21}))\CG\CF(\Delta\gamma)\CF^{-1}.\end{equation}
One may verify that this is the identity for the particular $\CG,\gamma, T_{\CF}$, as it must by Proposition~\ref{Xipreserved}.
\end{proof}

\begin{remark}\label{twicoc} Note that there is nothing stopping us taking $\varphi=1$ in this theorem  as well as in the following corollaries, to define a canonical twisting induced by $F$ alone with $\gamma=1$ and $\CG=(*S\tens *S)(\CF_{21})\CF^{-1}$. Note also that $(*S\tens *S)(\CF)=\CF$ is a reality condition sometimes verified (e.g. for the Octonions category) in which case $\CG=\CF_{21}\CF^{-1}$.   \end{remark}

\begin{corol} If $H$ is a Hopf $*$ algebra and $\CF$ in the above theorem is moreover a 2-cocycle (i.e.\ $\phi=1\tens 1\tens 1$ in (\ref{sjhdcvsjhvg})) then the twisted Hopf algebra $H^{\CF}$ becomes a strong quasi-$*$ Hopf algebra for any invertible $\varphi$, with $\gamma,\CG$ as stated. \end{corol}
\begin{proof} This is  a special case  of the theorem where $H^{\CF}$ remains a Hopf algebra. One may verify that this indeed yields a strong quasi-$*$ Hopf algebra as defined in Definition~\ref{qudef}, as it must by general Tannaka-Krein considerations (see appendix).  It is only necessary to factor $T_{\CF}=*_{\CF}S_{\CF}$ to obtain an explicit formula for $*_{\CF}$ from the known one for $S_{\CF}$ in \cite{MajBook}. \end{proof}

\begin{cor} \label{kgkgk}
If $\CF$ is a cocycle then $H^{\CF}$ is a $*$-Hopf algebra with 
the twisting functor ${}_H\CM\to {}_{H^{\CF}}\CM$ a bar functor 
if and only if there is an invertible element 
$\varphi\in H$ such that 
\[ \epsilon(\varphi)=1,\quad (S\varphi)^*=\varphi^{-1},\quad (*S\otimes 
*S)(\CF_{21})=(\varphi\otimes\varphi)\CF(\Delta\varphi^{-1})\] 
Here 
\[ *_{\CF}(h)=V (h^*) V^{-1};\quad V=\varphi (\CF{\uo} S \CF\ut)^*\] 
\end{cor}
\begin{proof} This is the special case $\gamma=1$ and $\CG=1\tens 1$, since 
these are trivial for a usual $*$-Hopf algebra. The reality condition when $\varphi=1$ 
recovers the condition $(*S\tens *S)(\CF)=\CF_{21}$ known for this case in \cite[Prop~2.3.7]{MajBook}. \end{proof}

\medskip

Now we look at this from a different point of view. Suppose that we are just
given a quasi-Hopf algebra, without a description of how it was obtained by twisting
(i.e.\ we do not know ${\CF}$ or $\varphi$ in the discussion above). 
However we do know $\Delta_{\CF}$, as it is just the coproduct
in the quasi-Hopf algebra. We describe a $*$-quasi Hopf algebra 
via an antilinear algebra map $T:H\to H$, which for 
a usual $*$-Hopf algebra is just $T(h)=(Sh)^*$ and which for the examples
constructed above by twisting is the map $T_{\CF}$,

\begin{defin}\label{kjzshcnb}
A $*$-quasi Hopf algebra $(H,\Delta,\phi,T,\gamma,{\CG})$ is a quasi Hopf algebra
$(H,\Delta,\phi)$ with additional data
a conjugate-linear algebra map $T:H\to H$, invertible $\gamma\in H$ and 
invertible ${\CG}\in H\tens H$ satisfying
\begin{eqnarray*}
T^2(h) &=& \gamma h \gamma^{-1}\ ,\cr
\tau\Delta(T h)&=&\CG^{-1}((T\tens T)\,\tau\Delta h)\CG\ ,\cr
(\epsilon\tens\id){\CG}\,=\,(\id\tens\epsilon){\CG} &=& 1\ ,\cr
\gamma &=& T(\gamma)\ , \cr
 ((T\tens T\tens T)\phi_{321})\, (1\tens {\CG})\,
((\id\tens\Delta){\CG}) &=& ({\CG}\tens 1)\,
((\Delta \tens\id){\CG})\, \phi^{-1}\ .
\end{eqnarray*}
The $*$-quasiHopf algebra is called {\em strong} if $(\gamma\tens\gamma)\Delta\gamma^{-1}=((T\tens T)(\CG_{21}))\CG$.
\end{defin}
The fifth condition here, for example, is extrapolated from (\ref{shjvccn}) but now taken to apply more generally than these twisting examples.

\begin{propos}\label{fgdfgdgyy}
If (resp. strong) $*$-quasi Hopf algebra then ${}_H{\CM}$ is a (resp. strong) bar category. The associator in the category is
\begin{eqnarray*}
\Phi((u\tens v)\tens w) \,=\, \phi^{(1)}\la u \tens (\phi^{(2)}\la v \tens \phi^{(3)}\la w)\ ,
\end{eqnarray*}
and the operations for the bar category are given by 
\[ h\la\overline{v}=\overline{T(h)\la v}, \quad \mathrm{bb}(x)=\overline{\overline{\gamma\la x}},\quad 
\Upsilon(\overline{x\tens y}) \,=\, \overline{{\CG}^{(2)}\la y} \tens 
\overline{{\CG}^{(1)}\la x}\ . \]
The coboundary natural isomorphism $\Xi$ is given by
\begin{eqnarray*}
\Xi(x\tens y) \,=\, \gamma^{-1}T(\tilde{{\CG}}^{(2)})\,{\CG}^{(1)}\gamma_{(1)}\la x \tens
\gamma^{-1}T(\tilde{{\CG}}^{(1)})\,{\CG}^{(2)}\gamma_{(2)}\la y\ .
\end{eqnarray*}
\end{propos}
\begin{proof} The condition for $\mathrm{bb}$ and
$\Upsilon$ to be morphisms are the first two
equations in the list in \ref{kjzshcnb}. The third equation gives
 condition (a) in \ref{barcatdef}, and the fourth equation gives
 condition (d). For condition (b), this is direct calculation as in the derivation of (\ref{shjvccn}). Similarly,  direct calculation gives the form of $\Xi$ as in the derivation (\ref{Xiformula}). \end{proof}

\begin{remark} Clearly one can again factor out the antipode $S$ of a quasi-Hopf algebra
and look at the formal properties of  $*=T\circ S^{-1}$ in line with usual $*$-Hopf algebras but
this is not particularly convenient in the quasi-Hopf algebra case.
\end{remark}

One may verify that $H^{\CF}$ in Theorem~\ref{mdvfmasm} is indeed a $*$-quasiHopf algebra and the special case of a $*$-quasiHopf algebra with $\phi=1$ is a quasi-$*$ Hopf algebra as previously defined. We will later give a concrete example of this theory applied to the Octonions.

\subsection{Example: The coset representative category}

Suppose that $G$ is a subgroup of a finite group $X$. 
The set $G\backslash X$ consists of the equivalence classes
of $X$ under left translation by elements of $G$, i.e.\ the left cosets.
The subgroup $G$ acts on $G\backslash X$ by right translation, which
we write as $\ra:G\backslash X\times G\to G\backslash X$
(i.e.\ $p\ra u\in G\backslash X$ for $p\in G\backslash X$ and $u\in G$).

\begin{defin}
The category ${\CC}_{\bowtie}$ consists of all finite dimensional vector spaces
over $\mathbb{C}$, whose objects are right representations of the group 
$G$ and possess $G\backslash X$-gradings, 
i.e.\ an object $V$ decomposes as a direct sum of subspaces
 $V=\oplus_{p\in G\backslash X} V_p$, with a compatibility
condition between the action and the grading. If $\xi\in V_p$ for some
$p\in G\backslash X$ we say that $\xi$ is a  homogenous element
of $V$, with grade $\langle \xi\rangle=p$. 

 We write the action for the representation as $\bar\ra:V\times G\to V$. 
The compatibility condition is
$\langle \xi\bar\ra u\rangle=\langle \xi\rangle\ra u$. 

 The morphisms are linear maps which preserve both the grading and the action, 
i.e.\ for a morphism $\phi:V\to W$ we have $\langle
\phi(\xi)\rangle =\langle \xi\rangle $ and $\phi(\xi)\bar\ra
u=\theta(\xi\bar\ra u)$ for all homogenous $\xi\in V$ and $u\in G$.
\end{defin}

The category has not yet been given a tensor product operation,
and to do this we need to make a choice of coset representatives. This is
a subset 
 $M\subset X$ such that for every $x\in X$ there is a unique $s\in M$ so that
$x\in Gs$. We shall call the decomposition $x=us$ for $u\in G$ and $s\in M$ the 
unique factorisation of $x$. For convenience, 
we identify the coset representatives with the cosets.
 The identity in $X$ will
be denoted $e$, and we assume that $e\in M$.

\begin{defin} Given $s,t\in M$, define $\tau(s,t)\in G$ and $s\cdot t\in M$
by the unique factorisation $st=\tau(s,t)(s\cdot t)$ in $X$. We also define functions
$\la:M\times G\to G$ and $\ra:M\times G\to M$ by the unique factorisation
$su=(s\la u)(s\ra u)$ for $s,s\ra u\in M$ and $u,s\la u\in G$. 
\end{defin}

The assumption that $e\in M$ gives the binary operation $(M,\cdot)$ a 2-sided identity.
Also there is a unique left inverse $t^L$ for every $t\in M$, satisfying the equation
$t^L\cdot t=e$. We will assume that there is also a unique right inverse $t^R$,
satisfying $t\cdot t^R=e$. For the identities relating the binary operation
and the `cochain' $\tau$, see \cite{cosrep}.

\begin{propos} Given a set of coset representatives $M$,
we can make ${\CC}_{\bowtie}$
 into a tensor category by taking $V\tens W$ to be the usual vector space
tensor product, with actions and gradings given by
\[
\langle \xi\tens \eta\rangle =\langle \xi\rangle \cdot\langle  \eta\rangle \quad{ \mathrm{and}}\quad
(\xi\otimes\eta)\bar\ra u=\xi\bar\ra(\langle \eta\rangle \la u)\otimes \eta\bar\ra u\ .
\]
For morphisms $\theta:V\to \tilde V$ and $\phi:W\to \tilde W$ we define
the morphism $\theta\tens\phi:V\tens W\to \tilde V\tens \tilde W$ by $(\theta\tens\phi)(\xi\tens\eta)=
\theta(\xi)\tens\phi(\eta)$, which is just the usual vector space formula.

The identity for the tensor operation is just the
vector space $\mathbb{C}$ with trivial $G$-action and grade $e\in M$.  For any
object $V$ the morphisms $l_V:V\to V\tens \mathbb{C}$ and $r_V:V\to \mathbb{C}\tens V$
are given by the formulae $l_V(\xi)=\xi\tens 1$ and $r_V(\xi)=1\tens
\xi$, where $1$ is the multiplicative identity in ${\CK}$.

This tensor product is not consistent
with the usual idea of ignoring the brackets in
tensor products.  The non-trivial associator $\Phi_{UVW}:(U\otimes V)\otimes
W\to U\otimes (V\otimes W)$ is given by
$$
\Phi((\xi\otimes\eta)\otimes\zeta)\ =\ \xi\bar\ra\tau(\langle \eta\rangle ,\langle \zeta\rangle )
\otimes(\eta\otimes\zeta)\ .
$$
\end{propos} 

We can make ${\CC}_{\bowtie}$ into a bar category as follows. Again we
identify objects $V$ and $\overline{V}$ as sets. The grading and action on 
 $\overline V$ are given by
\begin{eqnarray*}
\<\overline{v}\> \,=\, \<v\>^R\ ,\quad \overline{v}\bar\ra u \,=\, 
\overline{v\bar\ra (\<v\>^R\la u)}\ .
\end{eqnarray*}
It is left to the reader to show that this formula gives a right $G$-action, and that
the compatibility condition is satisfied. Then the morphism
$\mathrm{bb}_V:V\to \overline{\overline{V}}$ is given by
\begin{eqnarray*}
\mathrm{bb}_V(v) \,=\, \overline{\overline{v\bar\ra\tau(\<v\>^L,\<v\>)^{-1}}}\ .
\end{eqnarray*}
We also have the following formula for $\Upsilon^{-1}_{WV}:
\overline{V}\tens \overline{W} \to \overline{W\tens V}$,
\begin{eqnarray*}
\Upsilon^{-1}_{WV}(\overline{v}\tens \overline{w}) &=&
\overline{ w\bar\ra \tau(\<v\>\ra \tau(\<v\>^R,\<w\>^R), \<v\>^R\!\cdot\!
\<w\>^R )^{-1}\tens 
v\ra \tau(\<v\>^R,\<w\>^R)}\cr
&=& \overline{ (w\bar\ra \tau(\<v\>,\<v\>^R)^{-1}\tens v)\ra \tau(\<v\>^R,\<w\>^R)
}\ .
\end{eqnarray*}
With some effort it can be shown that $\Xi$ is the identity, so the bar category is strong. 

\section{Braided bar categories}

We recall that monoidal category is braided if there is a natural equivalence $\Psi:\tens\to \tens^{\rm op}$ obeying two natural `hexagon coherence identities' with respect to $\tens$ of each factor. These are sufficient to ensure the Yang-Baxter or representation of the braid relations between any three objects. There is also compatibility with $l,r$. We refer to \cite{Mac,JS} for details and \cite{MajBook} for our notations. 

\subsection{Real and hermitian braidings}

If the bar category ${\CC}$ is also braided, we can ask how
the braiding $\Psi$ behaves under the bar functor. In particular we can 
highlight two cases:

\begin{defin} \label{psdhvsd1}
If the braiding in a bar category ${\CC}$ fits the following diagram:

\begin{picture}(100,80)(-100,6)

\put(60,63){\vector(1,0){45}}
\put(60,23){\vector(1,0){45}}
\put(35,53){\vector(0,-1){20}}
\put(127,53){\vector(0,-1){20}}
\put(20,60){$\overline{X}\tens\overline{Y}$}
\put(20,20){$\overline{Y\tens X}$}
\put(111,60){$\overline{Y}\tens\overline{X}$}
\put(111,20){$\overline{X\tens Y}$}

\put(70,70){$\Psi$}
\put(67,30){$\overline{\Psi^{\pm1}}$}
\put(12,40){$\Upsilon^{-1}$}
\put(133,40){$\Upsilon^{-1}$}
\end{picture}

\noindent
and the sign is $+1$, we call
the braiding real. If the sign is $-1$, we call
the braiding antireal. 
\end{defin}

Clearly the two notions coincide in the symmetric case where $\Psi^2=\id$. In this case we say that the bar category is symmetric. Clearly {\bf Vect} is an elementary example with $\Psi$ the usual flip on tensor products.

\begin{example}
Consider modules ${}_H{\CM}$ over a $*$-Hopf algebra $H$ (see \ref{modbar}),
where $H$ is also quasitriangular, with ${\CR}\in H\tens H$. 
If we write ${\CR}={\CR}^{(1)}\tens {\CR}^{(2)}$
(summation implicit), then the braiding in ${}_H{\CM}$ is defined as
\begin{eqnarray*}
\Psi_{V,W}(v\tens w)\,=\, {\CR}^{(2)}\la w \tens {\CR}^{(1)}\la v\ .
\end{eqnarray*}
In the literature there are two cases for 
${\CR}^{(1)*}\tens {\CR}^{(2)*}$ singled out, and
we will consider the effect on the braiding of both of these possibilities. The reader should
 remember that $(S\tens S){\CR}={\CR}$. 

\medskip \noindent Case 1 - the real case - the braiding is real: 
${\CR}^{(1)*}\tens {\CR}^{(2)*}={\CR}^{(2)}\tens {\CR}^{(1)}$. 
\begin{eqnarray*}
\Psi_{\overline{V},\overline{W}}(\overline{v}\tens \overline{w}) &=& 
\overline{{\CR}^{(2)*}\la w} \tens \overline{{\CR}^{(1)*}\la v } \cr
&=& 
\overline{{\CR}^{(1)}\la w} \tens \overline{{\CR}^{(2)}\la v }\ .
\end{eqnarray*}

\medskip \noindent Case 2 - the antireal case - the braiding is antireal: 
${\CR}^{(1)*}\tens {\CR}^{(2)*}={\CR}^{-1}$. 
\begin{eqnarray*}
\Psi_{\overline{V},\overline{W}}(\overline{v}\tens \overline{w}) &=& 
\overline{{\CR}^{(2)*}\la w} \tens \overline{{\CR}^{(1)*}\la v } \cr
&=& 
\overline{({\CR}^{-1})^{(2)}\la w} \tens \overline{({\CR}^{-1})^{(1)}\la v }\ .
\end{eqnarray*}
\end{example}

\begin{example}
 For Yetter Drinfeld modules over $*$-Hopf algebras (see \ref{YDbar}), 
 there is a braiding given by
 \begin{eqnarray*}
\Psi(v\tens w) \,=\, w_{[0]}\tens v\ra w_{[1]}\ .
\end{eqnarray*}
This is antireal, as can be seen by the following calculation:
 \begin{eqnarray*}
\Psi(\overline{v}\tens \overline{w}) &=& \overline{w}_{[0]}\tens 
\overline{v}\ra \overline{w}_{[1]} \cr
&=& \overline{w_{[0]}}\tens \overline{v}\ra w_{[1]}^*   \cr
&=& \overline{w_{[0]}}\tens \overline{v\ra S^{-1}(w_{[1]}) }\ , \cr
\Upsilon^{-1}\Psi(\overline{v}\tens \overline{w}) &=& \overline{
v\ra S^{-1}(w_{[1]}) \tens w_{[0]} } \cr
&=& \overline{\Psi^{-1}(w\tens v)} \cr
&=& \overline{\Psi^{-1}}\Upsilon^{-1}(\overline{v}\tens \overline{w})
\end{eqnarray*}

\end{example}

\begin{example} 
For the example of ${u_q(su_2)}$, where $q$ is a primitive $l$th
root of unity (for odd integer $l>1$) (see \ref{hjsdgfc}) the braiding
defined by 
\begin{eqnarray*}
\Psi_{V,W}(v\tens w)\,=\, {\CR}^{(2)}\la w \tens {\CR}^{(1)}\la v
\end{eqnarray*}
is antireal. Note that $\Upsilon$ is non-trivial in this case.
\end{example}

\subsection{Real forms of bar categories}

\begin{defin} \label{bardefinhh}
A star bar category is a bar category ${\CC}$ 
such that

\noindent a)\quad There is
 a natural transformation
$\star$ from the identity to the $\mathrm{bar}$ functor, making every object into a
star object (see \ref{jkhdsz}). (For the identity object this is the same
as the $\star$ in the definition of bar category.)

\noindent b)\quad  The following diagram commutes:

\begin{picture}(100,80)(-100,6)

\put(60,63){\vector(1,0){45}}
\put(60,23){\vector(1,0){45}}
\put(28,53){\vector(0,-1){20}}
\put(137,53){\vector(0,-1){20}}
\put(-7,60){$(X\tens Y)\tens Z$}
\put(-7,20){$X\tens(Y\tens Z)$}
\put(111,60){$\overline{(X\tens Y)\tens Z}$}
\put(111,20){$\overline{X\tens(Y\tens Z)}$}

\put(57,70){$\star_{(X\tens Y)\tens Z}$}
\put(57,30){$\star_{X\tens (Y\tens Z)}$}
\put(-6,40){$\Phi_{X,Y,Z}$}
\put(143,40){$\overline{\Phi_{X,Y,Z}}$}
\end{picture}

\end{defin}

 Let us pause to unpack this definition. 
Firstly, the choice of the morphisms $\star$ is natural in the sense that,
 for any morphism $\phi:X\to Y$ we have $\overline{\phi}\circ \star_X=\star_Y\circ \phi:
 X\to \overline{Y}$. Secondly, we have two obvious morphisms from $X$ to 
 $\overline{\overline{X}}$, $\mathrm{bb}_X$ and $\star_{\bar X} \star_X$, and we insist that
 these are the same. In the example of the bimodule category, 
 we could write  $\star_E:E\to \overline{E}$ as
 $e\mapsto \overline{e^*}$, and we would then have $e^{**}=e$, 
 which is a rather more familiar way of writing an involution. 
There are some immediate deductions which we can make from the definition. 

\begin{propos} \label{skjd}
 $\star_{\bar X}=\overline{\star_X}:\overline{X}\to 
\overline{\overline{X}}$. 
\end{propos}
\begin{proof} Begin with the morphism $\star_X:X\to\overline{X}$, 
and use the fact that $\star$ is a natural transformation between
the identity functor and bar to see that $\star_{\bar X}\,\star_X=\overline{\star_X}\,\star_X$. 
Since $\star_X$ is invertible, we get the result.\end{proof}

\medskip It will be convenient to make the following definition:

\begin{defin}\label{skjd2}
 There is an invertible natural transformation $\Gamma$ between
the functors $\tens$ and $\tens\circ\mathrm{flip}$ from ${\CC}
\times {\CC}\to {\CC}$ so that the following diagram commutes:

\begin{picture}(100,80)(-100,6)

\put(60,63){\vector(1,0){45}}
\put(60,23){\vector(1,0){45}}
\put(28,33){\vector(0,1){20}}
\put(127,53){\vector(0,-1){20}}
\put(17,60){$X\tens Y$}
\put(17,20){$Y\tens X$}
\put(111,60){$\overline{X}\tens\overline{Y}$}
\put(111,20){$\overline{Y\tens X}$}

\put(68,70){$\star\tens\star$}
\put(76,30){$\star$}
\put(-3,40){$\Gamma_{Y,X}$}
\put(133,40){$\Upsilon^{-1}$}
\end{picture}
\end{defin}

\begin{propos}\label{pty111} In a star bar category ${\CC}$,
 $\overline{\Gamma}\, \Upsilon^{-1}=
\Upsilon^{-1}\Xi\,\Gamma^{-1}:\overline{X} \tens \overline{Y}\to
\overline{X\tens Y}$. 
\end{propos}
\begin{proof} We can rewrite the equation
$\overline{\star_{X\tens Y}}\,\star_{X\tens Y}=\mathrm{bb}_{X\tens Y}$,
using the definition of $\Gamma$ on one side and the definition
of bar category on the other, to give
\begin{eqnarray*}
\overline{\Upsilon^{-1}(\star\tens\star)\,\Gamma}\, \Upsilon^{-1}(\star\tens\star)\,\Gamma
&=& \overline{\Upsilon^{-1}}\, \Upsilon^{-1} \, 
(\overline{\star}\tens\overline{\star})\, (\star\tens\star)\,\Xi\ .
\end{eqnarray*}
 From the naturality of $\Gamma$, $(\star\tens\star)\,\Gamma=
\Gamma\,(\star\tens\star)$, and similarly for $\Xi$, so 
\begin{eqnarray*}
\overline{\Upsilon^{-1}(\star\tens\star)\,\Gamma}\, \Upsilon^{-1}\Gamma\,(\star\tens\star)
&=& \overline{\Upsilon^{-1}}\, \Upsilon^{-1} \, 
(\overline{\star}\tens\overline{\star})\,\Xi\, (\star\tens\star)\ ,
\end{eqnarray*}
and some cancellation gives
\begin{eqnarray*}
\overline{(\star\tens\star)\,\Gamma}\, \Upsilon^{-1}\Gamma
&=& \Upsilon^{-1} \, 
(\overline{\star}\tens\overline{\star})\,\Xi\ .
\end{eqnarray*}
{}From the naturality of $\Upsilon$, $\Upsilon^{-1} \, 
(\overline{\star}\tens\overline{\star})=\overline{(\star\tens\star)}\,
\Upsilon^{-1}$, and using this we get the result.\end{proof}

\begin{propos}\label{pty1117} In a star bar category ${\CC}$,
$\Xi=\Gamma^2$. 
\end{propos}
\begin{proof} As $\star$ is a natural transformation, we have
\begin{eqnarray*}
\star=\overline{\Gamma}\,\star\,\Gamma^{-1}:X\tens Y\to \overline{X\tens Y}\ .
\end{eqnarray*}
From the definition of $\Gamma$ in \ref{skjd2},
\begin{eqnarray*}
\star=\overline{\Gamma}\,\Upsilon^{-1}\,(\star\tens\star):X\tens Y\to \overline{X\tens Y}\ ,
\end{eqnarray*}
and using \ref{pty111} this gives
\begin{eqnarray*}
\star=\Upsilon^{-1}\Xi\,\Gamma^{-1}\,(\star\tens\star):X\tens Y\to \overline{X\tens Y}\ .
\end{eqnarray*}
Using the definition of $\Gamma$ again gives
\begin{eqnarray*}
\Upsilon^{-1}(\star\tens\star)\Gamma=\Upsilon^{-1}\Xi\,\Gamma^{-1}\,(\star\tens\star):X\tens Y\to \overline{X\tens Y}\ ,
\end{eqnarray*}
and the naturality of $\Gamma$ gives the answer.\end{proof}

\begin{propos}\label{pty222}  Given the rest of the definition,
the following condition is equivalent to \ref{bardefinhh} part (b):
\begin{eqnarray*}
\Phi\,(\Gamma_{Y,Z}\tens\id)\,\Gamma_{X,Y\tens Z}\,=\,
(\id\tens\Gamma_{X,Y})\,\Gamma_{X\tens Y, Z}\,\Phi^{-1}: X\tens (Y \tens Z) \to
Z\tens (Y \tens X)\ . 
\end{eqnarray*}
\end{propos}
\begin{proof} First consider $\star_{(X\tens Y)\tens Z}$. From the definition of $\Gamma$ we can rewrite this as
\begin{eqnarray*}
\star_{(X\tens Y)\tens Z} &=& 
\Upsilon^{-1}_{X\tens Y,Z}\,(\star_{Z}\tens\star_{X\tens Y})\,\Gamma_{X\tens Y,Z}\ ,
\end{eqnarray*}
and applying this again,
\begin{eqnarray*}
\star_{(X\tens Y)\tens Z} &=& 
\Upsilon^{-1}_{X\tens Y,Z}\,(\id\tens\Upsilon^{-1}_{X,Y})\,
(\star_{Z}\tens(\star_{Y}\tens \star_{X}))\,
(\id\tens\Gamma_{X,Y})\,\Gamma_{X\tens Y,Z}\ .
\end{eqnarray*}
Now we rewrite $\star_{X\tens(Y\tens Z)}$ using $\Gamma$, as
\begin{eqnarray*}
\star_{X\tens(Y\tens Z)} &=& 
\Upsilon^{-1}_{X,Y\tens Z}\,(\Upsilon^{-1}_{Y,Z}\tens\id)\,
((\star_{Z}\tens\star_{Y})\tens \star_{X})\,(\Gamma_{Y,Z}\tens\id)\,\Gamma_{X,Y\tens Z}
\end{eqnarray*}
{}From the definition of bar category, we have
\begin{eqnarray*}
\overline{\Phi_{X,Y,Z}}\,\star_{(X\tens Y)\tens Z} &=& 
\Upsilon^{-1}_{X,Y\tens Z}\,(\Upsilon^{-1}_{Y,Z}\tens\id)\,
\Phi^{-1}_{\bar Z,\bar Y,\bar X}\, 
(\star_{Z}\tens(\star_{Y}\tens \star_{X}))\,
(\id\tens\Gamma_{X,Y})\,\Gamma_{X\tens Y,Z}\ ,
\end{eqnarray*}
and by naturality of $\Phi$,
\begin{eqnarray*}
\overline{\Phi_{X,Y,Z}}\,\star_{(X\tens Y)\tens Z} &=& 
\Upsilon^{-1}_{X,Y\tens Z}\,(\Upsilon^{-1}_{Y,Z}\tens\id)\,
((\star_{Z}\tens\star_{Y})\tens \star_{X})\,\Phi^{-1}_{Z,Y,X}\, 
(\id\tens\Gamma_{X,Y})\,\Gamma_{X\tens Y,Z}\ . 
\end{eqnarray*}
If we substitute this into the condition in \ref{bardefinhh} part (b),
we can cancel terms to give the equivalence. \end{proof}

From the categorical point of view, $\Gamma=\del(\star)\circ\mathrm{flip}$ and the above is equivalent to $\del(\star)$ a cocycle, see Section~\ref{subsectfunc}. 

\section{$\star$-Algebras and $\star$-Hopf algebras in (braided) bar categories}

We already introduced the notions of algebras and coalgebras in bar
categories. In this section we give some nontrivial examples and, in the braided case, bring the two notions together to the study of Hopf algebras in braided bar categories. We recall that
a Hopf algebra in a braided category means an object $B$ which is both
a unital algebra, a counital algebra, $\Delta$ is multiplicative for the braided tensor product algebra on $B\tens B$, $\eps$ is is multiplicative for which there exists an antipode $S:B\to B$ defined as for
Hopf algebras. We refer to \cite{MajBook} for details.

\subsection{Braided $*$-Hopf algebras}

\begin{defin}\label{brstardef}
A braided Hopf algebra $B$ in a braided category ${\CC}$
which is also a star category with antireal braiding is a
braided $*$-Hopf algebra if:

\noindent 1)\quad $B$ is a star algebra in ${\CC}$.

\noindent 2)\quad The coproduct satisfies 
\begin{eqnarray*}
\overline{\Delta}\,\star\,=\, \overline{\Psi}\,\Upsilon^{-1}\, 
(\star\tens\star)\Delta
:B\to \overline{B\tens B}\ .
\end{eqnarray*}
\end{defin}

\begin{propos}
If $B$ is a braided $*$-Hopf algebra, then for the antipode $S$,
$S=\star_B^{-1}\overline{S^{-1}}\star_B$. 
\end{propos}
\begin{proof} We use the uniqueness of the antipode. First,
from the definition of star algebra,
\begin{eqnarray*}
\mu(\star_B^{-1}\overline{S^{-1}}\star_B\tens\id)\Delta &=&
\star_B^{-1}\overline{\mu}\Upsilon^{-1}
(\overline{S^{-1}}\star_B\tens\star_B)\Delta \cr
&=& \star_B^{-1}\overline{\mu(\id\tens S^{-1})}\Upsilon^{-1}
(\star_B\tens\star_B)\Delta\ .
\end{eqnarray*}
Now, from the coproduct rule in \ref{brstardef},
\begin{eqnarray*}
\mu(\star_B^{-1}\overline{S^{-1}}\star_B\tens\id)\Delta 
&=& \star_B^{-1}\overline{\mu(\id\tens S^{-1})\Psi^{-1}\Delta}
\star_B\ .
\end{eqnarray*}
Now, in a purely braided Hopf algebra calculation,
\begin{eqnarray*}
S\mu(\id\tens S^{-1})\Psi^{-1}\Delta &=& \mu
\Psi(S\tens S)(\id\tens S^{-1})\Psi^{-1}\Delta \cr
&=& \mu
\Psi(S\tens \id)\Psi^{-1}\Delta \cr
&=& \mu(\id\tens S)\Delta \,=\, 1_B . \epsilon\ .
\end{eqnarray*}
{} From this,
\begin{eqnarray*}
\mu(\star_B^{-1}\overline{S^{-1}}\star_B\tens\id)\Delta 
&=& \star_B^{-1}\overline{1_B . \epsilon}\,
\star_B\cr
&=& 1_B . \epsilon\, \star_B^{-1}\star_B\,=\, 1_B . \epsilon\ ,
\end{eqnarray*}
where we have used the fact that $\star_B$ is a natural transformation
and that $1_B . \epsilon$ is a morphism from $B$ to $B$.
\end{proof}

\medskip

We shall give some examples of this definition. For the first
we refer back to \ref{hjsdgfc}, the example of ${u_q(su_2)}$, 
where $q$ is a primitive $l$th
root of unity (for odd integer $l=2k+1>1$). We shall look at 
$\mathbb{C}^2_q$ in the category of left ${u_q(su_2)}$
modules where the quantum group is viewed as a quasi-$*$-Hopf algebra.

\begin{lemma} In \ref{kfkfk} a $\star$-algebra in the category of modules of a quasi-$*$ Hopf algebra means a module algebra $A$ and a map $\star:A\to \overline{A}$
(written $a\mapsto \overline{a^*}$) such that
\[ (h\la a)^*=(Sh)^* \la a^*,\quad a^{**}=\gamma\la a,\quad 
(ab)^*=({\CG}^{-(1)}\la b^*)({\CG}^{-(2)}\la a^*)\]
 for all $a,b\in A$. 
\end{lemma} 

\begin{propos}\label{bnxcghcv}
The algebra $\mathbb{C}^2_q$ with generators $x,y$
with relations $yx=qxy$ is a $\star$-algebra in the category of
$u_q(su_2)$ modules, where
\begin{eqnarray*}
(x^m y^n)^* \, =\, q^{-(m+n)(m+n-1)(k+1)/2}\, q^{(n-m)/2}\, q^{nm}
x^m y^n\ .
\end{eqnarray*}
\end{propos}
\begin{proof} 
The algebra $\mathbb{C}^2_q$ has generators $x,y$
with relations $yx=qxy$. We define a left action of $u_q(su_2)$ by
\begin{eqnarray*}
K\la x = q\,x\ ,\ K\la y = q^{-1}\,y\ ,\ E\la x= 0=F\la y\ ,\ 
E\la y=q^{-1/2}\,x\ ,\ F\la x= q^{1/2}\,y\ .
\end{eqnarray*}
Here the square root of $q$ is defined by $q^{k+1}$ because $l=2k+1$ is odd
and $1\div 2=k+1$ modulo $l$.
The braiding in the category
${}_G{\CM}$ is given by $\Psi(v\tens w)={\CR}^{(2)}\la w \tens 
{\CR}^{(1)}\la v$, which in this case gives
\begin{eqnarray}\label{cvbjhzscvjh}
&& \Psi(x\tens x)=q^{k+1}\,x\tens x\ ,\ \Psi(y\tens y)=q^{k+1}\,y\tens y\ ,\ 
\Psi(x\tens y)=q^{k}\,y\tens x\ ,\cr
&& \Psi(y\tens x) = q^{k}\, x\tens y
+q^{k+1}(1-q^{-2})\, y\tens x\ .
\end{eqnarray}
We can make a star object from $\mathbb{C}^2_q$ by
\begin{eqnarray*}
x^* \,=\, z\, q^{-1/2}\,y\ ,\quad y^* \,=\, z\, q^{1/2}\,x\ ,
\end{eqnarray*}
where $z$ is a complex parameter with $|z|=1$. We immediately set $z=1$ for simplicity!
Of course, this only defines star on the generators, we need to define it on all 
monomials $x^m y^n$. The hint for this is that the cocycle giving
$\Upsilon$ is the same as ${\CR}$ which gives the braiding,
so the condition for a star algebra can be rewritten as
\begin{eqnarray*}
(ab)^*\,=\, \mu\,\Psi^{-1}(a^*\tens b^*)\ .
\end{eqnarray*} 
This is combined with $\mu\Psi^{-1}(a\tens b)=q^{-|a||b|(k+1)}\, a.b$, where
$|a|$ and $|b|$ are the orders of $a$ and $b$ respectively. 
Then, supposing that we can extend to a star algebra, we get the rule
on a normal ordered basis:
\begin{eqnarray*}
(x^m y^n)^* \, =\, q^{-(m+n)(m+n-1)(k+1)/2}\, q^{(n-m)/2}\, q^{nm}
x^m y^n\ .
\end{eqnarray*}
A quick induction on order shows that this is consistent with the rule
$(g\la v)^*=(Sg)^*\la v^*$, in that if it works for $v=a$ and $v=b$, then it works for $v=ab$.
Thus we get a genuine star operation and star algebra. 
\end{proof}

A much more interesting question is whether $\mathbb{C}_q^2$ can be made into
a braided Hopf algebra, with some sort of star structure. The answer without 
the star structure has been known for some time \cite{MajBook}. There is a coproduct
\begin{eqnarray*}
\Delta(x^my^n) \,=\, \sum_{}^{} 
\left[\begin{array}{c}m \\ r\end{array}\right]_q
\left[\begin{array}{c}n \\ s\end{array}\right]_q \,q^{s(m-r)}\,
x^ry^s\tens x^{m-r}y^{n-s}\ ,
\end{eqnarray*}
where the $q$-binomial coefficients are defined in terms of $q$-factorials by
\begin{eqnarray*}
\left[\begin{array}{c}m \\ r\end{array}\right]_q \,=\,
\frac{[m]_q!}{[r]_q!\,[m-r]_q!}\ ,
\end{eqnarray*}
and the $q$-factorials are recursively defined in terms of $q$-integers
$[m]_q=(q^{2m}-1)/(q^2-1)$ by 
\begin{eqnarray*}
[0]_q!\,=\,1\ ,\quad [n+1]_q!\,=\, [n+1]_q\,[n]_q!\ .
\end{eqnarray*}
If we normalise the quasitriangular structure to $\tilde{\CR}$, giving braiding
\begin{eqnarray}\label{hncvhbvgb}
&& \tilde\Psi(x\tens x)=q^{2}\,x\tens x\ ,\ \tilde\Psi(y\tens y)=q^{2}\,y\tens y\ ,\ 
\tilde\Psi(x\tens y)=q\,y\tens x\ ,\cr
&& \tilde\Psi(y\tens x) = q\, x\tens y
+q^{2}(1-q^{-2})\, y\tens x\ ,
\end{eqnarray}
then this gives a braided Hopf algebra. First we note that for $q$ a primitive 
cube root of unity ($k=1$), no renormalisation is necessary. For the other roots,
the best way to carry out this renormalisation is to introduce the central extension
of  ${u_q(su_2)}$ described in \ref{hjsdgfcent}. We suppose that
$\tilde K$ acts as $q^\beta$ on both $x$ and $y$. Then, as explained in \ref{anyhfhfh}, 
the extra factor ${\CR}_{\tilde K}$ introduces an extra factor of
$q^{\beta^2\div 2}$ into (\ref{cvbjhzscvjh}), so to get 
(\ref{hncvhbvgb}) we need $\beta^2=3$ mod $l$. This is obviously satisfied by 
$\beta=0$ if $l=3$, so as already noted no extension is necessary. For
$l>3$ the existence of a $\beta$ with $\beta^2=3$ mod $l$ is well studied in 
number theory. If there is such a $\beta$, then $3$ is said to be a quadratic residue for $l$.
For $l$ prime, this is summarised in the Legendre symbol:
\begin{eqnarray*}
\left(\frac{3}{l}\right) \,=\, 
\left\{\begin{array}{cc}+1 & \exists \beta:\beta^2=3
\ \mathrm{mod}\ l \\-1 & \nexists \beta:\beta^2=3
\ \mathrm{mod}\ l\end{array}\right.
\end{eqnarray*}
The list of those primes in the first hundred prime numbers 
for which $3$ is a quadratic residue is:
\begin{eqnarray*}
&& 11  ,\,13  ,\,23  ,\,37  ,\,47  ,\,59  ,\,61  ,\,71  ,\,73 
 ,\,83  ,\,97  ,\,107  ,\,109  ,\,131  ,\,157  ,\,167  ,\,179 
 ,\,181  ,\,191  ,\,193 \cr  && 227  ,\,229  ,\,239  ,\,241  
,\,251  ,\,263  ,\,277  ,\,311  ,\,313  ,\,337  ,\,347  ,\,349 
 ,\,359  ,\,373  ,\,383  ,\,397  ,\,409  \cr && 419  ,\,421  
,\,431  ,\,433  ,\,443  ,\,457  ,\,467  ,\,479  ,\,491  ,\,503 
 ,\,541 \ .
\end{eqnarray*}
For example, in the first case on the list, $5^2=3$ mod $11$.

\begin{theorem}
In the following cases the coproduct
\begin{eqnarray*}
\Delta(x^my^n) \,=\, \sum_{}^{} 
\left[\begin{array}{c}m \\ r\end{array}\right]_q
\left[\begin{array}{c}n \\ s\end{array}\right]_q \,q^{s(m-r)}\,
x^ry^s\tens x^{m-r}y^{n-s}
\end{eqnarray*}
makes $\mathbb{C}^2_q$ into a $\star$-braided Hopf algebra, with braiding
\begin{eqnarray}
&& \tilde\Psi(x\tens x)=q^{2}\,x\tens x\ ,\ \tilde\Psi(y\tens y)=q^{2}\,y\tens y\ ,\ 
\tilde\Psi(x\tens y)=q\,y\tens x\ ,\cr
&& \tilde\Psi(y\tens x) = q\, x\tens y
+q^{2}(1-q^{-2})\, y\tens x\ .
\end{eqnarray}
and the star operation
\begin{eqnarray*}
(x^m y^n)^* \, =\, q^{-(m+n)(m+n-1)}\, q^{(n-m)/2}\, q^{nm}
x^m y^n\ .
\end{eqnarray*}

\medskip\noindent \textbf{Case 1:} $l=3$ and we use the usual ${u_q(su_2)}$. 

\medskip \noindent\textbf{Case 2:} We use the centrally extended $\widetilde{u_q(su_2)}$,
and restrict to the case where $3$ is a quadratic residue mod $l$. The central generator
$\tilde K$ acts as $q^\beta$ on $x$ and $y$, where $\beta^2=3$ mod $l$.

\end{theorem}
\begin{proof} Explicit calculation, using Examples
\ref{hjsdgfc} and \ref {hjsdgfcent}.\end{proof}

\subsection{Transmutation of Hopf algebras}

We recall that the transmutation of a quasitriangular Hopf algebra $H$ is a braided 
Hopf algebra in the braided category ${}_H{\CM}$
of left $H$-modules. It is defined by $BH$ being $H$ as a left
$H$ module (under the left adjoint action
$h\la b=h_{(1)}bS(h_{(2)})$, and has the same algebra, unit and counit as $H$. 
For the braided coalgebra, we have
\begin{eqnarray*}
\underline{\Delta} b \,=\, b_{(1)}\, S({\CR}^{(2)})\tens 
{\CR}^{(1)}\la b_{(2)}\ ,\quad
\underline{S}(b)\,=\, {\CR}^{(2)}\, S({\CR}^{(1)}\la b)\ .
\end{eqnarray*}
We shall see that, under some quite common restrictions, the transmutation
of a Hopf algebra with star operation gives a braided $*$-Hopf algebra.

\begin{theorem}
Suppose that $H$ is a quasitriangular $*$-Hopf algebra, and that
${\CR}^{*\tens *}={\CR}^{-1}$. Then
if we define $\star:BH\to \overline{BH}$ by $b\mapsto \overline{b^*}$
(using the star operation from $H$), then $BH$ is a 
braided $*$-Hopf algebra.
\end{theorem}
\begin{proof} The reader should remember that in this case
$*$ preserves the order of coproduct
and switches $S$ and $S^{-1}$. 
First check that 
$\star:BH\to \overline{BH}$ is a morphism in ${}_H{\CM}$:
\begin{eqnarray*}
\star(h\la b) &=& \overline{(h\la b)^*}\,=\, \overline{(h_{(1)} b S(h_{(2)}))^*} \cr
&=& \overline{ S(h_{(2)})^* b^* h_{(1)}^* } \ ,\cr 
h\la (\star b)&=& \overline{ S(h)^* \la b^*  }\,=\, 
\overline{ S(h)_{(1)}^*  b^* S(S(h)_{(2)}^*) } \cr
&=& \overline{ S(h_{(2)})^*  b^* S(S(h_{(1)})^*) } \cr
&=& \overline{ S(h_{(2)})^*  b^* S^{-1}(S(h_{(1)}))^* } \ .
\end{eqnarray*}
Next for the coproduct,
\begin{eqnarray*}
(\star\tens\star)\underline{\Delta} b &=& \overline{S({\CR}^{(2)})^*b_{(1)}^*}
\tens \overline{({\CR}^{(1)}\la b_{(2)})^*} \cr
&=&  \overline{S({\CR}^{(2)})^*b_{(1)}^*}
\tens \overline{S({\CR}^{(1)})^*\la b_{(2)}^*} \cr
&=&  \overline{{\CR}^{(2)*}b_{(1)}^*}
\tens \overline{{\CR}^{(1)*}\la b_{(2)}^*} \cr
\Upsilon^{-1}(\star\tens\star)\underline{\Delta} b &=& 
\overline{{\CR}^{(1)*}\la b_{(2)}^* \tens {\CR}^{(2)*}b_{(1)}^* } \cr
&=& 
\overline{{\CR}^{-(1)}\la b_{(2)}^* \tens {\CR}^{-(2)}b_{(1)}^* } \cr
\overline{\Psi}\Upsilon^{-1}(\star\tens\star)\underline{\Delta} b &=& 
\overline{  {\CR}^{(2)}\la ({\CR}^{-(2)}b_{(1)}^*)
\tens {\CR}^{(1)} {\CR}^{-(1)}\la b_{(2)}^* }\cr
&=& \overline{  {\CR}^{(2)}\la ({\CR}^{-(2)}b_{(1)}^*)
\tens {\CR}^{(1)} {\CR}^{-(1)}\la b_{(2)}^*} \cr
&=& \overline{  {\CR}^{(2)}_{\phantom{(2)}(1)}
{\CR}^{-(2)}b_{(1)}^* S({\CR}^{(2)}_{\phantom{(2)}(2)})
\tens {\CR}^{(1)} {\CR}^{-(1)}\la b_{(2)}^*}\ .
\end{eqnarray*}
Now from the definition of quasitriangular structure use $(\id\tens\Delta){\CR}
={\CR}_{13}{\CR}_{12}$ to get 
\begin{eqnarray*}
\overline{\Psi}\Upsilon^{-1}(\star\tens\star)\underline{\Delta} b &=&
\overline{ b_{(1)}^* S({\CR}^{(2)})
\tens {\CR}^{(1)} \la b_{(2)}^*}\ .
\end{eqnarray*}
Finally we check the antipode,
\begin{eqnarray*}
\underline{S}(b) &=& {\CR}^{(2)}\, S({\CR}^{(1)}\la b) \cr
&=& {\CR}^{(2)}\, S^2({\CR}^{(1)}_{\phantom{(1)}(2)})\, S(b)
\, S({\CR}^{(1)}_{\phantom{(1)}(1)})\ ,
\end{eqnarray*}
and using $(\Delta\tens\id){\CR}
={\CR}_{13}{\CR}_{23}$,
\begin{eqnarray*}
\underline{S}(b) 
&=&\tilde {\CR}^{(2)}\,{\CR}^{(2)}\, 
S^2({\CR}^{(1)})\, S(b)
\, S(\tilde{\CR}^{(1)}) \cr
&=&\tilde {\CR}^{(2)}\,u^{-1}\, S(b)
\, S(\tilde{\CR}^{(1)})   \cr
&=& {\CR}^{-(2)}\,u^{-1}\, S(b)
\, {\CR}^{-(1)}\ ,\cr
(\underline{S}b)^*
&=& {\CR}^{(1)}\, S(b)^*\, u^{-*}\, {\CR}^{(2)}\ ,\cr
\underline{S}((\underline{S}b)^*)
&=& \tilde{\CR}^{(2)}\, S^2(\tilde{\CR}^{(1)}_{\phantom{(1)}(2)})\, 
S({\CR}^{(1)}\, S(b)^*\, u^{-*}\, {\CR}^{(2)})
\, S(\tilde{\CR}^{(1)}_{\phantom{(1)}(1)}) \cr
&=& \tilde{\CR}^{(2)}\, S^2(\tilde{\CR}^{(1)}_{\phantom{(1)}(2)})\, 
S({\CR}^{(2)})\, S(u^{-*})\, b^*\, S({\CR}^{(1)})
\, S(\tilde{\CR}^{(1)}_{\phantom{(1)}(1)}) \cr
&=& \tilde{\CR}^{(2)}\, \hat{\CR}^{(2)}\, 
S^2(\hat{\CR}^{(1)})\, 
S({\CR}^{(2)})\, S(u^{-*})\, b^*\, S({\CR}^{(1)})
\, S(\tilde{\CR}^{(1)}) \cr
&=& \tilde{\CR}^{(2)}\, u^{-1}\, 
S({\CR}^{(2)})\, S(u^{-*})\, b^*\, S({\CR}^{(1)})
\, S(\tilde{\CR}^{(1)}) \cr
&=& \tilde{\CR}^{(2)}\, 
S^{-1}({\CR}^{(2)})\, u^{-1}\, S(u^{-*})\, b^*\, S({\CR}^{(1)})
\, S(\tilde{\CR}^{(1)}) \cr
&=& \tilde{\CR}^{(2)}\, 
S^{-1}({\CR}^{(2)})\, u^{-1}\, S(u^{-*})\, b^*\,
 S(\tilde{\CR}^{(1)}\, {\CR}^{(1)})   \cr
&=& \tilde{\CR}^{(2)}\, 
{\CR}^{-(2)}\, u^{-1}\, S(u^{-*})\, b^*\,
 S(\tilde{\CR}^{(1)}\, {\CR}^{-(1)}) \cr
&=& u^{-1}\, S(u^{-*})\, b^*\ .
\end{eqnarray*}\end{proof}

\begin{theorem}
Suppose that $H$ is a quasitriangular flip-$*$-Hopf algebra, and that
${\CR}^{*\tens *}={\CR}_{21}^{-1}$. Then
if we define $\star:BH\to \overline{BH}$ by 
\begin{eqnarray*}
\star b \,=\, \overline{{\CR}^{-(1)} b^* S^{-1}({\CR}^{-(2)})S^{-1}(u)}\ .
\end{eqnarray*}
(using the star operation from $H$), then $BH$ is a 
braided $*$-Hopf algebra. 
\end{theorem}
\begin{proof} 
The reader should remember that in this case
$*$ flips the order of the coproduct
and  $S(h^*)=S(h)^*$. First  
we check that $\star$ is a morphism: 
\begin{eqnarray*}
h\la \star b &=& 
\overline{ S(h)^*_{\phantom{*}(1)}
{\CR}^{-(1)} b^* S^{-1}({\CR}^{-(2)})S^{-1}(u)
S(S(h)^*_{\phantom{*}(2)})} \cr
&=& 
\overline{ S(h^*_{\phantom{*}(2)})
{\CR}^{-(1)} b^* S^{-1}({\CR}^{-(2)})S^{-1}(u)
S^2(h^*_{\phantom{*}(1)})} \cr
&=& 
\overline{ S(S^{-1}({\CR}^{-(1)})h^*_{\phantom{*}(2)})
 b^* S^{-1}({\CR}^{-(2)})h^*_{\phantom{*}(1)} S^{-1}(u)} \cr
 &=& 
\overline{ S({\CR}^{-(1)}h^*_{\phantom{*}(2)})
 b^* {\CR}^{-(2)}h^*_{\phantom{*}(1)} S^{-1}(u)} \cr
  &=& 
\overline{ S(h^*_{\phantom{*}(1)}{\CR}^{-(1)})
 b^*h^*_{\phantom{*}(2)} {\CR}^{-(2)}S^{-1}(u)} \cr
   &=& 
\overline{ S({\CR}^{-(1)})
S(h^*_{\phantom{*}(1)})
 b^*h^*_{\phantom{*}(2)} S^{-1}(S({\CR}^{-(2)}))S^{-1}(u)} \cr
    &=& 
\overline{ {\CR}^{-(1)}
S(h^*_{\phantom{*}(1)})
 b^*h^*_{\phantom{*}(2)} S^{-1}({\CR}^{-(2)})S^{-1}(u)} \cr
     &=& 
\overline{ {\CR}^{-(1)}
S(h^{\phantom{(2)}*}_{(2)})
 b^*h^{\phantom{(1)}*}_{(1)} S^{-1}({\CR}^{-(2)})S^{-1}(u)} \cr
      &=& 
\overline{ {\CR}^{-(1)}(h\la b)^*
 S^{-1}({\CR}^{-(2)})S^{-1}(u)} \ ,
\end{eqnarray*}
where we have used $(S\tens S){\CR}={\CR}$. 

Next, we do $\star$ twice:
\begin{eqnarray*}
\star\star b &=&  \star(\overline{{\CR}^{-(1)} b^* 
S^{-1}({\CR}^{-(2)})S^{-1}(u)}) \cr
&=&
\overline{\overline{\tilde{\CR}^{-(1)}
S^{-1}(u^*) S^{-1}({\CR}^{-(2)*}) b {\CR}^{-(1)*} S^{-1}
(\tilde{\CR}^{-(2)})S^{-1}(u)}} \cr
&=&
\overline{\overline{\tilde{\CR}^{-(1)}
S^{-1}(u^{-1}) S^{-1}({\CR}^{(1)}) b {\CR}^{(2)} S^{-1}
(\tilde{\CR}^{-(2)})S^{-1}(u)}} \cr
&=&
\overline{\overline{S^{-1}(u^{-1})   S^{-2}(\tilde{\CR}^{-(1)})
S^{-1}({\CR}^{(1)}) b {\CR}^{(2)} S^{-1}
(\tilde{\CR}^{-(2)})S^{-1}(u)}} \cr
&=&
\overline{\overline{S^{-1}(u^{-1})   
S^{-1}({\CR}^{(1)} S^{-1}(\tilde{\CR}^{-(1)})) 
b {\CR}^{(2)} S^{-1}
(\tilde{\CR}^{-(2)})S^{-1}(u)}} \cr
&=&
\overline{\overline{S^{-1}(u^{-1})   
S^{-1}({\CR}^{(1)} \tilde{\CR}^{-(1)}) 
b {\CR}^{(2)} \tilde{\CR}^{-(2)}S^{-1}(u)}} \cr
&=&
\overline{\overline{S^{-1}(u^{-1})  b S^{-1}(u)}} 
\end{eqnarray*}
using $u^*=u^{-1}$ and $(S\tens S){\CR}={\CR}$. As $\gamma=v^{-1}\nu$
is grouplike, we get
\begin{eqnarray*}
\gamma\la b &=& \gamma b \gamma^{-1} \cr
&=& v^{-1}\nu b \nu^{-1}v\,=\, v^{-1} b v\ ,
\end{eqnarray*}
using $\nu$ central. Now use $S^{-1}(u)=v$.

Now check that we get a star algebra. The following expression should be
$\star(ab)$:
\begin{eqnarray*}
&&\bar\mu\Upsilon^{-1}(
\overline{\hat{\CR}^{-(1)} a^* S^{-1}(\hat{\CR}^{-(2)})S^{-1}(u)}
\tens 
\overline{\tilde{\CR}^{-(1)} b^* S^{-1}(\tilde{\CR}^{-(2)})S^{-1}(u)}) \cr
&=& \overline{
(R^{-(1)}\la(\tilde{\CR}^{-(1)} b^* S^{-1}(\tilde{\CR}^{-(2)})S^{-1}(u)))\ 
(R^{-(2)}\la(\hat{\CR}^{-(1)} a^* S^{-1}(\hat{\CR}^{-(2)})S^{-1}(u)))
} \cr
&=&
 \overline{
R^{-(1)}_{\phantom{-(1)}(1)}
\tilde{\CR}^{-(1)} b^* S^{-1}(\tilde{\CR}^{-(2)})S^{-1}(u)
S(R^{-(1)}_{\phantom{-(1)}(2)})
R^{-(2)}_{\phantom{-(2)}(1)}
\hat{\CR}^{-(1)} a^* S^{-1}(\hat{\CR}^{-(2)})S^{-1}(u)
S({\CR}^{-(2)}_{\phantom{-(2)}(2)})
} \cr
&=&
 \overline{
R^{-(1)}_{\phantom{-(1)}(1)}
\tilde{\CR}^{-(1)} b^* S^{-1}(\tilde{\CR}^{-(2)})S^{-1}(u)
S(R^{-(1)}_{\phantom{-(1)}(2)})
R^{-(2)}_{\phantom{-(2)}(1)}
\hat{\CR}^{-(1)} a^* S^{-1}( {\CR}^{-(2)}_{\phantom{-(2)}(2)}
\hat{\CR}^{-(2)})S^{-1}(u)
} 
\end{eqnarray*}
{}From the usual rules for the quasitriangular structure we get
\begin{eqnarray*}
(\Delta\tens\Delta){\CR} \,=\, {\CR}_{14}
{\CR}_{13}{\CR}_{24}{\CR}_{23}\ ,
\end{eqnarray*}
so 
\begin{eqnarray*}
&& {\CR}^{-(1)}_{\phantom{-(1)}(1)} \tens
S^{-1}(u)
S({\CR}^{-(1)}_{\phantom{-(1)}(2)})
{\CR}^{-(2)}_{\phantom{-(2)}(1)} \tens 
 {\CR}^{-(2)}_{\phantom{-(2)}(2)} \cr
 &=& \breve{{\CR}}^{-(1)}{\CR}^{-(1)} 
 \tens S^{-1}(u)S(  \ddot{{\CR}}^{-(1)}  \dot{{\CR}}^{-(1)})
\ddot{{\CR}}^{-(2)} \breve{{\CR}}^{-(2)}
  \tens \dot{{\CR}}^{-(2)}  {\CR}^{-(2)} \cr
   &=& \breve{{\CR}}^{-(1)}{\CR}^{-(1)} 
 \tens S^{-1}(u)S(    \dot{{\CR}}^{-(1)}) S^{-1}(u^{-1})
 \breve{{\CR}}^{-(2)}
  \tens \dot{{\CR}}^{-(2)}  {\CR}^{-(2)} \cr
     &=& \breve{{\CR}}^{-(1)}{\CR}^{-(1)} 
 \tens S^{-1}(    \dot{{\CR}}^{-(1)}) 
 \breve{{\CR}}^{-(2)}
  \tens \dot{{\CR}}^{-(2)}  {\CR}^{-(2)} \cr
       &=& \breve{{\CR}}^{-(1)}{\CR}^{-(1)} 
 \tens S^{-1}(  S( \breve{{\CR}}^{-(2)})  \dot{{\CR}}^{-(1)}) 
  \tens \dot{{\CR}}^{-(2)}  {\CR}^{-(2)} \cr
         &=& \breve{{\CR}}^{(1)}{\CR}^{-(1)} 
 \tens S^{-1}(   \breve{{\CR}}^{(2)} \dot{{\CR}}^{-(1)}) 
  \tens \dot{{\CR}}^{-(2)}  {\CR}^{-(2)} \cr
\end{eqnarray*}
Now, substitute this back into the expression we hope is $\star(ab)$ to get:
\begin{eqnarray*}
&& \overline{
\breve{{\CR}}^{(1)}{\CR}^{-(1)} 
\tilde{\CR}^{-(1)} b^* S^{-1}(\tilde{\CR}^{-(2)})S^{-1}(   \breve{{\CR}}^{(2)} \dot{{\CR}}^{-(1)}) 
\hat{\CR}^{-(1)} a^* S^{-1}(  \dot{{\CR}}^{-(2)}  {\CR}^{-(2)}
\hat{\CR}^{-(2)})S^{-1}(u)} \cr
&=&  \overline{
\breve{{\CR}}^{(1)}{\CR}^{-(1)} 
\tilde{\CR}^{-(1)} b^* 
S^{-1}(   \breve{{\CR}}^{(2)} \dot{{\CR}}^{-(1)} \tilde{\CR}^{-(2)}) 
\hat{\CR}^{-(1)} a^* S^{-1}(  \dot{{\CR}}^{-(2)}  {\CR}^{-(2)}
\hat{\CR}^{-(2)})S^{-1}(u)} \ .
\end{eqnarray*}
{} From this we have, using the braid relations,
\begin{eqnarray*}
 \breve{{\CR}}_{12}    \dot{{\CR}}_{23}^- 
 {{\CR}}_{13}^-
   \tilde{{\CR}}_{12}^- \,=\,  {{\CR}}_{13}^-  \dot{{\CR}}_{23}^- 
\end{eqnarray*}
which gives, as required,
\begin{eqnarray*}
&&  \overline{{\CR}^{-(1)}  b^* 
S^{-1}(  \dot{{\CR}}^{-(1)} ) 
\hat{\CR}^{-(1)} a^* S^{-1}(  {\CR}^{-(2)} \dot{{\CR}}^{-(2)} 
\hat{\CR}^{-(2)})S^{-1}(u)} \cr
&=&   \overline{{\CR}^{-(1)}  b^* 
  \dot{{\CR}}^{(1)} 
\hat{\CR}^{-(1)} a^* S^{-1}(  {\CR}^{-(2)} \dot{{\CR}}^{(2)} 
\hat{\CR}^{-(2)})S^{-1}(u)} \cr
&=&   \overline{{\CR}^{-(1)}  b^* 
 a^* S^{-1}(  {\CR}^{-(2)})S^{-1}(u)} \cr
 &=& \star(ab)\ .
\end{eqnarray*}

Now we check the rule for the coproduct: Using $Q={\CR}_{21}{\CR}$,
\begin{eqnarray*}
\Upsilon\overline{\Psi^{-1}\underline{\Delta}}\star a &=& 
\overline{
{\CR}^{-(1)}_{\phantom{-(2)}(1)} a^*_{\phantom{*}(1)} 
S^{-1}({\CR}^{-(2)})_{(1)} S^{-1}(u)_{(1)} S(\tilde{{\CR}}^{(2)})} \cr
&& \tens
\overline{ \tilde{{\CR}}^{(1)}\la(
{\CR}^{-(1)}_{\phantom{-(2)}(2)} a^*_{\phantom{*}(2)} 
S^{-1}({\CR}^{-(2)})_{(2)} S^{-1}(u)_{(2)} 
)}    \cr
 &=& 
\overline{
{\CR}^{-(1)}_{\phantom{-(2)}(1)} a^*_{\phantom{*}(1)} 
S^{-1}({\CR}^{-(2)}_{\phantom{-(2)}(2)})
 S^{-1}(Q^{-(2)}) S^{-1}(u) 
S(\tilde{{\CR}}^{(2)})} \cr
&& \tens
\overline{ \tilde{{\CR}}^{(1)}_{\phantom{(1)}(1)}
{\CR}^{-(1)}_{\phantom{-(2)}(2)} a^*_{\phantom{*}(2)} 
S^{-1}({\CR}^{-(2)}_{\phantom{-(2)}(1)}) S^{-1}(Q^{-(1)}) S^{-1}(u) 
S (\tilde{{\CR}}^{(1)}_{\phantom{(1)}(2)})}  \cr
 &=& 
\overline{
{\CR}^{-(1)}_{\phantom{-(2)}(1)} a^*_{\phantom{*}(1)} 
S^{-1}({\CR}^{-(2)}_{\phantom{-(2)}(2)})
 S^{-1}(Q^{-(2)})
S^{-1}(\tilde{{\CR}}^{(2)}) S^{-1}(u)}  \cr
&& \tens
\overline{ \tilde{{\CR}}^{(1)}_{\phantom{(1)}(1)}
{\CR}^{-(1)}_{\phantom{-(2)}(2)} a^*_{\phantom{*}(2)} 
S^{-1}({\CR}^{-(2)}_{\phantom{-(2)}(1)}) S^{-1}(Q^{-(1)}) 
S^{-1} (\tilde{{\CR}}^{(1)}_{\phantom{(1)}(2)}) S^{-1}(u)  }  \cr
 &=& 
\overline{
{\CR}^{-(1)}_{\phantom{-(2)}(1)} a^*_{\phantom{*}(1)} 
S^{-1}(\tilde{{\CR}}^{(2)}
Q^{-(2)}{\CR}^{-(2)}_{\phantom{-(2)}(2)})
S^{-1}(u)}  \cr
&& \tens
\overline{ \tilde{{\CR}}^{(1)}_{\phantom{(1)}(1)}
{\CR}^{-(1)}_{\phantom{-(2)}(2)} a^*_{\phantom{*}(2)} 
S^{-1}( \tilde{{\CR}}^{(1)}_{\phantom{(1)}(2)}
Q^{-(1)}{\CR}^{-(2)}_{\phantom{-(2)}(1)}) 
S^{-1}(u)  }  \cr
 &=& 
\overline{
{\CR}^{-(1)}_{\phantom{-(2)}(1)} a^*_{\phantom{*}(1)} 
S^{-1}(\tilde{{\CR}}^{(2)}
{\CR}^{-(2)}_{\phantom{-(2)}(2)}Q^{-(2)})
S^{-1}(u)}  \cr
&& \tens
\overline{ \tilde{{\CR}}^{(1)}_{\phantom{(1)}(1)}
{\CR}^{-(1)}_{\phantom{-(2)}(2)} a^*_{\phantom{*}(2)} 
S^{-1}( \tilde{{\CR}}^{(1)}_{\phantom{(1)}(2)}
{\CR}^{-(2)}_{\phantom{-(2)}(1)}Q^{-(1)}) 
S^{-1}(u)  }  \ ,
\end{eqnarray*}
where we have used the fact that $Q$ commutes with coproducts. Now we
use the braid relations again:
\begin{eqnarray*}
(1\tens(\Delta\tens\id){\CR})((\Delta\tens\Delta){\CR}^{-1})
(1\tens 1\tens Q^{-1}) &=&
{\CR}_{24} {\CR}_{34} {\CR}_{23}^-
{\CR}_{24}^-
{\CR}_{13}^-{\CR}_{14}^-{\CR}_{34}^-
{\CR}_{43}^- \cr
&=&
{\CR}_{23}^-{\CR}_{34}
{\CR}_{13}^-{\CR}_{14}^-{\CR}_{34}^-
{\CR}_{43}^- \cr
&=&
\tilde{\CR}_{23}^-{\CR}_{14}^- \hat{\CR}_{13}^-
\check{{\CR}}_{43}^- \ .
\end{eqnarray*}
On substitution we get
\begin{eqnarray*}
\Upsilon\overline{\Psi^{-1}\underline{\Delta}}\star a &=& 
\overline{
{\CR}^{-(1)} \hat{\CR}^{-(1)} a^*_{\phantom{*}(1)} 
S^{-1}({{\CR}}^{-(2)} \check{\CR}^{-(1)})
S^{-1}(u)}  \cr
&& \tens
\overline{ \tilde{{\CR}}^{-(1)}  a^*_{\phantom{*}(2)} 
S^{-1}( \tilde{{\CR}}^{-(2)} \hat{{\CR}}^{-(2)}
\check{\CR}^{-(2)})
S^{-1}(u)  }  \cr
&=& 
\overline{
{\CR}^{-(1)} \hat{\CR}^{-(1)} a^*_{\phantom{*}(1)} 
 \check{\CR}^{-(1)} S^{-1}({{\CR}}^{-(2)})
S^{-1}(u)}  \cr
&& \tens
\overline{ \tilde{{\CR}}^{-(1)}  a^*_{\phantom{*}(2)} 
\check{\CR}^{-(2)} S^{-1}( \hat{{\CR}}^{-(2)}
) S^{-1}( \tilde{{\CR}}^{-(2)})
S^{-1}(u)  }  \cr
&=& 
\overline{
{\CR}^{-(1)} \hat{\CR}^{-(1)} 
 \check{\CR}^{-(1)} a^*_{\phantom{*}(2)}  S^{-1}({{\CR}}^{-(2)})
S^{-1}(u)}  \cr
&& \tens
\overline{ \tilde{{\CR}}^{-(1)} 
\check{\CR}^{-(2)}  a^*_{\phantom{*}(1)}  S^{-1}( \hat{{\CR}}^{-(2)}
) S^{-1}( \tilde{{\CR}}^{-(2)})
S^{-1}(u)  }  \ .
\end{eqnarray*}
Now apply the inverse of $\star$:
\begin{eqnarray*}
(\star^{-1}\tens \star^{-1})
\Upsilon\overline{\Psi^{-1}\underline{\Delta}}\star a 
&=& 
( \hat{\CR}^{-(1)} 
 \check{\CR}^{-(1)} a^*_{\phantom{*}(2)})^*   \tens
(\check{\CR}^{-(2)}  a^*_{\phantom{*}(1)}  S^{-1}( \hat{{\CR}}^{-(2)}))^* \cr
&=& a_{(1)} \check{\CR}^{(2)} \hat{\CR}^{(2)} \tens 
S^{-1}( \hat{{\CR}}^{(1)}) a_{(2)} \check{\CR}^{(1)} \cr
&=& a_{(1)} S(\check{\CR}^{(2)}) S(\hat{\CR}^{(2)}) \tens 
 \hat{{\CR}}^{(1)} a_{(2)} S(\check{\CR}^{(1)}) \cr
 &=& a_{(1)}  S(\hat{\CR}^{(2)}\check{\CR}^{(2)}) \tens 
 \hat{{\CR}}^{(1)} a_{(2)} S(\check{\CR}^{(1)}) \cr
  &=& a_{(1)}  S({\CR}^{(2)}) \tens 
 {\CR}^{(1)}_{\phantom{(1)}(1)} a_{(2)} S({\CR}^{(1)}_{\phantom{(1)}(1)}) \cr
 &=& \underline{\Delta} a\ .
\end{eqnarray*}\end{proof}

\subsection{Example: The octonions}

We construct these as a $\star$ algebra using the setting of Theorem~\ref{mdvfmasm}
and the following immediate corollary:

\begin{propos}
In the setting of Theorem \ref{mdvfmasm}, if $A$ is a 
star algebra in the category ${}_{H}{\CM}$ then 
$F(A)$ is a star algebra $A$ in the category ${}_{H^{{\CF}}}{\CM}$,
where we set $\star:F(a)\mapsto \overline{F(\varphi^{-1}\la a^*)}$. 
\end{propos}
\begin{proof} Using Theorem \ref{mdvfmasm}, we apply
the twisting functor to any commutative diagram in ${}_H\CM$ to obtain one in ${}_{H^{{\CF}}}\CM$ . \end{proof}

\medskip
In particular, consider $H=\C(G)$, the functions on a finite group. This has  
basis of delta-functions $\{\delta_a\}$ labelled by 
$a\in G$ and coproduct $\Delta \delta_a=\sum_{bc=a} \delta_b\tens  
\delta_c$, counit $\eps\delta_a=\delta_{a,e}$ and antipode $S\delta_a= 
\delta_{a^{-1}}$. Here $e$ is the group identity.  A cochain on $H$ is a suitable $\CF\in H\tens H$ i.e. a nowhere 
vanishing 2-argument function $\CF(a,b)$ on the group with value $1$ 
when either argument is the group identity $e$. Then 
\[ \phi(a,b,c)=\frac{\CF(b,c)\CF(a,bc)}{\CF(ab,c)\CF(a,b)}\] 
is the usual group-cohomology coboundary of $F$ and is a group 
3-cocycle. Then $H^{\CF}$ is the same algebra and coalgebra as $H$ but 
is viewed as a quasi-Hopf algebra with this $\phi$. 

We next take $A=\C G$, the group algebra of $G$. This has basis $\{e_a\}$  
labelled again by group elements. The product is just the product of  
$G$, so $e_ae_b=e_{ab}$.  This is covariant under $\C(G)$ with action 
\[ \delta_a\la  e_b=\delta_{a,b}e_b.\] 
Applying the twisting functor gives the twisted group 
algebra  $A_{\CF}$ with the new product 
\[ e_a\bullet e_b= \CF^{-1}(a,b) e_{ab}\]

An example is $G=\Z_2^3$ which we write additively as $3$-vectors $ 
\vec{a}$ with entries in $\Z_2$. 
We take 
\[ \CF(\vec a,\vec b)=(-1)^{\vec a{}^T 
\left(\begin{array}{ccc}1 & 1 & 1 \\0 & 1 & 1 \\0 & 0 & 1\end{array}\right)
 \vec b+ a_1b_2b_3+b_1a_2b_3 
+b_1b_2a_3},\quad \phi(\vec a,\vec b,\vec c)=(-1)^{\vec a\cdot( \vec b 
\times\vec c)} \] 
The new product \[ e_{\vec a}\bullet e_{\vec b}=\CF(\vec a,\vec 
b)e_{\vec a+\vec b}\] is that of the octonions $\Bbb O$ as 
explained in \cite{albmajoct}. This is the by-now well-known construction of the octonions as a quasi-algebra
obtained by twisiting. 

\begin{propos} The cochain $\CF$ for the octonions obeys the condition in Remark~\ref{twicoc}
and hence the category above of $\Z_2^3$-graded spaces is a bar category with 
\[\mathrm{bb}_V(v)=\overline{\overline{v}},\quad  \Upsilon_{V,W}(\overline{v\tens w})=\bar w\tens \bar v
\left\{\begin{array}{cc}-1 & 0\ne  |v|\ne |w|\ne 0 \\1 & \mathrm{else}\end{array}\right.
\]
on elements of homogeneous degree in objects $V,W$. Moreover, the Octonions $\mathbb{O}$ are
a $\star$-algebra in this category with $e_{\vec a}^*=e_{-\vec a}$ in the above basis.
\end{propos}
\begin{proof} The element $\CG=\CF_{21}\CF^{-1}$ is also the element that makes the category symmetric and was computed in \cite{albmajoct} to give the above. The $\varphi$ is trivial hence the $\star$ operation is the same as that on $\C\Z_2^3$. \end{proof}

\medskip
Clearly we can apply this twisting theory also to the action of $\C \Z_2^3$ on 
$\CH=l^2(\Z_2^3)$ obtaining in this way an inner product space of some kind 
in the bar category where the octonions live and an acton of the octonions 
on it is as a star-algebra, i.e. a C*-algebra in the bar category.

\subsection{A C*-algebra in the coset representative category}
There is an algebra $A_{\bowtie}$ in the category ${\CC}_{\bowtie}$
defined as follows. It has vector space basis $\delta_s\tens u$ for all
$s\in M$ and $u\in G$. The grade is defined to be the
unique solution to $s\cdot\langle \delta_s\tens u\rangle =s\ra u$,
and if we set $a=\langle \delta_s\tens u\rangle$, the action is given by
\begin{eqnarray}
(\delta_s\tens u)\bar\ra v\ =\ \delta_{s\ra(a\la v)}\tens (a\la v)^{-1}uv\ .\label{etal}
\end{eqnarray}
The formula for the product $\mu$ for $A_{\bowtie}$ consistent with the action
above is
\[
(\delta_s\tens u)(\delta_t\tens v)\ =\ \delta_{t,s\ra u}\ \delta_{s\ra\tau(a,b)}
\tens \tau(a,b)^{-1}uv\ ,
\]
where $a=\langle \delta_s\tens u\rangle $ and $b=\langle \delta_t\tens v\rangle $. 
There is a unit and counit defined by 
\[
I\ =\ \sum_{t}\delta_t\tens e\ ,\quad \epsilon(\delta_s\tens u)\ =\ \delta_{s,e}\ .
\]
The reader should refer to \cite{cosrep} for the details of the construction

\begin{propos}
The following morphism $\star:A_{\bowtie} \to\overline{A_{\bowtie}}$ 
makes $A_{\bowtie}$ into a star object in ${\CC}_{\bowtie}$,
where $a=\langle \delta_s\tens u\rangle$:
\begin{eqnarray} \label{Astardef}
*(\delta_s\tens u) \,=\, \overline{(\delta_s\tens u)^*}
\,=\, \overline{(\delta_{s\ra u}\tens u^{-1}\tau(a,a^R))\bar\ra \tau(a,a^R)^{-1}}\ .
\end{eqnarray}
\end{propos}
\begin{proof}
Begin with
\begin{eqnarray*}
(s\ra u) \!\cdot\! a^R \,=\, (s\!\cdot\! a) \!\cdot\! a^R \,=\, s\ra\tau(a,a^R)\!\cdot\! 
(a\!\cdot\! a^R)\,=\, s\ra\tau(a,a^R)\ .
\end{eqnarray*}
This means that $\delta_{s\ra u}\tens u^{-1}\tau(a,a^R)$ has grade $a^R$. Now we use
the formula $a^L=a^R\ra \tau(a,a^R)^{-1}$ to see that the following element,
which we define to be $(\delta_s\tens u)^*$, has grade
$a^L$:
\begin{eqnarray*}
(\delta_s\tens u)^* &=& (\delta_{s\ra u}\tens u^{-1}\tau(a,a^R))\ra \tau(a,a^R)^{-1}\cr
&=& \delta_{s\ra u\tau(a^L,a)^{-1}}\tens \tau(a^L,a)u^{-1} \ .
\end{eqnarray*}
To get the last equality we use $a^r\la \tau(a,a^R)^{-1}= \tau(a^L,a)^{-1}$. 
Now we define a map $\star:A_{\bowtie}\to \overline{A_{\bowtie}}$ by 
\begin{eqnarray} \label{Astardef675}
*(\delta_s\tens u) \,=\, \overline{(\delta_s\tens u)^*}
\,=\, \overline{(\delta_{s\ra u}\tens u^{-1}\tau(a,a^R))\bar\ra \tau(a,a^R)^{-1}}\ .
\end{eqnarray}
By construction this preserves grade
(since $(a^L)^R=a$), and now we show that it also
preserves the $G$-action, and is thus a morphism in ${\CC}_{\bowtie}$. 
\begin{eqnarray*}
*((\delta_s\tens u)\bar\ra v) &=& *(\delta_{s\ra(a\la v)}\tens (a\la v)^{-1}uv) \cr
&=& \overline{(\delta_{s\ra uv}\tens v^{-1}u^{-1}(a\la v)\tau(a\ra v,(a\ra v)^R))
\bar\ra \tau(a\ra v,(a\ra v)^R)^{-1}} \cr
&=& 
\overline{(\delta_{s\ra uv}\tens v^{-1}u^{-1}av(a\ra v)^R)
\bar\ra \tau(a\ra v,(a\ra v)^R)^{-1}}
\ ,\cr
(*(\delta_s\tens u))\bar\ra v &=& 
(\overline{(\delta_{s\ra u}\tens u^{-1}\tau(a,a^R))\bar\ra \tau(a,a^R)^{-1}})\bar\ra v \cr
&=& \overline{(\delta_{s\ra u}\tens u^{-1}\tau(a,a^R))\bar\ra \tau(a,a^R)^{-1}(a\la v)}\ .
\end{eqnarray*}
To show that these are the same, we first calculate
\begin{eqnarray*}
 \tau(a\ra v,(a\ra v)^R)^{-1}(a\la v)^{-1}\tau(a,a^R) 
&=& (a\ra v)^{R-1}(a\ra v)^{-1}(a\la v)^{-1}a a^R \cr
&=& (a\ra v)^{R-1}v^{-1} a^R\ ,
\end{eqnarray*}
so we have to show the following equality:
\begin{eqnarray} \label{kucvbcb}
\delta_{s\ra uv}\tens v^{-1}u^{-1}av(a\ra v)^R
\,=\, (\delta_{s\ra u}\tens u^{-1}\tau(a,a^R))\bar \ra
a^{R-1}v(a\ra v)^{R}\ .
\end{eqnarray}
Now we use the factorisation (remembering that $a^{R-1}v(a\ra v)^{R}\in G$)
\begin{eqnarray*}
a^Ra^{R-1}v(a\ra v)^{R}\,=\,v(a\ra v)^{R}\,=\,(a^R\la a^{R-1}v(a\ra v)^{R})
(a^R\ra a^{R-1}v(a\ra v)^{R})
\end{eqnarray*}
to apply the formula for the action on $A_{\bowtie}$,
\begin{eqnarray*}
 (\delta_{s\ra u}\tens u^{-1}\tau(a,a^R))\bar \ra
a^{R-1}v(a\ra v)^{R} &=& \delta_{s\ra uv}\tens v^{-1}u^{-1}\tau(a,a^R)a^{R-1}v(a\ra v)^{R}\ .
\end{eqnarray*}
This verifies (\ref{kucvbcb}), and completes the proof.\end{proof}

\begin{propos} 
$A_{\bowtie}$ is a star algebra in the category ${\CC}_{\bowtie}$.
\end{propos}
\begin{proof}
Set $x=\delta_s\tens u$ and $y=\delta_t\tens v$, and set $\<x\>=a$
and $\<y\>=b$. Now
\begin{eqnarray} \label{kjdhzdsv}
\overline{\mu}\,\Upsilon^{-1}(\overline{x^*}\tens \overline{y^*}) &=& 
\overline{\mu(y^*\bar\ra \tau(a^L,a)^{-1}\tens x^*) \bar\ra \tau(a,b)}\cr
&=& \overline{((\delta_{t\ra v\tau(b^L,b)^{-1}}\tens \tau(b^L,b)v^{-1})
\bar\ra \tau(a^L,a)^{-1})} \cr
&& \overline{(\delta_{s\ra u\tau(a^L,a)^{-1}}\tens \tau(a^L,a)u^{-1})
\bar\ra \tau(a,b)} \cr
&=& \overline{(\delta_{t\ra v\tau(b^L,b)^{-1}
(b^L\la \tau(a^L,a)^{-1})}\tens  (b^L\la \tau(a^L,a)^{-1})^{-1}
\tau(b^L,b)v^{-1} \tau(a^L,a)^{-1})} \cr
&& \overline{(\delta_{s\ra u\tau(a^L,a)^{-1}}\tens \tau(a^L,a)u^{-1})
\bar\ra \tau(a,b)}\ .
\end{eqnarray}
The grade of the first factor in this product we label as $c=b^L\ra \tau(a^L,a)^{-1}$.
Then we continue with the product:
\begin{eqnarray*}
\overline{\mu}\Upsilon^{-1}(\overline{x^*}\tens \overline{y^*}) &=&
 \delta_{s\ra u\tau(a^L,a)^{-1},t\ra \tau(a^L,a)^{-1}} \ \overline{
 (\delta_{t\ra v\tau(b^L,b)^{-1}
(b^L\la \tau(a^L,a)^{-1})\tau(c,a^L)}\tens  } \cr
&& \overline{\tau(c,a^L)^{-1}(b^L\la \tau(a^L,a)^{-1})^{-1}
\tau(b^L,b)v^{-1}u^{-1})\bar\ra \tau(a,b)}\cr
&=&
 \delta_{s\ra u,t} \ \overline{
 \delta_{t\ra v\tau(b^L,b)^{-1}
(b^L\la \tau(a^L,a)^{-1})\tau(c,a^L)w}\tens  } \cr
&& \overline{w^{-1}\tau(c,a^L)^{-1}(b^L\la \tau(a^L,a)^{-1})^{-1}
\tau(b^L,b)v^{-1}u^{-1} \tau(a,b)}\ ,
\end{eqnarray*}
where we have used $w=(c\!\cdot\! a^L)\la \tau(a,b)$. Next
\begin{eqnarray*}
&& \tau(b^L,b)^{-1}(b^L\la \tau(a^L,a)^{-1})\tau(c,a^L)w((c\!\cdot\! a^L)\ra 
\tau(a,b))\cr &=& 
\tau(b^L,b)^{-1}(b^L\la \tau(a^L,a)^{-1})c a^L\tau(a,b) \cr
&=& \tau(b^L,b)^{-1}b^L \tau(a^L,a)^{-1} a^L\tau(a,b) \cr
&=& b^{-1}a^{-1} \tau(a,b)\cr
&=& (\tau(a,b)(a\!\cdot\! b))^{-1}\tau(a,b)\cr
&=& (a\!\cdot\! b)^{-1} \cr
&=& \tau((a\!\cdot\! b)^L,a\!\cdot\! b)^{-1}\, (a\!\cdot\! b)^L\ .
\end{eqnarray*}
This means that
\begin{eqnarray*}
\overline{\mu}\Upsilon^{-1}(\overline{x^*}\tens \overline{y^*})
&=&
 \delta_{s\ra u,t} \ \overline{
 \delta_{t\ra v \tau((a\cdot b)^L,a\cdot b)^{-1}}\tens 
  \tau((a\!\cdot\! b)^L,a\!\cdot\! b)v^{-1}u^{-1} \tau(a,b)}\cr
&=&
 \delta_{s\ra u,t} \ \overline{
 \delta_{s\ra uv \tau((a\cdot b)^L,a\cdot b)^{-1}}\tens 
\tau((a\!\cdot\! b)^L,a\!\cdot\! b)v^{-1}u^{-1} \tau(a,b)}\ .
\end{eqnarray*}
We compare this with
\begin{eqnarray*}
*\mu(x\tens y ) &=&  \delta_{s\ra u,t} \ *(\delta_{s\ra \tau(a,b)}
\tens \tau(a,b)^{-1}uv) \ ,
\end{eqnarray*}
which gives exactly the same result. \end{proof}

\medskip
There is also an object $\mathbb{C}(M)$ with basis elements $\delta_s$
for all $s\in M$, and grade and action given by $\<\delta_s\>=s$
and $\delta_s\ra v=\delta_{s\ra v}$. This has a star operation
defined by $\star\delta_s=\overline{\delta_{s^L}}$. 
To check that this star is a morphism, we have
\begin{eqnarray*}
\overline{\delta_{s^L}} \bar\ra u \,=\, \overline{\delta_{s^L}\bar\ra(s\la u)} 
\,=\, \overline{\delta_{s^L\ra(s\la u)}} \,=\,\overline{\delta_{(s\ra u)^L}}\ .
\end{eqnarray*}

\begin{propos}
There is a nondegenerate inner product (see \ref{inndeff}) given by a morphism
$\<,\>_r:\mathbb{C}(M)\tens\overline{\mathbb{C}(M)}\to 1_{\mathbb{C}}$,
defined by  $\<\delta_t, \overline{\delta_s}\>_r=
\delta_{s, t}$. It is also symmetric.
\end{propos}
\begin{proof} The definition is designed to preserve grades. For the actions,
\begin{eqnarray*}
\<\delta_t\ra(\<\overline{\delta_s}\>\la u),\overline{\delta_s}\ra u\>_r \,=\,
\<\delta_{t\ra(s^R\la u)},\overline{\delta_{s\ra(s^R\la u)}}\>_r\,=\,
\delta_{t\ra(s^R\la u),s\ra(s^R\la u)}\,=\,\delta_{t,s}\ .
\end{eqnarray*}
To check the symmetry of the inner product on $\mathbb{C}(M)$:
\begin{eqnarray*}
\Upsilon^{-1}(\mathrm{bb}(\delta_s)\tens\overline{\delta_t}) &=& 
\Upsilon^{-1}(\overline{\overline{\delta_{s\ra\tau(s^L,s)^{-1}}}}\tens\overline{\delta_t})
\cr
&=& \overline{(\delta_{t\ra\tau(p,p^R)^{-1}}\tens 
\overline{\delta_{s\ra\tau(s^L,s)^{-1}}})\ra\tau(p^R,t^R)}
\end{eqnarray*}
where $p=\<\overline{\delta_{s\ra\tau(s^L,s)^{-1}}}\>=(s\ra\tau(s^L,s)^{-1})^R=s^L$. 
It follows that 
\begin{eqnarray*}
\overline{\<,\>_r}\,\Upsilon^{-1}(\mathrm{bb}(\delta_s)\tens\overline{\delta_t})
&=& \overline{\<\delta_{t\ra\tau(p,p^R)^{-1}},
\overline{\delta_{s\ra\tau(s^L,s)^{-1}}})\>_r} \cr
&=& \overline{\delta_{t\ra\tau(s^L,s)^{-1},  s\ra\tau(s^L,s)^{-1}  }} \cr
&=& \overline{\delta_{t,s}}\ . 
\end{eqnarray*}\end{proof}

\medskip
For any object $V$ in ${\CC}_{\bowtie}$, there is a right action
of $A_{\bowtie}$ on $V$ given by 
\begin{eqnarray*}
v\ra (\delta_s\tens u) \,=\, \delta_{\<v\>,s}\, v\ra u\ .
\end{eqnarray*}
Further $\ra:V\tens A_{\bowtie} \to V$ is a morphism in ${\CC}_{\bowtie}$.

\begin{propos}
The usual action $\ra:\mathbb{C}(M)\tens A_{\bowtie}\to \mathbb{C}(M)$
is compatible with the inner product (see \ref{dkjshbcvb}). 
\end{propos}
\begin{proof} 
Let $a=\<\delta_s\tens u\>$. 
\begin{eqnarray*}
(\bar\ra\tens\id)((\delta_p\tens(\delta_s\tens u))\tens\overline{\delta_q})
&=& \delta_{p,s}.\delta_{p\ra u}\tens \overline{\delta_q}\ ,\cr
(\id\tens(*\tens\id))\Phi((\delta_p\tens(\delta_s\tens u))\tens\overline{\delta_q})&=&
\delta_{p\ra\tau(a,q^R)} \tens( *(\delta_s\tens u)\tens\overline{\delta_q}) \cr
&=& \delta_{p\ra\tau(a,q^R)} \tens( (\overline{\delta_{s'}\tens u'})
\tens\overline{\delta_q})\ ,
\end{eqnarray*}
where $s'=s\ra u\tau(a^L,a)^{-1}$ and $u'=\tau(a^L,a)u^{-1}$. Following on from this,
\begin{eqnarray*}
\delta_{p\ra\tau(a,q^R)} \tens \overline{\ra}\Upsilon^{-1}( (\overline{\delta_{s'}\tens u'})
\tens\overline{\delta_q})
&=& \delta_{p\ra\tau(a,q^R)} \tens  \cr
&&\overline{\ra}(\overline{(\delta_{q\ra\tau(a^L,a)^{-1}}\tens(\delta_{s'}\tens u'))
\ra\tau(a,q^R)}) \cr
&=& \delta_{p\ra\tau(a,q^R)} \tens  \cr
&&(\overline{(\delta_{q\ra\tau(a^L,a)^{-1}}\ra(\delta_{s'}\tens u'))
\ra\tau(a,q^R)}) \cr
&=& \delta_{s',q\ra\tau(a^L,a)^{-1}}.  \delta_{p\ra\tau(a,q^R)} \tens \cr
&& \overline{\delta_{q\ra\tau(a^L,a)^{-1}u'}
\ra\tau(a,q^R)} \cr
&=& \delta_{s\ra u,q}.  \delta_{p\ra\tau(a,q^R)} \tens 
\overline{\delta_{q\ra u^{-1}\tau(a,q^R)}}
\end{eqnarray*}
Now apply the inner product.\end{proof}

\section{Duals}

\subsection{Rigid tensor categories}

The left dual is a contravariant 
functor from ${\CC}$ to ${\CC}$,
written as $X\mapsto X'$. For every object $X$ in ${\CC}$
there are morphisms $\ev^L_X:X'\tens X\to 1_{{\CC}}$
(evaluation) and 
$\coev^L_X:1_{{\CC}}\to X\tens X'$ (coevaluation) so that 
\begin{eqnarray*}
l^{-1}_X(\id\tens\ev_X)\Phi(\coev_X\tens \id)r_X  &=& \id_X:X\to X\ , \cr
r_{X'}^{-1}(\ev_X\tens\id)\Phi^{-1}(\id\tens\coev{X})l_{X'}  &=& \id_{X'}:X'\to X'\ .
\end{eqnarray*}
Further $\ev^L$ and $\coev^L$ obey the following properties for
every morphism $\phi:X\to Y$ (we are prevented from making a simpler statement
about being a natural transformation because of the contravariant nature of
dual).

\begin{picture}(100,80)(-100,6)

\put(60,63){\vector(1,0){45}}
\put(60,23){\vector(1,0){45}}
\put(38,55){\vector(0,-1){20}}
\put(127,55){\vector(0,-1){20}}
\put(20,60){$Y'\tens X$}
\put(20,20){$X'\tens X$}
\put(111,60){$Y'\tens Y$}
\put(121,20){$1_{\CC}$}

\put(65,70){$\id\tens\phi$}
\put(67,30){$\ev^L_X$}
\put(-1,43){$\phi'\tens\id$}
\put(133,43){$\ev^L_Y$}
\end{picture}

\begin{picture}(100,80)(-100,6)

\put(60,63){\vector(1,0){45}}
\put(60,23){\vector(1,0){45}}
\put(38,55){\vector(0,-1){20}}
\put(127,55){\vector(0,-1){20}}
\put(30,60){$1_{\CC}$}
\put(20,20){$X\tens X'$}
\put(111,60){$Y\tens Y'$}
\put(111,20){$Y\tens X'$}

\put(65,70){$\coev^L_Y$}
\put(67,30){$\phi\tens\id$}
\put(-1,43){$\coev^L_X$}
\put(133,43){$\id\tens\phi'$}
\end{picture}

We could also have a right dual, which would be a contravariant functor
 from ${\CC}$ to ${\CC}$,
written as $X\mapsto  X^\circ$. For every object $X$ in ${\CC}$
there are morphisms $\ev^R_X:X\tens X^\circ\to 1_{{\CC}}$
(evaluation) obeying the opposite versions of the rules for the left dual.

\begin{defin}
A tensor category is called left rigid if every object has a left dual, and right rigid
if every object has a right dual. 
\end{defin}

\subsection{Making right duals using bar} 

\begin{propos}
Suppose that we have a left rigid category which is also a bar category. Then
\begin{eqnarray*}
X^\circ &=& \overline{(\overline{X})'}\ ,\cr
\ev^R &=& \star^{-1}\, \overline{\ev^L_{\bar X}}\, \Upsilon^{-1}\,
(\mathrm{bb}\tens \id) : X\tens X^\circ \to X\ ,\cr
\coev^R &=& (\id\tens \mathrm{bb}^{-1})\, \Upsilon\, \overline{\coev^L_{\bar X}}\, \star :
1 \to X^\circ\tens X\ .
\end{eqnarray*}
makes the bar category right rigid.
\end{propos}
\begin{proof} We have to see if these obey the required properties for right rigidity,
beginning with the following composition being the identity:
\begin{eqnarray} \label{skhavchuas}
X \stackrel{\id\tens\coev_X^R}\rightarrow 
X\tens (X^\circ\tens X)    \stackrel{\Phi^{-1}}\rightarrow 
(X\tens X^\circ)\tens X    \stackrel{\ev^R_X}\rightarrow  X\ .
\end{eqnarray}
If we substitute in the various definitions, we have
\begin{eqnarray*}
&& X \stackrel{\mathrm{bb}}\longrightarrow \overline{\overline{X}}
 \stackrel{l} \longrightarrow\overline{\overline{X}} \tens 1 
  \stackrel{\id\tens \star}\longrightarrow \overline{\overline{X}} \tens \overline{1} 
\stackrel{\id\tens \coev^L_{\bar X}}\longrightarrow 
\overline{\overline{X}} \tens 
\overline{\overline{X} \tens (\overline{X} )'} 
\stackrel{\id\tens \Upsilon}\longrightarrow 
\overline{\overline{X}} \tens (
\overline{(\overline{X} )'} \tens \overline{\overline{X} } ) \cr
&&
\stackrel{\Phi^{-1}}\longrightarrow 
(\overline{\overline{X}} \tens 
\overline{(\overline{X} )'}) \tens \overline{\overline{X} } 
\stackrel{\Upsilon^{-1}\tens\id}\longrightarrow 
\overline{(\overline{X} )' \tens \overline{X} }
\tens \overline{\overline{X} } 
\stackrel{\star^{-1}\,\overline{\ev^L_{\bar X}}\tens\id}\longrightarrow 
1\tens \overline{\overline{X} }  
\stackrel{r^{-1}}\longrightarrow  \overline{\overline{X} }  
\stackrel{\mathrm{bb}^{-1}}\longrightarrow X\ .
\end{eqnarray*}
Now we use the left and right identity properties in the definition of bar
 category to rewrite this as
\begin{eqnarray*}
&& X \stackrel{\mathrm{bb}}\longrightarrow \overline{\overline{X}}
 \stackrel{\overline{r}} \longrightarrow
 \overline{1\tens\overline{X}} 
  \stackrel{\overline{\coev^L_{\bar X}\tens\id}} \longrightarrow
 \overline{(\overline{X}\tens (\overline{X})')\tens \overline{X}} 
  \stackrel{(\id\tens\Upsilon)\,\Upsilon}\longrightarrow 
\overline{\overline{X}} \tens (
\overline{(\overline{X} )'} \tens \overline{\overline{X} } ) \cr
&&
\stackrel{\Phi^{-1}}\longrightarrow 
(\overline{\overline{X}} \tens 
\overline{(\overline{X} )'}) \tens \overline{\overline{X} } 
\stackrel{\Upsilon^{-1}\,(\Upsilon^{-1}\tens\id)}\longrightarrow 
 \overline{\overline{X}\tens ((\overline{X})'\tens \overline{X})} 
\stackrel{\overline{\id\tens\ev^L_{\bar X}}}\longrightarrow 
\overline{\overline{X} \tens 1 }  
\stackrel{\overline{l^{-1}}}\longrightarrow  \overline{\overline{X} }  
\stackrel{\mathrm{bb}^{-1}}\longrightarrow X\ .
\end{eqnarray*}
Now we can use another rule in the definition of bar category to simplify this to
\begin{eqnarray*}
&& X \stackrel{\mathrm{bb}}\longrightarrow \overline{\overline{X}}
 \stackrel{\overline{r}} \longrightarrow
 \overline{1\tens\overline{X}} 
  \stackrel{\overline{\coev^L_{\bar X}\tens\id}} \longrightarrow
 \overline{(\overline{X}\tens (\overline{X})')\tens \overline{X}} 
  \stackrel{\overline{\Phi}}\longrightarrow 
 \overline{\overline{X}\tens ((\overline{X})'\tens \overline{X})} 
\stackrel{\overline{\id\tens\ev^L_{\bar X}}}\longrightarrow 
\overline{\overline{X} \tens 1 }  \cr
&&
\stackrel{\overline{l^{-1}}}\longrightarrow  \overline{\overline{X} }  
\stackrel{\mathrm{bb}^{-1}}\longrightarrow X\ ,
\end{eqnarray*}
and by right rigidity this simplifies to the identity on $X$. 

The other required result, that the following composition is the identity,
is left to the reader: 
\begin{eqnarray*}
X^\circ \stackrel{\coev^R_{X}\tens\id}\longrightarrow
(X^\circ\tens X)\tens X^\circ \stackrel{\Phi}\longrightarrow
X^\circ\tens (X\tens X^\circ)  \stackrel{\id\tens\ev^R_X}\longrightarrow
X^\circ\ .
\end{eqnarray*} 
\end{proof}

\begin{remark} Given left and right rigid structures, we can define a morphism
$X\to X^{\circ\prime}$ as follows:
\begin{eqnarray*} 
X \stackrel{l}\longrightarrow 
X \tens 1 \stackrel{\id\tens\coev^L_{X^\circ}}\longrightarrow 
X \tens (X^\circ \tens X^{\circ \prime}) 
\stackrel{\Phi^{-1}}\longrightarrow 
(X \tens X^\circ) \tens X^{\circ \prime}
\stackrel{\ev^R \tens\id}\longrightarrow
1\tens X^{\circ \prime}
\stackrel{r^{-1}}\longrightarrow X^{\circ \prime}\ ,
\end{eqnarray*}
and this morphism has inverse given by
\begin{eqnarray*} 
X^{\circ \prime} \stackrel{l}\longrightarrow 
X^{\circ \prime} \tens 1 \stackrel{\id\tens\coev^R_{X}}\longrightarrow 
X^{\circ \prime} \tens (X^\circ \tens X) 
\stackrel{\Phi^{-1}}\longrightarrow 
(X^{\circ \prime} \tens X^\circ) \tens X
\stackrel{\ev^L \tens\id}\longrightarrow
1\tens X
\stackrel{r^{-1}}\longrightarrow X\ .
\end{eqnarray*}
\end{remark}

\appendix

\section{Reconstruction theorem}

The reconstruction theory involves a representable functor
$F$ from a strong bar category to ${\bf Vect}_{\mathbb{C}}$ (with its usual
bar structure and transposition). There are two versions of this
theory, one which reconstructs a Hopf algebra via its comodules and which is technically
superior, and the other via its modules which is less useful but easier for most
readers. In keeping with the line taken in the paper we focus on the module case but
the other has identical proofs with arrows reversed as for example in
 \cite{majquasitann, MajBook}. 
 
Representable here means that there is a vector space $H$ with linear maps
$\mathrm{Lin}(V,H)$ in natural 1-1 correspendence
with natural transformations between the functors
$V\tens F$ and $F$
for every vector space $V$. [The functor $V\tens F$ maps $X$ to $V\tens F(X)$.] 
We begin with $\alpha_X:H\tens F(X)\to F(X)$, which corresponds to the
identity $\id\in \mathrm{Lin}(H,H)$, which may be more recognisable after
reconstruction as the action of $H$ on $F(X)$. Thus we also write
$\alpha_X$ as $h\tens x \mapsto h\la x$.
Using the bar structure,
we have the following natural transformation $\beta$ between
$\overline{H}\tens F$ and $F$: 
\begin{eqnarray*}
\overline{H}\tens F(X) \stackrel{ \id\tens\mathrm{bb} } \longrightarrow 
\overline{H}\tens \overline{\overline{ F(X)}}
\stackrel{ \overline{\Psi}\,\Upsilon^{-1} } \longrightarrow 
\overline{H\tens\overline{ F(X)}}
\stackrel{ \overline{\tilde\alpha_X }} \longrightarrow 
\overline{\overline{ F(X)}}
\stackrel{ \mathrm{bb}^{-1}} \longrightarrow F(X)\ ,
\end{eqnarray*}
where $\tilde\alpha_X $ is given by the composition
\begin{eqnarray*}
H\tens\overline{ F(X)} 
\stackrel{ \id\tens\mathrm{fb}_X} \longrightarrow 
H\tens F(\overline{ X})  
\stackrel{ \alpha_{\bar X}} \longrightarrow F(\overline{ X}) 
\stackrel{ \mathrm{fb}_X^{-1}} \longrightarrow \overline{ F(X)}\ .
\end{eqnarray*}
(The reader should note that we have only used the braiding in 
${\bf Vect}_{\mathbb{C}}$, i.e.\ transposition.)
The natural transformation $\beta$ will correspond to a map
$T:\overline{H}\to H$ so that:

\begin{picture}(100,80)(-100,6)

\put(60,63){\vector(1,0){45}}
\put(60,23){\vector(1,0){45}}
\put(38,55){\vector(0,-1){20}}
\put(127,55){\vector(0,-1){20}}
\put(10,60){$\overline{H}\tens F(X)$}
\put(10,20){$ H\tens F(X)$}
\put(111,60){$F(X)$}
\put(111,20){$F(X)$}

\put(75,70){$\beta_X$}
\put(75,30){$\alpha_X$}
\put(-1,43){$T\tens\id$}
\put(133,43){$\id$}


\end{picture}

As we are operating in ${\bf Vect}$ where both the braiding $\Psi$ and
$\Upsilon$ are transpositions (with appropriate bars in the case of $\Upsilon$),
\begin{eqnarray*}
\tilde\alpha_X(h\tens\overline{x})\,=\, \overline{T(\overline{h})\la x}\ ,
\end{eqnarray*}
and
the definition of $T$ simplifies to $T(\overline{h})\la$ being the following
 composition:
\begin{eqnarray} \label{recobasic}
F(X) \stackrel{\mathrm{bb}} \longrightarrow 
\overline{\overline{F(X)}} \stackrel{\overline{\mathrm{fb}}} \longrightarrow 
\overline{F(\overline{X})} \stackrel{\overline{h\la}} \longrightarrow 
\overline{F(\overline{X})} \stackrel{\overline{\mathrm{fb}^{-1}}}
 \longrightarrow 
 \overline{\overline{F(X)}}  \stackrel{\mathrm{bb}^{-1}} \longrightarrow F(X)\ .
\end{eqnarray}

\begin{propos} \label{propiii1}
$T(hh')\la=T(h)T(h')\la$.
\end{propos}
\begin{proof} By applying (\ref{recobasic}) twice.\end{proof}

\begin{propos}  \label{propiii2}
$T(\overline{T(\overline{h})}\la=h\la$.
\end{propos}
\begin{proof}
The following composition gives $T(\overline{T(\overline{h})}\la$:
\begin{eqnarray*}
F(X) \stackrel{\mathrm{bb}} \longrightarrow 
\overline{\overline{F(X)}} \stackrel{\overline{\mathrm{fb}}} \longrightarrow 
\overline{F(\overline{X})} \stackrel{\overline{T(\overline{h})\la}} \longrightarrow 
\overline{F(\overline{X})} \stackrel{\overline{\mathrm{fb}^{-1}}}
 \longrightarrow 
 \overline{\overline{F(X)}}  \stackrel{\mathrm{bb}^{-1}} \longrightarrow F(X)\ ,
\end{eqnarray*}
and if we use the definition of $T$ again we get
\begin{eqnarray*}
&& F(X) \stackrel{\mathrm{bb}} \longrightarrow 
\overline{\overline{F(X)}} \stackrel{\overline{\mathrm{fb}}} \longrightarrow 
\overline{F(\overline{X})} \stackrel{\overline{\mathrm{bb}}} \longrightarrow 
\overline{\overline{\overline{F(\overline{X})} }} 
\stackrel{\overline{\overline{\mathrm{fb}}}} \longrightarrow 
\overline{\overline{F(\overline{\overline{X}})}} 
\stackrel{\overline{\overline{h\la}}} \longrightarrow 
\overline{\overline{F(\overline{\overline{X}})}}
\stackrel{\overline{\overline{\mathrm{fb}^{-1}}}} \longrightarrow 
\overline{\overline{\overline{F(\overline{X})} }} 
\stackrel{\overline{\mathrm{bb}^{-1}}} \longrightarrow  
\overline{F(\overline{X})} \cr
&& \stackrel{\overline{\mathrm{fb}^{-1}}}
 \longrightarrow 
 \overline{\overline{F(X)}}  \stackrel{\mathrm{bb}^{-1}} \longrightarrow F(X)\ .
\end{eqnarray*}
By using Proposition \ref{kxzj}, this becomes
\begin{eqnarray*}
 F(X) \stackrel{\mathrm{bb}} \longrightarrow 
\overline{\overline{F(X)}} 
\stackrel{\overline{\overline{F(\mathrm{bb})}}} \longrightarrow 
\overline{\overline{F(\overline{\overline{X}})}} 
\stackrel{\overline{\overline{h\la}}} \longrightarrow 
\overline{\overline{F(\overline{\overline{X}})}}
\stackrel{\overline{\overline{F(\mathrm{bb}^{-1})}}} \longrightarrow 
 \overline{\overline{F(X)}}  \stackrel{\mathrm{bb}^{-1}} \longrightarrow F(X)\ ,
\end{eqnarray*}
and by the functorial property of $h\la$ (as $\alpha$ is functorial), we get
\begin{eqnarray*}
 F(X) \stackrel{\mathrm{bb}} \longrightarrow 
\overline{\overline{F(X)}} 
\stackrel{\overline{\overline{h\la}}} \longrightarrow 
\overline{\overline{F(X)}} 
  \stackrel{\mathrm{bb}^{-1}} \longrightarrow F(X)\ ,
\end{eqnarray*}
and as $\mathrm{bb}$ is a natural transformation, this is just
\begin{eqnarray*}
 F(X) 
\stackrel{h\la} \longrightarrow F(X)\ .
\end{eqnarray*}\end{proof}

\medskip

The coproduct is reconstructed from the following diagram,
where we have used $\Delta h\la(x\tens y)=h_{(1)}\la x \tens h_{(2)}\la y$:

\begin{picture}(100,80)(-100,6)

\put(60,63){\vector(1,0){45}}
\put(70,23){\vector(1,0){35}}
\put(38,55){\vector(0,-1){20}}
\put(127,55){\vector(0,-1){20}}
\put(10,60){$F(X\tens Y)$}
\put(0,20){$ F(X)\tens F(Y)$}
\put(111,60){$F(X\tens Y)$}
\put(111,20){$F(X)\tens F(Y)$}

\put(75,70){$h\la$}
\put(72,30){$\Delta h\la$}
\put(-1,43){$f^{-1}$}
\put(133,43){$f^{-1}$}

\end{picture}

\begin{propos}  \label{propiii3}
$T(\overline{h})_{(1)}\la \tens T(\overline{h})_{(2)}\la =
T(\overline{h_{(2)}})\la\tens T(\overline{h_{(1)}})\la$.
\end{propos}
\begin{proof} 
The following composition
gives $T(\overline{h})\la$ on $F(X\tens Y)$:
\begin{eqnarray*}
F(X\tens Y) \stackrel{\mathrm{bb}} \longrightarrow 
\overline{\overline{F(X\tens Y)}} \stackrel{\overline{\mathrm{fb}}} \longrightarrow 
\overline{F(\overline{X\tens Y})} \stackrel{\overline{h\la}} \longrightarrow 
\overline{F(\overline{X\tens Y})} \stackrel{\overline{\mathrm{fb}^{-1}}}
 \longrightarrow 
 \overline{\overline{F(X\tens Y)}}  \stackrel{\mathrm{bb}^{-1}} \longrightarrow
  F(X\tens Y)\ .
\end{eqnarray*}
By functoriality this is the same as
\begin{eqnarray*}
&& F(X\tens Y) \stackrel{\mathrm{bb}} \longrightarrow 
\overline{\overline{F(X\tens Y)}} \stackrel{\overline{\mathrm{fb}}} \longrightarrow 
\overline{F(\overline{X\tens Y})}  \stackrel{\overline{F(\Upsilon)}} \longrightarrow 
\overline{F(\overline{Y}\tens \overline{X})} 
 \stackrel{\overline{f^{-1}}} \longrightarrow 
 \overline{F(\overline{Y})\tens F(\overline{X})} 
\stackrel{\overline{\Delta h\la}} \longrightarrow \cr
&&
 \overline{F(\overline{Y})\tens F(\overline{X})} 
  \stackrel{\overline{f}} \longrightarrow 
  \overline{F(\overline{Y}\tens \overline{X})} 
   \stackrel{\overline{F(\Upsilon^{-1})}} \longrightarrow 
\overline{F(\overline{X\tens Y})} \stackrel{\overline{\mathrm{fb}^{-1}}}
 \longrightarrow 
 \overline{\overline{F(X\tens Y)}}  \stackrel{\mathrm{bb}^{-1}} \longrightarrow
  F(X\tens Y)\ .
\end{eqnarray*}
Now we use the definition of $T$ to get
\begin{eqnarray*}
&& F(X\tens Y) \stackrel{\mathrm{bb}} \longrightarrow 
\overline{\overline{F(X\tens Y)}} \stackrel{\overline{\mathrm{fb}}} \longrightarrow 
\overline{F(\overline{X\tens Y})}  \stackrel{\overline{F(\Upsilon)}} \longrightarrow 
\overline{F(\overline{Y}\tens \overline{X})} 
 \stackrel{\overline{f^{-1}}} \longrightarrow 
 \overline{F(\overline{Y})\tens F(\overline{X})} 
\stackrel{\overline{\mathrm{fb}^{-1}\tens \mathrm{fb}^{-1}}} \longrightarrow \cr
&&
 \overline{\overline{F(Y)}\tens \overline{F(X)}} 
 \stackrel{\overline{\overline{T(\overline{h_{(1)}})\la}\tens  \overline{T(\overline{h_{(2)}})\la}  }} 
 \longrightarrow 
  \overline{\overline{F(Y)}\tens \overline{F(X)}} 
  \stackrel{\overline{\mathrm{fb}\tens \mathrm{fb}}} \longrightarrow 
 \overline{F(\overline{Y})\tens F(\overline{X})} \cr
 &&
  \stackrel{\overline{f}} \longrightarrow 
  \overline{F(\overline{Y}\tens \overline{X})} 
   \stackrel{\overline{F(\Upsilon^{-1})}} \longrightarrow 
\overline{F(\overline{X\tens Y})} \stackrel{\overline{\mathrm{fb}^{-1}}}
 \longrightarrow 
 \overline{\overline{F(X\tens Y)}}  \stackrel{\mathrm{bb}^{-1}} \longrightarrow
  F(X\tens Y)\ .
\end{eqnarray*}
By property (4) in Definition~\ref{recdef} this is the same as
\begin{eqnarray*}
&& F(X\tens Y) \stackrel{\mathrm{bb}} \longrightarrow 
\overline{\overline{F(X\tens Y)}} \stackrel{\overline{\overline{f^{-1}}}} \longrightarrow 
\overline{\overline{F(X)\tens F( Y)}} 
\stackrel{\overline{\Upsilon}} \longrightarrow 
\overline{\overline{F(Y)}\tens \overline{F(X)}} 
 \cr
&&
 \stackrel{\overline{\overline{T(\overline{h_{(1)}})\la}\tens  \overline{T(\overline{h_{(2)}})\la}  }} 
 \longrightarrow 
  \overline{\overline{F(Y)}\tens \overline{F(X)}} 
   \stackrel{\overline{\Upsilon^{-1}}} \longrightarrow 
\overline{\overline{F(X)\tens F( Y)}}  \stackrel{\overline{\overline{f}}}
 \longrightarrow 
 \overline{\overline{F(X\tens Y)}}  \stackrel{\mathrm{bb}^{-1}} \longrightarrow
  F(X\tens Y)\ .
\end{eqnarray*}
Using the functorial property of $\Upsilon$ gives
\begin{eqnarray*}
&& F(X\tens Y) \stackrel{\mathrm{bb}} \longrightarrow 
\overline{\overline{F(X\tens Y)}} \stackrel{\overline{\overline{f^{-1}}}} \longrightarrow 
\overline{\overline{F(X)\tens F( Y)}} 
 \cr
&&
 \stackrel{\overline{\overline{T(\overline{h_{(2)}})\la\tens T(\overline{h_{(1)}})\la}  }} 
 \longrightarrow 
\overline{\overline{F(X)\tens F( Y)}}  \stackrel{\overline{\overline{f}}}
 \longrightarrow 
 \overline{\overline{F(X\tens Y)}}  \stackrel{\mathrm{bb}^{-1}} \longrightarrow
  F(X\tens Y)\ .
\end{eqnarray*}
Now as $\mathrm{bb}$ is a natural transformation we get
\begin{eqnarray*}
&& F(X\tens Y) \stackrel{f^{-1}} \longrightarrow 
F(X)\tens F(Y) \stackrel{\mathrm{bb}} \longrightarrow 
\overline{\overline{F(X)\tens F( Y)}} 
 \cr
&&
 \stackrel{\overline{\overline{T(\overline{h_{(2)}})\la\tens T(\overline{h_{(1)}})\la}  }} 
 \longrightarrow 
\overline{\overline{F(X)\tens F( Y)}}  \stackrel{\mathrm{bb}^{-1}}
 \longrightarrow 
F(X)\tens F(Y)  \stackrel{f} \longrightarrow
  F(X\tens Y)\ ,
\end{eqnarray*}
which, for the same reason, becomes
\begin{eqnarray*}
&& F(X\tens Y) \stackrel{f^{-1}} \longrightarrow 
F(X)\tens F(Y) \stackrel{T(\overline{h_{(2)}})\la\tens T(\overline{h_{(1)}})\la} \longrightarrow 
F(X)\tens F(Y)  \stackrel{f} \longrightarrow
  F(X\tens Y)\ .
\end{eqnarray*}\end{proof}

\medskip The antipode is reconstructed using the left dual
structure from the following diagram for $S(h)\la$:
\begin{eqnarray}\label{sdefppp}
&&F(X) \stackrel{F(\mathrm{coev}^L)\tens\id}\longrightarrow 
F(X\tens X')\tens F(X) \stackrel{f^{-1}\tens\id}\longrightarrow 
F(X)\tens F(X')\tens F(X)  \stackrel{\id\tens h\la\tens\id}\longrightarrow \nonumber\\
&& F(X) \tens F(X')\tens F(X)   \stackrel{\id\tens f}\longrightarrow 
F(X)\tens F(X'\tens X)  \stackrel{\id\tens F(\mathrm{ev}^L)}\longrightarrow 
F(X)\ .\nonumber\\  & \  &
\end{eqnarray} 
{}From the right dual structure we also have the following diagram for $S^{-1}(h)\la$:
\begin{eqnarray}\label{sdefpppinv}
&& F(X) \stackrel{\id\tens F(\mathrm{coev}^R)}\longrightarrow 
F(X)\tens F(X^\circ\tens X) \stackrel{\id\tens f^{-1}}\longrightarrow 
F(X)\tens F(X^\circ)\tens F(X)  \stackrel{\id\tens h\la\tens\id}\longrightarrow \nonumber\\
&& F(X)\tens F(X^\circ)\tens F(X)   \stackrel{f\tens\id}\longrightarrow 
F(X\tens X^\circ)\tens F(X)  \stackrel{F(\mathrm{ev}^R\tens\id}\longrightarrow 
F(X)\ .\nonumber \\ &\ &
\end{eqnarray} That these really are mutually inverse can be seen either by explicit calculation,
using the functoriality of $\alpha$, or by quoting the uniqueness of antipode
in a Hopf algebra.

\begin{propos}  \label{propiii4} 
$S(T(\overline{h}))\la$ is the same as 
$p\la:F(\overline{X}) \to
F(\overline{X}) $ where $\overline{p}=T^{-1}(S^{-1}(h))$.
\end{propos}
\begin{proof} 
If we replace $X$ by $\overline{X}$ in (\ref{sdefppp}), we get (using the functoriality
of $\alpha$):
\begin{eqnarray*}
&& F(\overline{X}) \stackrel{f^{-1}\,F(\mathrm{coev})\tens\id}\longrightarrow 
F(\overline{X})\tens F((\overline{X})')\tens F(\overline{X})
\stackrel{\id\tens\mathrm{fb}^{-1}\, F(\mathrm{bb})\tens\id}\longrightarrow
F(\overline{X})\tens \overline{F(X^\circ)} \tens F(\overline{X}) \cr
&&  \stackrel{\id\tens  \overline{T(\overline{h}) \la}\tens\id}\longrightarrow 
F(\overline{X})\tens \overline{F(X^\circ)} \tens F(\overline{X}) 
\stackrel{\id\tens F(\mathrm{bb}^{-1})\,\mathrm{fb}\tens\id }\longrightarrow
F(\overline{X})\tens F((\overline{X})')\tens F(\overline{X}) 
 \stackrel{\id\tens F(\mathrm{ev})\,f}\longrightarrow 
F(\overline{X}) \ , \nonumber\\
\end{eqnarray*}
where we have used the fact that $\mathrm{bb}:(\overline{X})'\to 
\overline{X^\circ}$. 
Now, using the fact that $T(\overline{T(\overline{h})}\la=h\la$, we see that
$S(T(\overline{h}))\la$ is given by the composition
\begin{eqnarray*} 
&& F(\overline{X}) \stackrel{f^{-1}\,F(\mathrm{coev}_{\bar X})\tens\id}\longrightarrow 
F(\overline{X})\tens F((\overline{X})')\tens F(\overline{X})
\stackrel{\id\tens\mathrm{fb}^{-1}\, F(\mathrm{bb})\tens\id}\longrightarrow
F(\overline{X})\tens \overline{F(X^\circ)} \tens F(\overline{X}) \cr
&&  \stackrel{\id\tens  \overline{h \la}\tens\id}\longrightarrow 
F(\overline{X})\tens \overline{F(X^\circ)} \tens F(\overline{X}) 
\stackrel{\id\tens F(\mathrm{bb}^{-1})\,\mathrm{fb}\tens\id }\longrightarrow
F(\overline{X})\tens F((\overline{X})')\tens F(\overline{X}) 
 \stackrel{\id\tens F(\mathrm{ev}_{\bar X})\,f}\longrightarrow 
F(\overline{X})\ .
\end{eqnarray*}
We can rewrite this as 
\begin{eqnarray*} 
&& F(\overline{X}) \stackrel{r}\longrightarrow 1 \tens F(\overline{X})
\stackrel{f^{-1}\,F(\mathrm{coev}^L_{\bar X})\tens\id}\longrightarrow 
F(\overline{X})\tens F((\overline{X})')\tens F(\overline{X})
\stackrel{\mathrm{fb}^{-1}\tens\mathrm{fb}^{-1}\, F(\mathrm{bb})\tens \mathrm{fb}^{-1}}
\longrightarrow \cr
&&
\overline{F(X)}\tens \overline{F(X^\circ)} \tens \overline{F(X)} 
  \stackrel{\id\tens  \overline{h \la}\tens\id}\longrightarrow 
\overline{F(X)}\tens \overline{F(X^\circ)} \tens \overline{F(X)}
\stackrel{\mathrm{fb}\tens F(\mathrm{bb}^{-1})\,\mathrm{fb}\tens\mathrm{fb} }\longrightarrow \cr 
&&
F(\overline{X})\tens F((\overline{X})')\tens F(\overline{X}) 
 \stackrel{\id\tens F(\mathrm{ev}^L_{\bar X})\,f}\longrightarrow 
F(\overline{X}) \tens 1
\stackrel{l^{-1}}\longrightarrow
F(\overline{X})
\end{eqnarray*}
Using the commutativity of the following diagram (from property (4) in Definition~\ref{recdef})
and the definition of the right dual structures,

\begin{picture}(100,80)(-110,12)

\put(60,63){\vector(1,0){50}}
\put(60,23){\vector(1,0){80}}
\put(23,55){\vector(0,-1){20}}
\put(157,55){\vector(0,-1){20}}
\put(-17,60){$\overline{F(X^\circ)}\tens \overline{F(X)}$}
\put(-17,20){$ \overline{F(X)\tens F(X^\circ)}$}
\put(119,60){$F((\overline{X})') \tens F(\overline{X})$}
\put(155,20){$\overline{1}$}

\put(52,73){$F(\mathrm{bb}^{-1})\,\mathrm{fb}\tens \mathrm{fb}$}
\put(77,30){$\overline{F(\ev^R_X)\,f}$}
\put(-6,43){$\Upsilon^{-1}$}
\put(163,43){$\star\,F(\ev^L_{\bar X})\,f$}

\end{picture}

\noindent
this becomes
\begin{eqnarray*} 
&& F(\overline{X}) \stackrel{r}\longrightarrow 1 \tens F(\overline{X})
\stackrel{f^{-1}\,F(\mathrm{coev}^L_{\bar X})\tens\id}\longrightarrow 
F(\overline{X})\tens F((\overline{X})')\tens F(\overline{X})
\stackrel{\mathrm{fb}^{-1}\tens\mathrm{fb}^{-1}\, F(\mathrm{bb})\tens \mathrm{fb}^{-1}}
\longrightarrow \cr
&&
\overline{F(X)}\tens \overline{F(X^\circ)} \tens \overline{F(X)} 
  \stackrel{\id\tens  \overline{h \la}\tens\id}\longrightarrow 
\overline{F(X)}\tens \overline{F(X^\circ)} \tens \overline{F(X)}
\stackrel{\mathrm{fb}\tens \Upsilon^{-1}}\longrightarrow \cr 
&&
F(\overline{X})\tens \overline{F(X) \tens F(X^{\circ})}
 \stackrel{\id\tens \overline{F(\mathrm{ev}^R_{X})\,f}}\longrightarrow 
F(\overline{X}) \tens \overline{1} \stackrel{\id\tens \star^{-1}}\longrightarrow
F(\overline{X}) \tens 1
\stackrel{l^{-1}}\longrightarrow
F(\overline{X})
\end{eqnarray*}
Using the commutativity of the following diagram (from property (4) in \ref{recdef})
and the definition of the right dual structures,

\begin{picture}(100,80)(-120,6)

\put(45,63){\vector(1,0){60}}
\put(60,23){\vector(1,0){45}}
\put(28,55){\vector(0,-1){20}}
\put(127,55){\vector(0,-1){20}}
\put(25,60){$1$}
\put(-15,20){$ \overline{F(X^\circ)\tens F(X)}$}
\put(111,60){$ F(\overline{X}) \tens F((\overline{X})') $}
\put(111,20){$ \overline{F(X^\circ)\tens F(X)}$}

\put(42,70){$f^{-1}\, F(\coev^L_{\bar X})$}
\put(77,30){$\Upsilon$}
\put(-50,43){$\overline{f^{-1}\, F(\coev^R_X)}\,\star$}
\put(133,43){$\mathrm{fb}^{-1}\tens \mathrm{fb}^{-1}F(\mathrm{bb})$}


\end{picture}

\noindent so we get
\begin{eqnarray*} 
&& F(\overline{X}) \stackrel{(\id\tens\star)\,r}\longrightarrow 
\overline{1} \tens F(\overline{X})
\stackrel{\overline{f^{-1}\,F(\mathrm{coev}^R_{ X})}\tens\id}\longrightarrow 
\overline{F(X^\circ)\tens F(X)}\tens F(\overline{X})
\stackrel{\Upsilon\tens \mathrm{fb}^{-1}}
\longrightarrow \cr
&&
\overline{F(X)}\tens \overline{F(X^\circ)} \tens \overline{F(X)} 
  \stackrel{\id\tens  \overline{h \la}\tens\id}\longrightarrow 
\overline{F(X)}\tens \overline{F(X^\circ)} \tens \overline{F(X)}
\stackrel{\mathrm{fb}\tens \Upsilon^{-1}}\longrightarrow \cr 
&&
F(\overline{X})\tens \overline{F(X) \tens F(X^{\circ})}
 \stackrel{\id\tens \overline{F(\mathrm{ev}^R_{X})\,f}}\longrightarrow 
F(\overline{X}) \tens \overline{1} \stackrel{\id\tens \star^{-1}}\longrightarrow
F(\overline{X}) \tens 1
\stackrel{l^{-1}}\longrightarrow
F(\overline{X})
\end{eqnarray*}
{} From the axiom giving the behaviour of left and right identities
under bar, this is
\begin{eqnarray*} 
&& F(\overline{X}) \stackrel{\mathrm{fb}^{-1}}\longrightarrow
\overline{F(X)} \stackrel{\overline{l}}\longrightarrow
\overline{F(X)\tens 1}
 \stackrel{\Upsilon}\longrightarrow
\overline{1} \tens \overline{F(X)}
\stackrel{\overline{f^{-1}\,F(\mathrm{coev}^R_{ X})}\tens\id}\longrightarrow 
\overline{F(X^\circ)\tens F(X)}\tens \overline{F(X)} \cr
&&
\stackrel{\Upsilon\tens \id} \longrightarrow 
\overline{F(X)}\tens \overline{F(X^\circ)} \tens \overline{F(X)} 
  \stackrel{\id\tens  \overline{h \la}\tens\id}\longrightarrow 
\overline{F(X)}\tens \overline{F(X^\circ)} \tens \overline{F(X)}
\stackrel{\id\tens \Upsilon^{-1}}\longrightarrow \cr 
&&
\overline{F(X)}\tens \overline{F(X) \tens F(X^{\circ})}
 \stackrel{\id\tens \overline{F(\mathrm{ev}^R_{X})\,f}}\longrightarrow 
\overline{F(X)} \tens \overline{1} \stackrel{\Upsilon^{-1}}\longrightarrow
\overline{1\tens F(X)}
\stackrel{\overline{r^{-1}}}\longrightarrow
\overline{F(X)} \stackrel{\mathrm{fb}}\longrightarrow
F(\overline{X})
\end{eqnarray*}
The functoriality of $\Upsilon$ gives
\begin{eqnarray*} 
&& F(\overline{X}) \stackrel{\mathrm{fb}^{-1}}\longrightarrow
\overline{F(X)} \stackrel{\overline{l}}\longrightarrow
\overline{F(X)\tens 1}
\stackrel{\id\tens f^{-1}\,F(\mathrm{coev}^R_{ X})}\longrightarrow
 \overline{F(X)\tens F(X^\circ) \tens F(X)}\cr
&&
\stackrel{(\Upsilon\tens \id)\,\Upsilon} \longrightarrow 
\overline{F(X)}\tens \overline{F(X^\circ)} \tens \overline{F(X)} 
  \stackrel{\id\tens  \overline{h \la}\tens\id}\longrightarrow 
\overline{F(X)}\tens \overline{F(X^\circ)} \tens \overline{F(X)}
\stackrel{\Upsilon^{-1}\,(\id\tens \Upsilon^{-1})}\longrightarrow \cr 
&&
 \overline{F(X)\tens F(X^\circ) \tens F(X)}
 \stackrel{\overline{F(\mathrm{ev}^R_{X})\,f\tens\id}}\longrightarrow 
\overline{1\tens F(X)}
\stackrel{\overline{r^{-1}}}\longrightarrow
\overline{F(X)} \stackrel{\mathrm{fb}}\longrightarrow
F(\overline{X})
\end{eqnarray*}
and by more functoriality and cancelling,
\begin{eqnarray*} 
&& F(\overline{X}) \stackrel{\mathrm{fb}^{-1}}\longrightarrow
\overline{F(X)} \stackrel{\overline{l}}\longrightarrow
\overline{F(X)\tens 1}
\stackrel{\id\tens f^{-1}\,F(\mathrm{coev}^R_{ X})}\longrightarrow
 \overline{F(X)\tens F(X^\circ) \tens F(X)}\cr
&&
  \stackrel{ \overline{\id\tens (h \la)\tens\id}}\longrightarrow 
 \overline{F(X)\tens F(X^\circ) \tens F(X)}
 \stackrel{\overline{F(\mathrm{ev}^R_{X})\,f\tens\id}}\longrightarrow 
\overline{1\tens F(X)}
\stackrel{\overline{r^{-1}}}\longrightarrow
\overline{F(X)} \stackrel{\mathrm{fb}}\longrightarrow
F(\overline{X})
\end{eqnarray*}
By the definition of $S^{-1}$, this is
\begin{eqnarray*} 
&& F(\overline{X}) \stackrel{\mathrm{fb}^{-1}}\longrightarrow
\overline{F(X)}
 \stackrel{\overline{S^{-1}(h)\la}}\longrightarrow
\overline{F(X)} \stackrel{\mathrm{fb}}\longrightarrow
F(\overline{X})\ .
\end{eqnarray*}
By the definition of $T$, this is the same as $p\la:F(\overline{X}) \to
F(\overline{X}) $ where $\overline{p}=T^{-1}(S^{-1}(h))$. \end{proof}

\begin{theorem}
If we define $*$ by $T\, S^{-1}$ (allowing antilinear maps by dropping the bars),
then: 

\noindent (1)\quad $(a\,b)^*=b^*\, a^*$

\noindent (2)\quad $a^{**}=a$

\noindent (3)\quad $\Delta(h^*)=\Delta(h)^{*\tens *}$

\noindent (4)\quad $S(h)^*=S^{-1}(h^*)$

\noindent I.e.\ this $*$ makes $H$ into a $*$-Hopf algebra. 
\end{theorem}
\begin{proof} By the reconstruction, $h,h'\in H$ are equal
if and only if $h\la$ and $h'\la$ are identical on all objects, so the previous propositions
give these results (again, dropping the bars)

\noindent (a) \quad $T(h\,h')=T(h)\, T(h')$  (see \ref{propiii1}).

\noindent (b) \quad $T^2=\id$  (see \ref{propiii2}).

\noindent (c) \quad $T$ reverses order in the coproduct  (see \ref{propiii3}).

\noindent (d) \quad $ST=T^{-1}S^{-1}$  (see \ref{propiii4}).

\noindent 
As $S$ reverses order in the product, (a) shows that $*$ reverses order in the product (1).

\noindent 
As $S$ reverses order in the coproduct, (c) shows that $*$ preserves order in the coproduct
(3).

\noindent 
As (b) shows that $T^{-1}=T$, from (d) we get $T\, S^{-1}\,T\, S^{-1}=\id$ (2).

\noindent 
This is sufficient to show that we have a $*$-Hopf algebra, and (4) follows automatically.
\end{proof}

We have described here reconstruction for the `basic case' of $*$-Hopf algebra from its strong bar category of modules, which confirms that our definition
of a bar category and functor is reasonable (otherwise the reconstruction would not have worked). 
Critical here was  Definition~\ref{recdef} where we defined the properties of the functor we were to use. Using a weaker notion in which $\mathrm{fb}$ is not required to be monoidal would be expected to give a strong quasi-$*$ Hopf algebra. Having this with a general bar category is expected to give a general quasi-$*$-Hopf algebra. Similarly if $F$ is multiplicative but not monoidal (see \cite{majquasitann}), then
we would expect to be able to reconstruct a $*$-quasi Hopf algebra.


\begin{thebibliography}{ggghhh}


\bibitem{albmajoct}
{Albuquerque H.\  \& Majid S.,} {\sl Quasialgebra Structure of the Octonions,  
J. Algebra 220 (1999) 188-224.}





\bibitem{shombeggs}
{ Al-Shomrani M. M. \&  Beggs  E. J.,} {\sl Making
nontrivially associated modular categories from finite groups.
International Journal of Mathematics and Mathematical Science, vol
2004, no 42, 2231-2264, 2004.}

\bibitem{shombeggs3}
{ Al-Shomrani M. M. \&  Beggs  E. J.,} {\sl Further results on the coset representative
category, to appear.}


\bibitem{barr}
{Barr, M.,} {*-Autonomous categories, Springer Lecture Notes in Mathematics 752,
Berlin, 1979 }


\bibitem{baez2hilbert}
{Baez, J.C.,} {Higher-dimensional algebra. II. 2-Hilbert spaces, Adv. Math. 127 (1997), no. 2, 125--189.}


\bibitem{cosrep}
{ Beggs  E. J.,} {\sl Making non-trivially associated tensor
categories from left coset representatives. Journal of Pure and
Applied Algebra, vol 177 , 5 - 41, $2003$.}


\bibitem{BHMnonass}
{Bouwknegt\, P., Hannabuss\, K.C., Mathai\, V.,}
{Nonassociative tori and applications to T-duality, Commun.\ Math.\ Phys.\
264 (2006) 41-69.}

\bibitem{BrzWis}
Brzezinski, T. and Wisbauer, R, {Corings and comodules} LMS lecture Notes 309, CUP (2003).


\bibitem{connesbook}
{Connes A.,} { Connes, Alain Noncommutative geometry. 
Academic Press, Inc., San Diego, CA, 1994. }

\bibitem{Dri}
Drinfeld V.G, {\sl Quantum Groups}, Proc. ICM., AMS, 1986.

\bibitem{driquasi}
{Drinfeld V.G.,} {\sl Quasi Hopf algebras. Leningrad Math.\ J.\ 1 (1990) 1419-1457.}

\bibitem{JS}
Joyal, A \& Street, R., {\sl Braided tensor categories} Adv. Maths. 102 (1993) 20-78.

\bibitem{Mac}
S. Mac Lane, Categories for the Working Mathematician.




\bibitem{MajBook}
{Majid\,S.,} {Foundations of quantum group theory, CUP}



\bibitem{majstarbr}
{Majid S.,} {\sl *-Structures on Braided Spaces, J. Math. Phys. 36 (1995)  
4436-4449.}


\bibitem{majquastpoin}
{Majid S.,} {\sl Quasi-* Structure on q-Poincar\'e Algebras, J. Geom. Phys.22  
(1997) 14-58.}


\bibitem{majquasitann}
{Majid S.,} {\sl Tannaka-Krein Theorem For QuasiHopf Algebras and Other  
Results, Contemp. Math. 134 (1992) 219-232.}


\bibitem{woron87}
{Woronowicz, S. L.,} {\sl Twisted ${\rm SU}(2)$ group. An example of a noncommutative differential calculus.
Publ. Res. Inst. Math. Sci. 23 (1987), no. 1, 117--181.}


 

\end{thebibliography}
\end{document}